\documentclass[10pt,a4paper,reqno]{amsart}

\usepackage[british]{babel}
\usepackage{csquotes}
\usepackage{epsfig}
\usepackage{amsmath}
\usepackage{bbm}
\numberwithin{equation}{section}
\usepackage[hidelinks]{hyperref}

\usepackage[left=2.5cm,right=2.5cm,top=2.5cm,bottom=2.5cm,footskip=1.5cm,headsep=1cm]{geometry}
\usepackage{appendix}
\usepackage{thmtools}
\newtheorem{claim}{Claim}
\usepackage{ascmac}

\newcommand{\interior}[1]{%
  {\kern0pt#1}^{\mathrm{o}}%
}

\usepackage[T1]{fontenc} 
\usepackage{lmodern} 
\rmfamily 
\DeclareFontShape{T1}{lmr}{b}{sc}{<->ssub*cmr/bx/sc}{}
\DeclareFontShape{T1}{lmr}{bx}{sc}{<->ssub*cmr/bx/sc}{}
\usepackage{lipsum}
\usepackage{mathtools}
\usepackage{amsmath}
\usepackage{amssymb}
\usepackage{tikz}
\usetikzlibrary{matrix,positioning,decorations.pathreplacing,calc}
\usetikzlibrary{
decorations.text,%
decorations.markings,%
shadows}
\usepackage{adjustbox}
\usepackage{graphicx} 

\usepackage{caption}
\usepackage{subcaption}

\usepackage{dsfont} 
\usepackage[format=hang]{caption}

\usepackage[normalem]{ulem}

\usepackage{bm}

\usepackage[
sorting=nty,
maxnames=99,
maxalphanames=5,
natbib=true,
sortcites]{biblatex}

\DeclareNameAlias{default}{family-given}

\graphicspath{{figures/}}

\AtEveryBibitem{
 \clearfield{url}
 \clearfield{issn}
 \clearfield{isbn}
 \clearfield{urldate}
 
 \ifentrytype{book}{
  \clearfield{pages}}{%
 }
}

\usepackage{xargs}                      
\usepackage[prependcaption,textsize=tiny,textwidth=3cm]{todonotes}
\newcommandx{\unsure}[2][1=]{\todo[linecolor=red,backgroundcolor=red!25,bordercolor=red,#1]{#2}}
\newcommandx{\change}[2][1=]{\todo[linecolor=blue,backgroundcolor=blue!25,bordercolor=blue,#1]{#2}}
\newcommandx{\info}[2][1=]{\todo[linecolor=OliveGreen,backgroundcolor=OliveGreen!25,bordercolor=OliveGreen,#1]{#2}}
\newcommandx{\improvement}[2][1=]{\todo[linecolor=black,backgroundcolor=black!25,bordercolor=black,#1]{#2}}
\newcommandx{\thiswillnotshow}[2][1=]{\todo[disable,#1]{#2}}
\usepackage{amsthm}
\usepackage[noabbrev,capitalize]{cleveref}
\crefname{equation}{}{}

\newtheorem{theorem}{Theorem}[section]
\newtheorem{lemma}[theorem]{Lemma}
\newtheorem{proposition}[theorem]{Proposition}
\newtheorem{corollary}[theorem]{Corollary}

\theoremstyle{definition}
\newtheorem{definition}[theorem]{Definition}

\newtheorem{remark}[theorem]{Remark}

\crefname{assumption}{Assumption}{Assumptions}
\crefname{definition}{Definition}{Definitions}
\crefname{corollary}{Corollary}{Corollaries}
\crefname{enumi}{item}{items}

\usepackage{fancyhdr}

\fancypagestyle{plain}{%
\fancyhead[C]{}
\fancyfoot[C]{}

}
\fancypagestyle{nosection}{%
\fancyhead[CE]{}
\fancyhead[CO]{}
\fancyfoot[CE,CO]{\thepage}
}
\newcommand{\nm}{\noalign{\smallskip}}
\newcommand{\ds}{\displaystyle}

\DeclareMathOperator{\R}{\mathbb{R}}

\renewcommand{\tilde}{\widetilde}

\renewcommand{\qedsymbol}{$\blacksquare$}

\usepackage{fancyhdr}

\fancypagestyle{plain}{%
\fancyhead[C]{}
\fancyfoot[C]{}

}
\fancypagestyle{nosection}{%
\fancyhead[CE]{}
\fancyhead[CO]{}
\fancyfoot[CE,CO]{\thepage}
}

\renewcommand{\epsilon}{\varepsilon}

\let\emptyset\varnothing

\renewcommand{\tilde}{\widetilde}

\usepackage{tcolorbox}
\usetikzlibrary{tikzmark,calc}

\usepackage{enumerate}
\usepackage{upgreek}

\begin{document}

\title[]{Sharp Quantitative Stability of Extremizable Brascamp-Lieb Inequalities}

\author[B. Ammari]{Bader Ammari}
\address[]{Department of Mathematics, ETH Z\"{u}rich, R\"{a}mistrasse 101, CH-8092 Z\"{u}rich, Switzerland}
\email{bader.ammari@math.ethz.ch}

\begin{abstract}
    We prove a sharp quantitative version of the general Euclidean forward Brascamp-Lieb inequalities that correspond to extremizable data, without any restriction on the exponents. More specifically, for any extremizable datum, we show that the Brascamp-Lieb deficit controls the squared Hellinger distance to the family of normalized extremizing tuples, including those with non-Gaussian factors. 
\end{abstract}

\maketitle

\vspace{0.2cm}

\vspace{0.2cm}
\vspace{0.5cm}

\section{Introduction}

Recently, there has been a growing interest in understanding quantitative stability for functional and geometric inequalities. For instance, this has led to several works investigating the stability of Sobolev and Sobolev-type inequalities \cite{Bianchi_Egnell,Figalli_Zhang,Dolbeault_Esteban_Figalli_Frank_Loss}, isoperimetric inequalities \cite{Fusco_Maggi_Pratelli,Figalli_Maggi_Pratelli,Figalli_Zhang2}, the Brunn-Minkowski inequality \cite{Figalli_Jerison,VanHintum_Spink_Tiba,Figalli_Tiba_VH,Figalli_Tiba_VH2} or the variance Brascamp-Lieb inequality \cite{Machado_Ramos,Karoly_Yaozhong_Cyril}. 

\hfill \\
Following this line of research, we shall investigate the stability of extremizers of the Euclidean forward Brascamp-Lieb inequalities.

\subsection{Brascamp-Lieb inequalities}
\label{Sec1.1}
The forward Brascamp-Lieb inequalities are of the form $$\int_{\mathbb{R}^n}\prod_{j=1}^m f_j(C_jx)^{p_j}\, dx \le \mathrm{BL}(\textbf{C},p)\prod_{j=1}^m \bigg(\int f_j\, dx\bigg)^{p_j},$$ where $C_j:\mathbb{R}^n \rightarrow \mathbb{R}^{n_j}$ are linear surjective maps, $0<p_j \le 1$ and $f_j$ are non-negative real-valued functions on $\mathbb{R}^{n_j},$ or equivalently, $$\int_{\mathbb{R}^n}\prod_{j=1}^m f_j(C_jx)\, dx \le \mathrm{BL}(\textbf{C},p)\prod_{j=1}^m \|f_j\|_{L^{1/p_j}(dx_j)}.$$ The crucial observation of Lieb \cite{Lieb} is that the optimal constant $$0<\mathrm{BL}(\mathbf{C},p) \le +\infty,$$ usually referred to as the Brascamp-Lieb constant, is saturated by Gaussians. Moreover, the necessary and sufficient conditions for the Brascamp-Lieb constant to be finite or for its value to be attained were established in \cite{BCCT}: indeed, $\mathrm{BL}(\textbf{C},p)$ is finite if and only if $$n=\sum_{j=1}^m p_j n_j \text{ and } \dim(V) \le \sum_{j=1}^m p_j\dim(C_jV) \text{ for any } V \text{ linear subspace of } \mathbb{R}^n.$$ In addition, Valdimarsson in \cite{Vald} settled the question of the form of extremizers of forward Brascamp-Lieb data that are extremizable, \textit{i.e.} whose Brascamp-Lieb constants are attained.

\hfill \\
Let us recall the notion of equivalence between Brascamp-Lieb data from \cite[Definition 3.1]{BCCT}. 

\begin{definition}
    Two Brascamp-Lieb data $(\textbf{C},p)$ and $(\textbf{C'},p')$ are equivalent if there exist invertible linear transformations $D:H' \rightarrow H$ and $D_j:H_j' \rightarrow H_j$ such that $C_j'=D_j^{-1}C_jD$ for all $1 \le j \le m$ and $p=p'.$ We refer to $D$ and $D_j$ as the intertwining transformations.
\end{definition}

\hfill \\ 
The Brascamp-Lieb constant of two equivalent Brascamp-Lieb data are closely related thanks to \cite[Lemma 3.8]{BCCT}. 

\begin{lemma} \label{Equivalence_BL_Constants_Lemma}
    Suppose that $(\textbf{C},p)$ and $(\textbf{C'},p')$ are two equivalent Brascamp-Lieb data, with intertwining maps $D:H' \rightarrow H$ and $D_j:H_j' \rightarrow H_j.$ Then, $$\mathrm{BL}(\textbf{C'},p')=\frac{\prod_{j=1}^m |\det D_j|^{p_j}}{|\det D|}\mathrm{BL}(\textbf{C},p).$$ 
\end{lemma}

\hfill \\
Finally, we remind the reader of the definition of extremizable and Gaussian-extremizable Brascamp-Lieb data. 

\begin{definition}
    Let $(\textbf{C},p)$ be a Brascamp-Lieb datum. 
    $(\textbf{C},p)$ is extremizable if there exist functions $(g_j)_{1 \le j \le m}$, with $g_j: H_j \rightarrow [0;+\infty)$ and $0<\int_{H_j}g_j<+\infty$, so that $$\int_H \prod_{j=1}^m g_j(C_jx)^{p_j}\, dx =\mathrm{BL}(\textbf{C},p) \prod_{j=1}^m \Big(\int_{H_j}g_j\Big)^{p_j}.$$ Moreover, $(\textbf{C},p)$ is Gaussian-extremizable if there exist $A_j \in \mathrm{Sym}_{>0}(H_j), 1 \le j \le m$ such that $$\mathrm{BL}(\textbf{C},p)=\bigg(\frac{\prod_{j=1}^m (\det A_j)^{p_j}}{\det \big(\sum_{j=1}^m p_jC_j^*A_jC_j\big)} \bigg)^{1/2}.$$
\end{definition}
\hfill \\
In fact, by \cite[Proposition 6.5]{BCCT}, extremizable Brascamp-Lieb data are Gaussian-extremizable and therefore, by \cite[Proposition 3.6]{BCCT}, equivalent to a geometric Brascamp-Lieb datum, \textit{i.e.} $$\sum_{j=1}^m p_j C_j^*C_j=\mathrm{id}_H, \quad C_iC_i^*=\mathrm{id}_{H_i} \quad \forall 1 \le i \le m,$$ which has Brascamp-Lieb constant $1.$ Furthermore, the independent subspaces of a geometric Brascamp-Lieb datum $(\textbf{C},p)$ are defined as the non-zero subspaces of $H$ of the form $$\bigcap_{i \in A}H_i \cap \bigcap_{i \in A^c}H_i^\perp, \quad A \subset \{1,\cdots,m\},$$ identifying $H_i$ as subspaces of $H$ thanks to the fact that $C_i^*C_i$ is a projection for any $1 \le i \le m$. We may always assume such an identification when working with geometric Brascamp-Lieb data.

\hfill \\
More recently, the author revisited the optimal transport proof of the forward Brascamp-Lieb inequality (see \cite{Barthe}) to show a more general regularized version in \cite{Regularized_BL}. More precisely, if $f_j$ are more log-convex than $g_{Q_j}$ for any $1 \le j \le m$, $$\int_{\mathbb{R}^n}\prod_{j=1}^m f_j(C_jx)^{p_j}\, dx \le \mathrm{BL}(\textbf{C},p,\textbf{Q})\prod_{j=1}^m \bigg(\int f_j \, dx\bigg)^{p_j},$$ where the optimal constant $$0<\mathrm{BL}(\textbf{C},p,\textbf{Q}) \le \mathrm{BL}(\textbf{C},p),$$ is again saturated by Gaussians that satisfy such a log-convexity property. In addition, a structural study was conducted in order to determine the necessary and sufficient conditions for finiteness, and to find the conditions for a given datum to be extremizable in the general case of mixed conditions of log-convexity and log-concavity as well as the extremizers.

\subsection{The stability question} Unlike the well-known special reverse Brascamp-Lieb inequalities, that are Pr\'ekopa-Leindler inequality at the functional level and Brunn-Minkowski inequality at the geometric level (see \cite{Karoly_Alessio_Joao,Figalli_Jerison,Figalli_Tiba_VH,Figalli_Tiba_VH2,Figalli_Tiba_VH3}), the stability of forward Brascamp-Lieb inequalities has not been, until recently, studied in a systematic manner. One of the main difficulties stems through the richness of extremizers that were fully characterized by \cite{Vald} and that heavily depend on the structure of the considered Brascamp-Lieb datum. 

\hfill \\
As the forward Brascamp-Lieb inequality unifies central inequalities, as H\"older's, Young or Loomis-Whitney inequality, it is not surprising that the first stability result of a special case of a forward Brascamp-Lieb inequality was shown for the H\"older's inequality in \cite{Aldaz}. Another special case of a forward Brascamp-Lieb inequality, whose stability was studied, is Young's inequality $$\int_{\mathbb{R}^n \times \mathbb{R}^n}f_1(x)f_2(x-y)f_3(y)\, dx \, dy\le A_q^n\prod_{i=1}^3 \|f_i\|_{L^{q_i}},$$ where $$\sum_{i=1}^3 \frac1{q_i}=2, \quad A_q=\prod_{i=1}^3 \Big(\frac{q_i^{1/q_i}}{(q_i')^{1/q_i'}}\Big)^{1/2},$$   essentially through the works of Christ (see \cite{Christ_Young}, \cite{Christ_Young2}). More precisely, in \cite{Christ_Young}, qualitative precompactness is shown without any restriction on the exponents, while \cite{Christ_Young2} proves quantitative stability for exponents $1<q_j \le 2$ \footnote{The constant may be chosen uniformly, if the exponents $(q_j)_{1 \le j \le 3}$ belong to a compact subset of $(1;2)^3$.} for $1 \le j \le 3$, always under the homogeneity assumption $\sum_{i=1}^3 q_i^{-1}=2.$ Let us further note that the stability of the usual dual version of Young's inequality $$\|f_1 *f_2\|_{L^{q_3'}} \le A_q^n \|f_1\|_{L^{q_1}}\|f_2\|_{L^{q_2}}$$ is shown, without needing $q_3$ to be inside the interval $(1;2).$  Finally, a discrete stability version of the Loomis-Whitney inequality had been explored in \cite{Ellis_Friedgut_Kindler_Yehudayoff} via information-theoretic tools. For the continuous case, stability properties of near-extremizers for the Loomis-Whitney inequality were very recently investigated in \cite{Bennett_Oh}. 

\hfill \\
Note that, in this subsection, the $L^p$ norms are taken with respect to the Lebesgue measure. 

\hfill \\
In \cite{Bennett_Christ}, the authors prove the following quantitative stability statement with respect to (uncentered) Gaussians: if $\mathrm{BL}(\textbf{C},\frac1q)$ is finite, $1<q_j<2$ for all $j$, then there are constants $c_j=c(q_j,n_j)>0$ such that \begin{equation} \label{Partial_Gaussian_Stability}\frac{\ds \int_{\mathbb{R}^n}\prod_{j=1}^m f_j(C_jx)\, dx}{\ds \prod_{j=1}^m \|f_j\|_{L^{q_j}}}  \le \mathrm{BL}\Big(\textbf{C},\frac1q \Big)\prod_{j=1}^m \Biggl[ 1-c_j \Bigl( \frac{\mathrm{dist}_{q_j}(f_j,\mathbbm{G}(n_j))}{\|f_j\|_{L^{q_j}}} \Bigl)^2 \Biggl].\end{equation} Here, $$\mathrm{dist}_{q_j}(f_j,\mathbbm{G}(n_j))=\inf_{g \text{ uncentered Gaussian}}\|f_j-g\|_{L^{q_j}}.$$

\hfill \\
The proof is based on the Fourier-invariance of the Brascamp-Lieb constants (see \cite[Proposition 2.3]{Bennett_Bez_NonlinearBL}) and a quantitative stability result developed in \cite{Christ_HausdorffYoung}. It is crucial here that the exponents $q_j$ all belong to $(1;2)$ to be able to employ the stability of Hausdorff-Young and as the above result imposes all the extremizers of the Brascamp-Lieb datum $(\textbf{C},\frac1q)$ to be Gaussian, which is indeed the case under such an assumption on the exponents. 

\hfill \\
Moreover, by \cite[Proposition 3.4, 3.6]{Bennett_Christ}, Bennett and Christ show that, for Gaussian extremizable data, under the only assumption that $q_j>1$ for all $j$, the sharp decay exponent at any term $$\frac{\mathrm{dist}_{q_j}(f_j,\mathbbm{G}(n_j))}{\|f_j\|_{L^{q_j}}},$$ is $\max\{2,q_j\},$ for any $1 \le j \le m.$ In particular, \eqref{Partial_Gaussian_Stability} fails for any sequence of positive constants $c_j$ if $q_k>2$ for some $k$, by using the phenomenon of escaping perturbation at infinity. This may strongly remind the reader of the situation encountered for the stability of Sobolev inequality in \cite{Figalli_Zhang}.  

\hfill \\
Finally, they prove, that under the assumption of simplicity of $(\textbf{C},p)$ and $1<q_j<2$ for any $1 \le j \le m$, a qualitative stability statement with respect to the extremizers, which consist, up to scaling and translation, of a unique centered Gaussian extremizer. In fact, \eqref{Partial_Gaussian_Stability} allows to reduce this statement to a qualitative stability statement among Gaussians, which is then proven by compactness.

\hfill \\
The goal of this manuscript is to completely answer the question of stability for extremizable Brascamp-Lieb datum. We believe this also unlocks the door to exciting further understanding of the stability for Brascamp-Lieb inequalities (or possibly, more ambitiously, reverse Brascamp-Lieb inequalities). For example, we are not yet able to answer positively the forward version\footnote{Note that in the paper, only the distance between $f$ and the standard Gaussian is considered but instead, it should be $\inf_{\sigma>0}\|f-\gamma_\sigma\|_{L^1}.$} of Conjecture 2.2 of \cite{Karoly_Fodor_Hug}, as a uniform control of the constants would be needed. 

\hfill \\ 
Our main result is as follows. 

\begin{theorem} \label{Complete_Gaussian_stability}
Let $(\textbf{C},p)$ be an extremizable Brascamp-Lieb datum. 
   There exists a constant $c=c(\textbf{C},p)>0$ such that, for any $(f_i)_{1 \le i \le m}$ integrable with $f_i \ge 0$ and $0<\int f_i<+\infty$, there exists an extremizer $(g_i)_{1 \le i \le m}$ of the Brascamp-Lieb datum $(\textbf{C},p)$ such that $\|g_i\|_{L^1}=1$ and 
    $$1-\frac{\ds \int_{\mathbb{R}^n} \prod_{i=1}^m f_i(C_ix)^{p_i}\, dx}{\ds \mathrm{BL}(\textbf{C},p)\prod_{i=1}^m \Big(\int_{\mathbb{R}^{n_i}}f_i \, dx\Big)^{p_i}} \ge c\sum_{i=1}^m p_i \Bigg\|\sqrt{\frac{f_i}{\|f_i\|_{L^1}}}-\sqrt{g_i}\, \Bigg\|_{L^2}^2.$$ 
\end{theorem}
\hfill \\
This complete result, with the sharp power, allows to give a direct answer to the question of Gaussian stability raised by Bennett and Tao in \cite{Bennett_Tao}, \textit{among extremizable data.} In our subsequent work with Figalli \cite{Bader_Alessio2}, we give a complete characterization of Gaussian-stable Brascamp-Lieb data.

\begin{corollary}
Let $(\textbf{C},p)$ be an extremizable Brascamp-Lieb datum. Then, it is Gaussian-stable if and only if all its non-zero extremizers are (uncentered) Gaussians or equivalently, a geometric Brascamp-Lieb datum that is equivalent to $(\textbf{C},p)$ has no independent subspaces. 
\end{corollary}

\begin{proof}
    By \cite[Def 7.1]{Bennett_Tao}, recall that $(\textbf{C},p)$ is Gaussian-stable if $\mathrm{BL}(\textbf{C},p)<+\infty$ and whenever $f_i:\mathbb{R}^{n_i} \rightarrow \mathbb{R}_+$ are integrable, such that $$\int_{\mathbb{R}^n}\prod_{i=1}^m f_i(C_ix)^{p_i}\, dx \ge (1-\epsilon)\mathrm{BL}(\textbf{C},p)\prod_{i=1}^m \Big(\int_{\mathbb{R}^{n_i}}f_i\Big)^{p_i},$$ for some $\epsilon>0$, then there exists (uncentered) Gaussians $ \tilde g_i:\mathbb{R}^{n_i} \rightarrow \mathbb{R}_+$ such that $$\|f_i-\tilde g_i\|_{L^1} \le o(\|f_i\|_{L^1}) \quad \text{ for any } 1 \le i \le m,$$ where $o(X)$ denotes a quantity bounded in magnitude by $c(\epsilon)X$ for some function $c(\epsilon)$, depending on the Brascamp-Lieb datum, that tends to zero as $\epsilon \rightarrow 0.$  
    
    \hfill \\
    On the one hand, it is clear that if $(\textbf{C},p)$ is extremizable and Gaussian-stable, all extremizers must be (uncentered) Gaussians. This implies that a geometric Brascamp-Lieb datum, that is equivalent to $(\textbf{C},p)$, has no independent subspaces by the characterization of extremizers of \cite{Vald}. 
    On the other hand, it is also worth recalling that, if $\tilde f_i=\frac{f_i}{\|f_i\|_{L^1}}$, $$\|\tilde f_i-g_i\|_{L^1}=\|(\sqrt{\tilde f_i}+\sqrt{g_i})(\sqrt{\tilde f_i}-\sqrt{g_i})\|_{L^1} \le \|\sqrt{\tilde f_i}+\sqrt{g_i}\|_{L^2}\|\sqrt{\tilde f_i}-\sqrt{g_i}\|_{L^2},$$ for any $1 \le i \le m,$ and $$\|\sqrt{\tilde f_i}+\sqrt{g_i}\|_{L^2}^2=2+2\int \sqrt{\tilde f_i}\sqrt{g_i} \le 4,$$ so that $$\|\tilde f_i-g_i\|_{L^1} \le 2\|\sqrt{\tilde f_i}-\sqrt{g_i}\|_{L^2}.$$ Thus, \cref{Complete_Gaussian_stability} implies \begin{equation} \label{BL_L1_Stability} 1-\frac{\int_{\mathbb{R}^n} \prod_{i=1}^m f_i(C_ix)^{p_i}\, dx}{\mathrm{BL}(\textbf{C},p)\prod_{i=1}^m \Big(\int_{\mathbb{R}^{n_i}}f_i \, dx\Big)^{p_i}} \ge \frac c4\sum_{i=1}^m p_i \Big(\frac{\|f_i-\tilde g_i\|_{L^1}}{\|f_i\|_{L^1}}\Big)^2,\end{equation} for $(\tilde g_i)_{1 \le i \le m}=(\|f_i\|_{L^1}g_i)_{1 \le i \le m}$ an extremizer of the Brascamp-Lieb datum $(\textbf{C},p)$, that is by assumption Gaussian.  
    Hence, if $$\int_{\mathbb{R}^n}\prod_{i=1}^m f_i(C_ix)^{p_i}\, dx \ge (1-\epsilon)\mathrm{BL}(\textbf{C},p)\prod_{i=1}^m \Big(\int_{\mathbb{R}^{n_i}}f_i\Big)^{p_i},$$ then $$\epsilon \ge \frac c4 p_i \Big(\frac{\|f_i-\tilde g_i\|_{L^1}}{\|f_i\|_{L^1}}\Big)^2 \Longrightarrow  \|f_i-\tilde g_i\|_{L^1} \le \sqrt{\frac4{p_ic}\epsilon} \cdot  \|f_i\|_{L^1},$$ for any $1 \le i \le m$, as desired.
\end{proof}

\hfill \\
Moreover, we retrieve a more general form of \eqref{Partial_Gaussian_Stability}, with no assumptions on the exponents. 

\begin{theorem}
    Let $(\textbf{C},\frac1q)$ be an extremizable Brascamp-Lieb datum. There exists a constant $\tilde c=\tilde c(\textbf{C},\frac1q)>0$ such that, for any $(f_i)_{1 \le i \le m}$ with $f_i \ge 0$ and $0<\|f_i\|_{L^{q_i}(\mathbb{R}^{n_i})}<+\infty$, there exists a tuple $(g_i)_{1 \le i \le m}$ such that $\int_{\mathbb{R}^n}\prod_{i=1}^m g_i(C_ix)=\mathrm{BL}(\textbf{C},\frac1q)\prod_{i=1}^m \|g_i\|_{L^{q_i}}$ and $$\mathrm{BL}\Big(\textbf{C},\frac1q\Big)\prod_{j=1}^m \Bigg[ 1-\tilde c \Bigg(\frac{\|f_j-g_j\|_{L^{q_j}}}{\|f_j\|_{L^{q_j}}}\Bigg)^{\max\{2,q_j\}}\Bigg] \ge \frac{\ds \int_{\mathbb{R}^n}\prod_{i=1}^m f_i(C_ix)\, dx}{\ds \prod_{i=1}^m \|f_i\|_{L^{q_i}}}.$$
\end{theorem}

\begin{proof}
    Let $$F=f^{1/p}, G=g^{1/p}, \|F\|_p=\|G\|_p=1,$$ with $p \ge 1.$ Then, we claim that $$\|\sqrt f-\sqrt g\|_{L^2}^2 \ge \min\big\{1,\frac {p^2}4\big\}\|F-G\|_{L^p}^{\max\{2,p\}}.$$ 
Indeed, if $p \ge 2$, $$|a^r-b^r| \ge |a-b|^r \text{ for } r \ge 1, a,b \ge 0,$$ then $$|F-G|^{p/2} \le |\sqrt f-\sqrt{g}| \Longrightarrow \|F-G\|_{L^p}^p \le \|\sqrt{f}-\sqrt{g}\|_{L^2}^2.$$  
If $1 \le p<2$, setting $r=\frac2p$, one has $F=\sqrt{f}^r, G=\sqrt{g}^r$ and $\|\sqrt{f}\|_{L^2}=\|\sqrt{g}\|_{L^2}=1.$ In addition, $$a^r-b^r=r(a-b)\int_0^1 \big( (1-t)b+ta\big)^{r-1}\, dt,$$ so that by H\"older's inequality, $$\begin{array}{lll}
\ds \|F-G\|_{L^p}=\|\sqrt{f}^r-\sqrt{g}^r\|_{L^p} &\le& \ds \frac2p \|\sqrt{f}-\sqrt{g}\|_{L^2} \Big\|\int_0^1 \big((1-t)\sqrt{g}+t\sqrt{f}\big)^{r-1}\, dt \Big\|_{L^{2p/(2-p)}} \\ \\
&\le& \ds \frac2p  \|\sqrt{f}-\sqrt{g}\|_{L^2}\int_0^1 \|(1-t)\sqrt{g}+t\sqrt{f}\|_{L^2}^{r-1}\, dt \\ \\
&\le& \ds \frac2p \|\sqrt{f}-\sqrt{g}\|_{L^2}.
\end{array}$$ 
\hfill \\
Applying \cref{Complete_Gaussian_stability} to the tuple $(f_i^{q_i})_{1 \le i \le m}$, with $p_i=\frac1{q_i}$ for any $1 \le i \le m$, we get that $$\begin{array}{lll}
\ds 1-\frac{\int_{\mathbb{R}^n}\prod_{i=1}^m f_i(C_ix)\, dx}{\mathrm{BL}\Big(\textbf{C},\frac1q \Big)\prod_{i=1}^m \|f_i\|_{L^{q_i}}} &\ge& \ds c\sum_{i=1}^m \frac1{q_i} \Big\| \sqrt{\frac{f_i^{q_i}}{\|f_i^{q_i}\|_{L^1}}}-\sqrt{g_i}\Big\|_{L^2}^2 \\ \\
&\ge& \ds c\sum_{i=1}^m \frac1{q_i}\min\big\{1,\frac {q_i^2}4\big\}\Big\|\frac{f_i}{\|f_i\|_{L^{q_i}}}-g_i^{1/q_i}\Big\|_{L^{q_i}}^{\max\{2,q_i\}} \\ \\
&\ge& \ds c\sum_{i=1}^m \frac1{q_i} \min\big\{1,\frac{q_i^2}4\big\} \Big(\frac{\|f_i-\tilde g_i\|_{L^{q_i}}}{\|f_i\|_{L^{q_i}}}\Big)^{\max\{2,q_i\}},
\end{array}$$ where $(\tilde g_i)_{1 \le i \le m}=(\|f_i\|_{L^{q_i}}g_i^{1/q_i})_{1 \le i \le m}$ satisfies that $$\int_{\mathbb{R}^n}\prod_{i=1}^m \tilde g_i(C_ix)\, dx =\mathrm{BL}\Big(\textbf{C},\frac1q\Big)\prod_{i=1}^m \|\tilde g_i\|_{L^{q_i}}.$$
Defining $$\tilde c=\min_{1 \le i \le m}\Big\{ \frac{c \min\{1,q_i^2/4\}}{q_i},2^{-\max\{2,q_i\}}\Big\},$$ we obtain that $$1-\frac{\int_{\mathbb{R}^n}\prod_{i=1}^m f_i(C_ix)\, dx}{\mathrm{BL}\Big(\textbf{C},\frac1q \Big)\prod_{i=1}^m \|f_i\|_{L^{q_i}}} \ge \tilde c \sum_{i=1}^m \Big(\frac{\|f_i-\tilde g_i\|_{L^{q_i}}}{\|f_i\|_{L^{q_i}}}\Big)^{\max\{2,q_i\}},$$ and $$0 \le \tilde c \Big(\frac{\|f_i-\tilde g_i\|_{L^{q_i}}}{\|f_i\|_{L^{q_i}}}\Big)^{\max\{2,q_i\}} \le 1 \text{ for any } 1 \le i \le m.$$
\hfill \\
To conclude, it suffices to use the standard inequality $$\prod_{j=1}^m (1-x_j) \ge 1-\sum_{j=1}^m x_j,$$ for any $x_1,\cdots,x_m \in [0;1].$ 
\end{proof}

\subsection{Structure of the paper}

In \cref{Sec:2}, we first prove the quantitative stability of the (regularized) Brascamp-Lieb inequalities in the restricted function space of (uncentered) Gaussians.
Its purpose is to highlight an elegant way of proving quantitative stability in this setting. Moreover, this section also serves to build intuition for establishing the spectral gap in \cref{Sec:3}, especially when we deal with Hermite polynomials.
\cref{Sec:3} is devoted to the study of the local stability of Brascamp-Lieb inequalities\footnote{In the general case, we will not pursue the stability of regularized Brascamp-Lieb inequalities to simplify exposition. However, the arguments used here should extend in a straightforward manner to obtain similar results in the regularized context.} for general functions, by carefully investigating the relevant quadratic form and establishing the spectral gap required for the proof of \cref{Complete_Gaussian_stability}, which is then postponed to \cref{Sec:4}. The sharpness is also addressed in this last section. 

\hfill \\
\textbf{Notation.} In our estimates, we often write positive constants as $C$ and $c$, with $C$ usually denoting a constant larger than $1$ and $c$ a constant less than $1.$ The constants may vary between appearances, even within a chain of inequalities. The notation $a \sim b$ indicates that both inequalities $a \le Cb$ and $b \le Ca$ hold.

\section{Stability in the Gaussian setting}
\label{Sec:2}

\subsection{Non-regularized case}

The main focus of this section is to establish a quantitative version for Gaussians of Brascamp-Lieb inequalities with respect to the Hilbert-Schmidt norm, which is contained in \eqref{BL_L1_Stability} by the next lemma. All constants, that will be used, may depend on the datum $(\textbf{C},p)$, or in the regularized case, $(\textbf{C},p,\textbf{Q})$, but not on the matrices considered.

\begin{lemma} \label{L1Gaussian_Behaviour}
For any $A,B \in \mathrm{Sym}_{>0}(n)$, $$\|g_A-g_B\|_{L^1} \sim \min \Big\{ 1, \|A^{-1/2}BA^{-1/2}-I_n\|_{\mathrm{HS}} \Big\},$$ where we will write $$g_A(x)=(2\pi)^{-n/2}\det(A)^{1/2}e^{-\frac12 \langle Ax,x\rangle}.$$
\end{lemma}

\begin{proof}
Although this lemma seems to be quite folklore, we provide a short self-contained argument to show that $$\frac1{3\pi}\min\Big\{ 1,\|A^{-1/2}BA^{-1/2}-I_n\|_{\mathrm{HS}}\Big\} \le \|g_A-g_B\|_{L^1} \le 2\min\Big\{1, \|A^{-1/2}BA^{-1/2}-I_n\|_{\mathrm{HS}}\Big\}.$$ It suffices to prove the previous inequality with $\|g_{I_n}-g_M\|_{L^1}$ for $M=\mathrm{diag}(\lambda_1,\cdots,\lambda_n)$, with $\lambda_i>0$ representing the eigenvalues of $A^{-1/2}BA^{-1/2}.$ Set $$\delta:=\|M-I_n\|_{\mathrm{HS}}=\Big(\sum_{i=1}^n (\lambda_i-1)^2 \Big)^{1/2}.$$ Then, $$\frac{\delta}{2+\delta} \le \Big(\sum_{i=1}^n \frac{(\lambda_i-1)^2}{(\lambda_i+1)^2}\Big)^{1/2} \le \delta,$$ so that $$\frac13 \min\{1,\delta\} \le \min \Big\{ 1, \Big(\sum_{i=1}^n \frac{(\lambda_i-1)^2}{(\lambda_i+1)^2}\Big)^{1/2} \Big\} \le \min\{1,\delta\}.$$ Thus, $$\alpha:=\int_{\mathbb{R}^n}\sqrt{g_M(x)g_{I_n}(x)}\, dx =\prod_{i=1}^n \Big(\frac{4\lambda_i}{(\lambda_i+1)^2}\Big)^{1/4}=\prod_{i=1}^n \Big(1-\frac{(\lambda_i-1)^2}{(\lambda_i+1)^2}\Big)^{1/4},$$ so that by H\"older's inequality, $$\|g_{I_n}-g_M\|_{L^1} \le \Big( \int (\sqrt{g_{I_n}}-\sqrt{g_M})^2 \Big)^{1/2} \Big( \int (\sqrt{g_{I_n}}+\sqrt{g_M})^2 \Big)^{1/2}=2\sqrt{1-\alpha^2}.$$ Now, for the upper bound, $$1-\alpha^2 \le \sum_{i=1}^n \bigg(1-\sqrt{1-\frac{(\lambda_i-1)^2}{(\lambda_i+1)^2}} \bigg) \le \sum_{i=1}^n \frac{(\lambda_i-1)^2}{(\lambda_i+1)^2}.$$ For the lower bound, if we consider $R$ a Gaussian random variable with density $\sqrt{g_M g_{I_n}}/\alpha$ and $Z=\log \frac{g_M}{g_{I_n}}$, one has $$\|g_{I_n}-g_M\|_{L^1}=2\alpha \mathbb{E}_R\Big| \sinh\Big(\frac Z2\Big)\Big| \ge \alpha \mathbb{E}_R|Z|.$$ Under $R$, $$Z=\frac12 \sum_{i=1}^n \log \lambda_i-\sum_{i=1}^n \frac{\lambda_i-1}{\lambda_i+1}G_i^2,$$ where $G_i \sim \mathcal N(0,1)$ are independent and we may consider an independent copy $Z'$ of $Z$ with $$Z'=\frac12 \sum_{i=1}^n \log \lambda_i-\sum_{i=1}^n \frac{\lambda_i-1}{\lambda_i+1}(G_i')^2 \text{ where } G_i' \sim \mathcal N(0,1) \text{ are independent}.$$ If we define $U_i=\frac{G_i+G_i'}{\sqrt 2}, V_i=\frac{G_i-G_i'}{\sqrt 2}$, all the $U_i$ and $V_i$ are independent standard normal random variables, such that $$Z-Z' =^d  -2\sum_{i=1}^n \frac{\lambda_i-1}{\lambda_i+1}U_iV_i.$$ Consequently, by independence, $$\Big(Z-Z'|(U_1,\cdots,U_n)=(u_1,\cdots,u_n)\Big) \sim \mathcal N \Big(0,4\sum_{i=1}^n \frac{(\lambda_i-1)^2}{(\lambda_i+1)^2}u_i^2\Big),$$ hence $$\mathbb{E}_R(|Z-Z'|\; | \; U)=2\sqrt{\frac2\pi}\Big(\sum_{i=1}^n \frac{(\lambda_i-1)^2}{(\lambda_i+1)^2}U_i^2\Big)^{1/2},$$ as $\mathbb{E}|\mathcal N(0,\sigma^2)|=\sigma \sqrt{2/\pi}.$ Taking the expectation, $$2\mathbb{E}_R(|Z|) \ge \mathbb{E}_R(|Z-Z'|)=2\sqrt{\frac2\pi}\mathbb{E} \Big[  \Big(\sum_{i=1}^n \frac{(\lambda_i-1)^2}{(\lambda_i+1)^2}U_i^2\Big)^{1/2}\Big] \ge \frac4\pi \Big(\sum_{i=1}^n \frac{(\lambda_i-1)^2}{(\lambda_i+1)^2}\Big)^{1/2},$$ as, by Cauchy-Schwarz, setting $w_i=\frac{(\lambda_i-1)^2}{(\lambda_i+1)^2 \sum_{j=1}^n \frac{(\lambda_j-1)^2}{(\lambda_j+1)^2}}$ for any $1 \le i \le n$, $$\Big(\sum_{i=1}^n w_i U_i^2\Big)^{1/2} \ge \sum_{i=1}^n w_i|U_i|.$$ In conclusion, $$\|g_{I_n}-g_M\|_{L^1} \ge \frac{2\alpha}{\pi} \Big(\sum_{i=1}^n \frac{(\lambda_i-1)^2}{(\lambda_i+1)^2}\Big)^{1/2},$$ while clearly, $$\|g_{I_n}-g_M\|_{L^1} \ge \int (\sqrt{g_M}-\sqrt{g_{Id_n}})^2=2(1-\alpha),$$ so that $$\|g_{I_n}-g_M\|_{L^1} \ge \frac1\pi\min \Big\{1, \Big(\sum_{i=1}^n \frac{(\lambda_i-1)^2}{(\lambda_i+1)^2}\Big)^{1/2}\Big\}.$$
\end{proof}

\begin{remark}
    From the proof of \cref{L1Gaussian_Behaviour}, we may conclude that $$\|g_A-g_B\|_{L^1} \sim \|\sqrt{g_A}-\sqrt{g_B}\|_{L^2} \text{ for any } A,B \in \mathrm{Sym}_{>0}(n).$$ As previously, it suffices to prove this equivalence between $\|g_{I_n}-g_M\|_{L^1}$ and $\|\sqrt{g_{I_n}}-\sqrt{g_M}\|_{L^2}.$ On the one hand, by H\"older's inequality, $$\|g_{I_n}-g_M\|_{L^1} \le 2\|\sqrt{g_{I_n}}-\sqrt{g_M}\|_{L^2}.$$ On the other hand, using the previous notations, $$\|\sqrt{g_{I_n}}-\sqrt{g_M}\|_{L^2}^2=2(1-\alpha) \le 2 \sum_{i=1}^n \frac{(\lambda_i-1)^2}{(\lambda_i+1)^2} \Longrightarrow \|g_{I_n}-g_M\|_{L^1} \ge \frac{\sqrt 2 \alpha}{\pi}\|\sqrt{g_{I_n}}-\sqrt{g_M}\|_{L^2}.$$ If $\alpha \ge \frac12$, we conclude that $$\|g_{I_n}-g_M\|_{L^1} \ge \frac{1}{\sqrt 2\pi}\|\sqrt{g_{I_n}}-\sqrt{g_M}\|_{L^2}.$$ Otherwise, $$\|\sqrt{g_{I_n}}-\sqrt{g_M}\|_{L^2}=\sqrt{2(1-\alpha)} \ge 1,$$ and thus, $$\|g_{I_n}-g_M\|_{L^1} \ge \|\sqrt{g_{I_n}}-\sqrt{g_M}\|_{L^2}^2 \ge \|\sqrt{g_{I_n}}-\sqrt{g_M}\|_{L^2},$$ using that $$|x-y| \ge |\sqrt x-\sqrt y|^2 \quad  \forall x,y>0.$$ 
    \hfill \\
    In particular, in order to prove the Gaussian analogue of \cref{Complete_Gaussian_stability}, we will work with the more familiar $L^1$ distance. 
\end{remark}

\hfill \\
The main result of this subsection is as follows.

\begin{theorem} \label{Gaussian_Stab_Thm}
Assume that $(\textbf{C},p)$ is a Gaussian extremizable Brascamp-Lieb datum. 
There exists $\epsilon_0>0$ depending only on the datum $(\textbf{C},p)$, such that for any $\epsilon \in [0;\epsilon_0]$, if
$$\frac{\prod_{j=1}^m \det(B_j)^{p_j/2}}{\det(\sum_{j=1}^m p_jC_j^*B_jC_j)^{1/2}}\ge (1-\epsilon)\mathrm{BL}(\textbf{C},p),$$ then there exists $(\tilde B_j)_{1 \le j \le m}$ such that $$\|g_{B_j}-g_{\tilde B_j}\|_{L^1} \le C\sqrt{\epsilon}.$$ In particular, $$\|{B_j}-{\tilde B_j}\|_{\mathrm{HS}} \le C\sqrt{\epsilon}\|B_j\|_{\mathrm{HS}},$$ where $$\mathrm{BL}(\textbf{C},p)=\bigg(\frac{\prod_{j=1}^m \det(\tilde B_j)^{p_j}}{\det(\sum_{j=1}^m p_jC_j^*\tilde B_j C_j)}\bigg)^{1/2}.$$ 
\end{theorem}

\begin{remark}
    The inequality $$\begin{array}{lll}
    \ds \|B_j-\tilde B_j\|_{\mathrm{HS}} &=& \|B_j^{1/2}(I_{n_j}-B_j^{-1/2}\tilde B_j B_j^{-1/2})B_j^{1/2}\|_{\mathrm{HS}} \\ \\
    &\le& \ds  \sqrt{n_j}\|B_j\|_{\mathrm{HS}}\|I_{n_j}-B_j^{-1/2}\tilde B_j B_j^{-1/2}\|_{\mathrm{HS}}, \end{array}$$ and \cref{L1Gaussian_Behaviour} explain the second part of the theorem. 
\end{remark}

\hfill \\
By change of variables, we may restrict ourselves to studying an arbitrary geometric Brascamp-Lieb datum $(\textbf{C},p)$, \textit{i.e.} $$C_jC_j^*=\mathrm{id}_{H_j} \quad \forall 1 \le j \le m \text{ and } \sum_{j=1}^m p_jC_j^*C_j=\mathrm{id}_H,$$ which has Brascamp-Lieb constant equal to $1.$ 
Note that we will freely use the Hilbert-space formalism, by writing $H$ and $H_j$ as general Euclidean spaces but of course, $H \cong \mathbb{R}^n, H_j \cong \mathbb{R}^{n_j}$ for any $1 \le j \le m.$
Moreover, define $$F(A_1,\cdots,A_m)=-\sum_{j=1}^m p_j\log \det(A_j)+\log \det \Bigl(\sum_{j=1}^m p_jC_j^*A_jC_j \Bigl) \ge 0,$$ as $(\textbf{C},p)$ is a geometric datum. In addition, $$F(\mathrm{id}_{H_1}+\epsilon A_1,\cdots,\mathrm{id}_{H_m}+\epsilon A_m)=\frac{\epsilon^2}2 \Biggl[\sum_{j=1}^m p_j\mathrm{Tr}(A_j^2)-\mathrm{Tr}\biggl(\Big(\sum_{j=1}^m p_jC_j^*A_jC_j\Big)^2\biggl) \Biggl]+o(\epsilon^2),$$ since $$F(\mathrm{id}_{H_1},\cdots,\mathrm{id}_{H_m})=0, \quad \nabla F(\mathrm{id}_{H_1},\cdots,\mathrm{id}_{H_m})=0.$$

\hfill \\
Let us denote $$Q(A)=\sum_{j=1}^m p_j\mathrm{Tr}(A_j^2)-\mathrm{Tr}\biggl(\Big(\sum_{j=1}^m p_jC_j^*A_jC_j\Big)^2\biggl),$$ for any $A=(A_j)_{1 \le j \le m} \subset \prod_{j=1}^m \mathrm{Sym}(n_j)$, as the quadratic form on matrices corresponding to the Brascamp-Lieb datum $(\textbf{C},p).$

\hfill \\
We first prove the precompactness for Gaussian extremizable datum.

\begin{theorem}[Precompactness] \label{Precompactness_Full}
    If $(\textbf{C},p)$ is a Gaussian-extremizable datum, then for every $\epsilon>0$, there exists $\delta>0$ such that for any $(A_j)_{1 \le j \le m}$ so that $\mathrm{BL}(\textbf{C},p;A) \ge (1-\delta)\mathrm{BL}(\textbf{C},p)$, there exists $(B_j)_{1 \le j \le m}$ an extremizer of $(\textbf{C},p)$ so that $$\|A_j-B_j\|_{\mathrm{HS}} \le \epsilon \|A_j\|_{\mathrm{HS}}.$$  
\end{theorem}

Sequentially, this means that if $(A_j^{(n)})_{1 \le j \le m}$ is such that $F(A^{(n)}) \rightarrow_{n \rightarrow +\infty} 0$, then there exists a sequence of extremizers $(B_j^{(n)})_{1 \le j \le m}$ so that $$\frac{\|A_j^{(n)}-B_j^{(n)}\|_{\mathrm{HS}}}{\|A_j^{(n)}\|_{\mathrm{HS}}} \rightarrow_{n \rightarrow +\infty} 0 \quad \forall 1 \le j \le m.$$

\begin{remark}
The main claim used in the proof of \cite[Lemma 4.1]{Bennett_Christ}, which proves \cref{Precompactness_Full} in the case of simple\footnote{$(\textbf{C},p)$ is simple if there does not exist non-trivial (\textit{i.e.} neither $\{0\}$ nor $H$) critical subspaces.} data, is the compactness of the set of near-extremizers matrices, \textit{i.e.} suppose that $(\textbf{C},p)$ is a simple datum, and that $p_j \in (0;1)$ for any $1 \le j \le m.$ Then, for every $\eta \in (0;1)$, the set $$\mathcal T_\eta=\left \{ A=(A_j)_{1 \le j \le m}: \det \Big(\sum_{j=1}^ m p_jC_j^*A_jC_j\Big)=1, \prod_{j=1}^m (\det A_j)^{p_j} \ge (1-\eta)\mathrm{BL}(\textbf{C},p)^2 \right\},$$ is compact.
Its generalization does not hold in the non-simple case. Indeed, consider the following geometric Brascamp-Lieb datum: $n=2,m=4, p_1=p_2=p_3=p_4=\frac12$ and $$C_1=C_2=P_1,C_3=C_4=P_2 \text{ and } A_1^{(n)}=A_2^{(n)}=n,A_3^{(n)}=A_4^{(n)}=\frac1n.$$ Then, $$\det \Big(\sum_{j=1}^4 p_jC_j^*A_j^{(n)}C_j \Big)=\det \Big(\begin{pmatrix}
        n & 0 \\
        0 & \frac1n 
    \end{pmatrix}\Big)=1,\; \prod_{j=1}^4 (\det A_j^{(n)})^{p_j}=1.$$ 
\end{remark}

\hfill \\
Thus, we need to argue differently. 
For a positive definite matrix $B$, there exists a unique real symmetric matrix $\tilde B:=\log(B)$ such that $$B=e^{\tilde B}.$$ In particular, this implies that $$\det(B)=e^{\mathrm{Tr}(\log B)}.$$

\hfill \\
First, a useful observation is as follows: set $K=\bigoplus_{j=1}^m H_j,$ and $$\mathcal B=\big\{X=X_1 \oplus \cdots \oplus X_m: X_j \in \mathrm{Sym}(n_j) \big\}, \quad V:H \rightarrow K, x \mapsto (p_1^{1/2}C_1x,\cdots,p_m^{1/2}C_mx).$$
Let us define a weighted inner product on $\mathcal B$ by $$\langle X,Y\rangle_p=\sum_{j=1}^m p_j\mathrm{Tr}(X_jY_j).$$

\hfill \\
Then, as $(\textbf{C},p)$ is a geometric datum, $P=VV^*$ is an orthogonal projection from $K$ onto $S=\mathrm{Ran}(V).$ Moreover, denoting $A=(A_j)_{1 \le j \le m}$, it is worth noticing that $$V^*AV=\sum_{j=1}^m p_jC_j^*A_jC_j,$$ and $$\mathrm{Tr}(V^* \log(A)V) =\sum_{j=1}^m p_j\mathrm{Tr}(C_j^*\log(A_j)C_j)=\sum_{j=1}^m p_j\mathrm{Tr}(\log(A_j))=\sum_{j=1}^m p_j \log \det(A_j).$$
Therefore, \begin{equation} \label{Logarithmic_Form_BLDeficit_Gaussian}
F(A_1,\cdots,A_m)=\mathrm{Tr}\big[\log(V^*AV)\big]-\mathrm{Tr}(V^*\log(A)V)=\mathrm{Tr}_S\big[\log(PAP)\big|_S\big]-\mathrm{Tr}(P\log(A)P|_S).\end{equation}

\hfill \\
\begin{lemma}
    Let $K$ be a finite-dimensional Hilbert space, $S \subset K$ and $P_S$ the orthogonal projection onto $S.$ If $A=(A_j)_{1 \le j \le m}>0$, then $$\log(P_SAP_S)|_S \ge P_S\log(A)P_S|_S.$$ Moreover, equality holds if and only if $$AS \subset S.$$
\end{lemma}

\begin{proof}
    Write $$A=\begin{pmatrix}
        B & C \\
        C^* & D
    \end{pmatrix}, B=P_SAP_S|_S,$$ with respect to the decomposition $K=S \oplus S^\perp.$ For $x>0$, $$\log(x)=\int_0^{+\infty}\Big( \frac1{1+t}-\frac1{x+t}\Big)\, dt.$$ Therefore, $$\log B=\int_0^{+\infty} \Big(\frac1{1+t}I_S-(B+tI_S)^{-1}\Big)\, dt,$$ and similarly, $$P_S(\log A)P_S|_S=\int_0^{+\infty}\Big( \frac1{1+t}I_S-P_S(A+tI_K)^{-1}P_S|_S \Big)\, dt,$$ using that $P_SI_KP_S|_S=I_S.$ Substracting both terms yields $$\log B-P_S(\log A)P_S|_S=\int_0^{+\infty} \Big(P_S(A+tI_K)^{-1}P_S|_S-(B+tI_S)^{-1}\Big)\, dt.$$ Moreover, $$A+tI_K=\begin{pmatrix}
        B+tI_S & C \\
        C^* & D+tI_{S^\perp}
    \end{pmatrix},$$ and by Schur's complement, $$P_S(A+tI_K)^{-1}P_S|_S=(B+tI_S-C(D+tI_{S^\perp})^{-1}C^*)^{-1} \ge (B+tI_S)^{-1},$$ for any $t \ge 0.$ This proves the first claim. If equality holds, then equality holds in $$P_S(A+tI_K)^{-1}P_S|_S \ge (B+tI_S)^{-1},$$ so that $$C(D+tI_{S^\perp})^{-1}C^*=0 \Longrightarrow C=0.$$ In other words, $A$ is block-diagonal with respect to the decomposition $S \oplus S^\perp.$ Conversely, if $C=0$, $\log A$ is also block-diagonal with respect to $S\oplus S^\perp$ and equality clearly holds.
\end{proof}
\hfill \\
Define $$\mathfrak{n}=\big\{ A \in \mathcal B :AS \subset S  \big\},$$ and $$\mathcal E=\big\{ E=E_1 \oplus \cdots \oplus E_m: E_i \in \mathrm{Sym}_{>0}(n_i), ES \subset S\big \}.$$  

\begin{remark}
    In particular, this shows another equivalent way of characterizing the Gaussian extremizers of a Brascamp-Lieb datum. Indeed, by \cref{Logarithmic_Form_BLDeficit_Gaussian} and the previous lemma, $E$ is a Gaussian extremizer if and only if $E(S) \subset S.$ In this case, $$EV=VM, \text{ where } M=V^*EV=\sum_{j=1}^m p_jC_j^*E_jC_j,$$ and component by component, this means that $$E_jC_j=C_j M \Longrightarrow E_jC_jM^{-1}C_j^*=\mathrm{id}_{H_j} \Longrightarrow C_jM^{-1}C_j^*=E_j^{-1} \quad \forall 1 \le j \le m,$$ which is the usual condition for $(E_j)_{1 \le j \le m}$ to be a Gaussian extremizer (see \cite[Proposition 3.6]{BCCT}). 
\end{remark}

\hfill \\
We claim that, if $E \in \mathcal E$, then for any $A=A_1 \oplus \cdots \oplus A_m>0$, $$F(E^{-1/2}AE^{-1/2})=F(A).$$ Indeed, $$EV=V(V^*EV) \text{ and thus, } E^{-1/2}V=V(V^*EV)^{-1/2},$$ which implies in particular that $$V^*E^{-1/2}AE^{-1/2}V=(V^*EV)^{-1/2}(V^*AV)(V^*EV)^{-1/2},$$ and $$\log \det(V^*E^{-1/2}AE^{-1/2}V)=\log \det(V^*AV)-\log \det(V^*EV).$$ 
Moreover, $$E^{-1/2}AE^{-1/2}=(E_1^{-1/2}A_1E_1^{-1/2})\oplus \cdots \oplus (E_m^{-1/2}A_mE_m^{-1/2}),$$ which yields $$\sum_{j=1}^m p_j \log \det(E_j^{-1/2}A_jE_j^{-1/2})=\sum_{j=1}^m p_j \log \det(A_j)-\sum_{j=1}^m p_j \log \det(E_j).$$ 
This shows that $$F(E^{-1/2}AE^{-1/2})=F(A)-F(E)=F(A),$$ since $F(E)=0$ by the previous lemma. 

\hfill \\
For positive definite operators $R,T$ on a finite-dimensional Hilbert space, we may define $$d(R,T)=\bigl\|\log(R^{-1/2}TR^{-1/2}) \bigl\|_{\mathrm{HS}}.$$ Note that for every invertible $G$, \begin{equation} \label{Invariance_Thompson_Distance} d(GRG^*,GTG^*)=d(R,T),\end{equation} since $(GRG^*)^{-1}GTG^*=G^{-*}(R^{-1}T)G^*$ is similar to $R^{-1}T$, itself similar to $R^{-1/2}TR^{-1/2}$; in particular, they have the same eigenvalues.  

\hfill \\
For $A,E \in \mathcal B$, define the weighted distance $$d_p(A,E)^2=\sum_{j=1}^m p_jd(A_j,E_j)^2.$$ 

\begin{lemma} \label{Orthogonal_Minimizer}
    For every $A=A_1 \oplus \cdots \oplus A_m>0$, there exists $E \in \mathcal E$ attaining $$\inf_{\tilde E \in \mathcal E}d_p(A,\tilde E).$$ If $E$ is such a minimizer and $$e^X=E^{-1/2}AE^{-1/2},$$ then $X \in \mathfrak{n}^{\perp_p}.$
\end{lemma}
\hfill \\
Here, $$\mathfrak{n}^{\perp_p}=\big \{ X \in \mathcal B : \langle X,Z\rangle_p=0 \quad \forall Z \in \mathfrak{n} \big\}.$$

\begin{proof}
    As $I_K \in \mathcal E$, it suffices to minimize over the set $$\big \{ \tilde E \in \mathcal E: d_p(A,\tilde E) \le 1+d_p(A,I_K)\big\}.$$ Moreover, if $d_p(A,\tilde E) \le R$, the eigenvalues of $\tilde E_j^{-1/2}A_j \tilde E_j^{-1/2}$ are in the interval $[e^{-Rp_j^{-1/2}};e^{Rp_j^{-1/2}}]$ for each $1 \le j \le m$; especially, $E$ lies in the compact subset of a positive cone and the condition $E(S) \subset S$ is stable under closedness. As the above set is compact and the distance is continuous, we indeed find a minimizer $E \in \mathcal E.$ 
    
\hfill \\
    Write $e^X=E^{-1/2}AE^{-1/2}.$ The goal is to prove that $\langle X,Z\rangle_p=0$ for any $Z \in \mathfrak{n}.$ First, $$E_s=E^{1/2}e^{sZ}E^{1/2} \in \mathcal{E},$$ as the block diagonal structure and the fact that $S$ is stable under those matrices is preserved. Then, by minimality of $E$, $$h(s)=d_p(E_s,A)^2 \text{ attains its minimum at } 0.$$ By \cref{Invariance_Thompson_Distance}, $$h(s)=d_p(e^{sZ},e^X)^2=\sum_{j=1}^m p_j \mathrm{Tr}\Big[ \log(e^{-sZ_j/2}e^{X_j}e^{-sZ_j/2})^2 \Big]$$ and we compute $h'(0)=0.$ As $$\frac d{ds}(e^{-sZ_j/2}e^{X_j}e^{-sZ_j/2})\big|_{s=0}=-\frac {Z_j}2 e^{X_j}-e^{X_j} \frac {Z_j}2,$$ and denoting $\Psi(Y)=\mathrm{Tr}[(\log Y)^2]$, $$D\Psi_{e^X}[H]=2\mathrm{Tr}(Xe^{-X}H),$$ we get that $$0=h'(0)=\sum_{j=1}^m p_j \Big( -\mathrm{Tr}(X_je^{-X_j}Z_je^{X_j})-\mathrm{Tr}(X_je^{-X_j}e^{X_j}Z_j)\Bigl)=-2\langle X,Z\rangle_p,$$ where we used that $X_j$ commutes with $e^{X_j}$ and the cyclicity of trace.  
\end{proof}

\begin{lemma}
    Let $(X_n)_{n \ge 1} \in \mathfrak{n}^{\perp_p}.$ If $$F(e^{X_n}) \rightarrow_{n \rightarrow +\infty} 0,$$ then $$X_n \rightarrow_{n \rightarrow +\infty} 0 \text{ in Hilbert-Schmidt norm.}$$
\end{lemma}

\begin{proof}
    For $X$ in block form, define $$\varphi_X(t)=F(e^{tX}), t \in \mathbb{R}.$$ It is quite standard that $\varphi_X$ is convex (see for example \cite[Section 4]{Barthe_Wolff} or \cite[Proposition 6]{Barthe}). Essentially, $$\varphi_X(t)=\log \det(V^*e^{tX}V)-\Big(\sum_{j=1}^m p_j\mathrm{Tr}(X_j)\Big)t,$$ the second term is then linear whereas the first is convex in $t$ by Cauchy-Binet formula.
    \hfill \\
    Set $$f(A_1,\cdots,A_m)=F(e^{A_1},\cdots,e^{A_m})=\log \det(V^*e^AV)-\sum_{j=1}^m p_j \mathrm{Tr}(A_j).$$ For fixed $A \in \mathcal B$, denoting $$M(t)=V^*e^{tA}V,$$ one has $$M(0)=\mathrm{Id}_H, M'(0)=V^*AV, M''(0)= V^*A^2V.$$ Since $$\frac d{dt}\log \det M(t)\Big|_{t=0}=\mathrm{Tr}(M^{-1}(0)M'(0))=\sum_{j=1}^m p_j \mathrm{Tr}(A_j),$$ the linear term cancels and gives $Df(0)=0.$ Moreover, the Hessian quadratic form is $$Q(A)=\mathrm{Tr}(V^*A^2V)-\mathrm{Tr}((V^*AV)^2),$$ which we already computed for local perturbations of the identity matrix. Using the fact that $P=VV^*$ is a projection, $$Q(A)=\mathrm{Tr}(PA^2)-\mathrm{Tr}(PAPA)=\mathrm{Tr}(PA(\mathrm{Id}-P)A)=\|PA(\mathrm{Id}-P)\|_{\mathrm{HS}}^2 \ge 0.$$ In particular, if $Q(A)=0$, $$PA(\mathrm{Id}-P)=0 \Longrightarrow A(S) \subset S.$$ In other words, $$\ker(Q)=\mathfrak{n}.$$ Consequently, there exists $\lambda_p>0$ such that $$Q(A) \ge \lambda_p \|A\|_p^2 \quad \forall A \in \mathfrak{n}^{\perp_p}.$$ 
As $$F(e^X)=\frac12 Q(X)+O(\|X\|_p^3) \text{ as } X \rightarrow 0,$$ by choosing $r>0$ small enough, we obtain that \begin{equation} \label{Local_Exp_SP_Gaussian} F(e^X) \ge \frac{\lambda_p}4 \|X\|_p^2 \quad \forall X \in \mathfrak{n}^{\perp_p} \text{ with } \|X\|_p \le r,\end{equation} by continuity of $F.$
    
    \hfill \\
    Since $\mathbb{R} \ni t \mapsto F(e^{tX})$ is convex and by \eqref{Local_Exp_SP_Gaussian}, we conclude that \begin{equation} \label{Global_Exp_SP_Gaussian} F(e^X) \ge c_p\min \Big \{ \|X\|_p, \|X\|_p^2 \Big\} \quad \forall X \in \mathfrak{n}^{\perp_p}.\end{equation} This proves the desired claim.  
\end{proof}

\begin{proof}[Proof of \cref{Precompactness_Full}]
    Let $A^{(n)}=A_1^{(n)} \oplus \cdots \oplus A_m^{(n)}$ such that $F(A_1^{(n)},\cdots,A_m^{(n)}) \rightarrow_{n \rightarrow +\infty} 0.$ 
    By \cref{Orthogonal_Minimizer}, choose $E^{(n)}$ minimizing $d_p(A^{(n)},E)$ over $E \in \mathcal E.$ As shown, we may write $$({E^{(n)}})^{-1/2}A^{(n)}(E^{(n)})^{-1/2}=e^{X_n},$$ with $X_n \in \mathfrak{n}^{\perp_p}.$ Then, $$F(e^{X_n})=F(({E^{(n)}})^{-1/2}A_n({E^{(n)}})^{-1/2})=F(A_n) \rightarrow_{n \rightarrow +\infty} 0.$$ By the previous lemma, $X_n \rightarrow_{n \rightarrow +\infty} 0,$ so that $$({E^{(n)}})^{-1/2}A^{(n)}({E^{(n)}})^{-1/2} \rightarrow_{n \rightarrow +\infty} I_K \text{ or in other words, } (E_j^{(n)})^{-1/2}A^{(n)}_j (E_j^{(n)})^{-1/2} \rightarrow_{n \rightarrow +\infty} \mathrm{id}_{H_j},$$ in Hilbert-Schmidt norm. This proves that $$\|g_{A_j^{(n)}}-g_{E_j^{(n)}}\|_{L^1} \rightarrow_{n \rightarrow +\infty} 0 \quad \forall 1 \le j \le m,$$ as wished. 
\end{proof}
\hfill \\
Based on the approach adopted to prove precompactness, it is easy to conclude the desired Gaussian stability statement.

\begin{proof}[Proof of \cref{Gaussian_Stab_Thm}]
By \cref{Orthogonal_Minimizer}, choosing $E=(E_j)_{1 \le j \le m}$ that minimizes $d_p(A,E)$, if $$E_j^{-1/2}A_jE_j^{-1/2}=e^{X_j}, X=X_1 \oplus \cdots \oplus X_m,$$ then $X \in \mathfrak{n}^{\perp_p}$ and $d_p(A,E)^2=\|X\|_p^2.$  Moreover, $$F(e^X)=F(A) \le -2\log(1-\epsilon) \sim \epsilon.$$ Since $X \in \mathfrak{n}^{\perp_p}$, \cref{Global_Exp_SP_Gaussian} gives $\|X\|_p \le r$ for sufficiently small $\epsilon>0$. In particular, this implies that $$C_p \epsilon \ge \|X\|_p^2=\sum_{j=1}^m p_j \big\|\log(E_j^{-1/2}A_jE_j^{-1/2}) \big\|_{\mathrm{HS}}^2,$$ using here again \eqref{Local_Exp_SP_Gaussian}. For $\|Y_j\|_{\mathrm{HS}}$ small, $\|e^{Y_j}-\mathrm{id}_{H_j}\|_{\mathrm{HS}} \sim \|Y_j\|_{\mathrm{HS}},$ and thus, considering $Y_j=\log(E_j^{-1/2}A_jE_j^{-1/2})$, $$\|g_{A_j}-g_{E_j}\|_{L^1} \le C_p \sqrt{\epsilon} \text{ for any } 1 \le j \le m,$$ as wished. 
\end{proof}

\begin{remark} \label{QF_Degree2_Rk}
We provide an alternative way of giving a complete description of the kernel of the quadratic form, which will be particularly useful in the next section, as it retrieves the form of Gaussian extremizers as stated in \cite{Vald}. 
Let us define the map $$T(A)=\sum_{j=1}^m p_jC_j^*A_jC_j \text{ where } A=(A_j)_{1 \le j \le m}.$$ Then, one may write $$Q(A)=\|A\|_p^2-\|T(A)\|_{\mathrm{HS}}^2.$$ Notice that $$\langle T(A),Y\rangle_{\mathrm{HS}}=\sum_{j=1}^m p_j\mathrm{Tr}(A_jC_jYC_j^*) \Longrightarrow T^*(Y)=(C_jYC_j^*)_{1 \le j \le m},$$ where the adjoint is taken with respect to the inner product $\langle \cdot,\cdot\rangle_p$ on $\mathcal B.$ Moreover, $$\|T^*(Y)\|_p^2=\sum_{j=1}^m p_j \mathrm{Tr}((C_jYC_j^*)^2)=\sum_{j=1}^m p_j\mathrm{Tr}(C_j^*C_jYC_j^*C_jY).$$ As $C_j^*C_j$ is a projection, $$\begin{array}{lll}
    \ds -\mathrm{Tr}(C_j^*C_jYC_j^*C_jY)+\mathrm{Tr}(C_j^*C_jY^2) &=& \ds \mathrm{Tr}(C_j^*C_jY(\mathrm{id}_H-C_j^*C_j)Y) \\ \\
    &=& \ds \|C_j^*C_jY(\mathrm{id}_H-C_j^*C_j)\|_{\mathrm{HS}}^2 \ge 0.
    \end{array}$$
In particular, $$\|T^*(Y)\|_p^2\le \sum_{j=1}^m p_j \mathrm{Tr}(C_j^*C_jY^2)=\mathrm{Tr}(Y^2)=\|Y\|_{\mathrm{HS}}^2.$$
This means that $$\|T(A)\|_{\mathrm{HS}}=\sup_{\|Y\|_{\mathrm{HS}}=1}\langle T(A),Y\rangle_{\mathrm{HS}}=\sup_{\|Y\|_{\mathrm{HS}}=1}\langle A,T^*(Y)\rangle_p \le \|A\|_p.$$
Now, if there is equality, $$T^*T(A)=A.$$ Indeed, $$\|A\|_p^2=\|T(A)\|_{\mathrm{HS}}^2=\langle A,T^*T(A)\rangle_{p} \le \|A\|_p\|T^*T(A)\|_p \le \|A\|_p\|T(A)\|_{\mathrm{HS}}=\|A\|_p^2,$$ and so, there is equality in Cauchy-Schwarz, which proves the desired claim.  
\hfill \\
If we define $Y=T(A)$, we thus obtain that $$\mathrm{Tr}(C_j^*C_jYC_j^*C_jY)=\mathrm{Tr}(C_j^*C_jY^2) \Longrightarrow C_j^*C_jY(\mathrm{id}_H-C_j^*C_j)=0$$ for any $1 \le j \le m.$
\hfill \\
Consequently, $$Q(A) \ge 0 \text{ and } Q(A)=0 \iff A_j=C_jYC_j^* \text{ for any } 1 \le j \le m,$$ where $$C_j^*C_jY(\mathrm{id}_H-C_j^*C_j)=0 \quad \forall 1 \le j \le m.$$ 
\hfill \\
By taking adjoints, $C_j^*C_jY=C_j^*C_jYC_j^*C_j=YC_j^*C_j$ and in particular, $C_j^*C_j(\mathrm{id}_H-C_j^*C_j)Y=C_j^*C_jY(\mathrm{id}_H-C_j^*C_j)=0.$ This implies that $$Y=
    \left(\begin{array}{ c | c }
    Y|_{H_j} & 0 \\
    \hline
    0 &  Y|_{H_j^\perp}
  \end{array} \right) \text{ for any } 1 \le j \le m.
$$ Since this holds true for any $1 \le j \le m$, $$Y=\left( \begin{array}{ c | c } 
    \nm
    Y|_{K_1}   & 0   \\ 
    \hline
    0 &  \left. \begin{array}{c | c} \nm \ddots & 0 \\
    \hline
    0 & \left. \begin{array}{c | c} \nm  Y|_{K_{k_0}} & 0 \\
    \nm
    \hline
     0 & Y|_{K_{\mathrm{dep}}}  \end{array} \right.  \end{array} \right. 
    \end{array} \right),$$ where $\oplus_{k=1}^{k_0}K_k \oplus K_{\mathrm{dep}}$ is the independent decomposition of $H.$ 
    
\hfill \\
Let $d_{k_0+1},\cdots,d_{k_1}$ be the eigenvalues of $Y$ on $K_{\mathrm{dep}}.$ Then, the eigenspace $K_k:=K_\mathrm{dep} \cap \ker(Y|_{K_{\mathrm{dep}}}-d_k\mathrm{id}_{K_{\mathrm{dep}}})$ is critical for any $k_0+1 \le k \le k_1$. Indeed, $P_{K_k}$ and $P_{H_j}$ commute and in particular, $$\mathrm{Tr}(P_{H_j}P_{K_k})=\dim(H_j \cap K_k) \text{ for any } 1 \le j \le m.$$ This implies that $$
\sum_{j=1}^m p_j\dim (H_j \cap K_k)=\sum_{j=1}^m p_j \mathrm{Tr}(C_j^*C_jP_{K_k})=\mathrm{Tr}(P_{K_k})=\dim(K_k).$$ 
Hence, we get a critical decomposition of $K_{\mathrm{dep}}$, \textit{i.e.} $$K_{\mathrm{dep}}=\bigoplus_{k=k_0+1}^{k_1}K_k.$$

\hfill \\
In conclusion, there exists matrices $(D_k)_{1 \le k \le k_0}$ and real numbers $(d_k)_{k_0+1 \le k \le k_1}$, so that $$Y=\sum_{k=1}^{k_0}P_{K_k}^* D_kP_{K_k}+\sum_{k=k_0+1}^{k_1}d_kP_{K_k}^*P_{K_k},$$ and $$A_j=C_jYC_j^*=\sum_{k=1}^{k_0}C_jP_{K_k}^*D_kP_{K_k}C_j^*+\sum_{k=k_0+1}^{k_1}d_kC_jP_{K_k}^*P_{K_k}C_j^*,$$ which is exactly the form of the Gaussian optimizers of \cite{Vald}. 
\end{remark}

\subsection{Regularized case}

We now would like to show how \cref{Gaussian_Stab_Thm} can be extended to a regularized version. More precisely, we actually may expect that one has a slight improvement compared to the non-regularized case, since the first-order term in the Taylor expansion generally does not need to vanish. 

\hfill \\
By change of variables, we may restrict to studying the stability of a generalized geometric Brascamp-Lieb datum around identity. By \cite[Definition 3.20]{Regularized_BL} , the geometricity conditions are in the regularized case as follows: $$C_jC_j^* \le \mathrm{id}_{H_j}, \sum_{j=1}^m p_jC_j^*C_j=\mathrm{id}_H,$$ and $$Q_j \ge \mathrm{id}_{H_j}, (Q_j-\mathrm{id}_{H_j})(\mathrm{id}_{H_j}-C_jC_j^*)=0.$$  
Define $$F(A_1,\cdots,A_m)=-\sum_{j=1}^m p_j\log \det(A_j)+\log \det \Bigl(\sum_{j=1}^m p_jC_j^*A_jC_j \Bigl) \ge 0,$$ for any $0<A_j \le Q_j$, since $(\textbf{C},p,\textbf{Q})$ is a generalized geometric datum (see \cite[Theorem 3.21]{Regularized_BL}). Then, for $0<\mathrm{id}_{H_j}+\epsilon A_j \le Q_j$, $$\begin{array}{lll}
\ds F(\mathrm{id}_{H_1}+\epsilon A_1,\cdots,\mathrm{id}_{H_m}+\epsilon A_m)=\epsilon \sum_{j=1}^m p_j \mathrm{Tr}\Big( (C_jC_j^*-\mathrm{id}_{H_j})A_j \Big) \\ 
\ds +\frac{\epsilon^2}2 \bigg[ \sum_{j=1}^m p_j \mathrm{Tr}(A_j^2)-\mathrm{Tr}\Big((\sum_{j=1}^m p_jC_j^*A_jC_j)^2\Big)  \bigg]+o(\epsilon^2).
\end{array}$$

\hfill \\
As we will often use the subspaces $\ker(C_jC_j^*-\mathrm{id}_{H_j})$ and $\ker(Q_j-\mathrm{id}_{H_j})$, set $$N_j:=\ker(C_jC_j^*-\mathrm{id}_{H_j}), \quad  V_j:=\ker(Q_j-\mathrm{id}_{H_j}) \quad \forall 1 \le j \le m.$$ 

\hfill \\
Set $(v_i^{(j)})_{1 \le i \le m_j}$ eigenvectors of $Q_j-\mathrm{id}_{H_j}$, that are associated to its non-zero eigenvalues and span $V_j^\perp$. Then, $$(C_jC_j^*-\mathrm{id}_{H_j})(Q_j-\mathrm{id}_{H_j})v_i^{(j)}=0 \Longrightarrow (C_jC_j^*-\mathrm{id}_{H_j})v_i^{(j)}=0 \quad \forall 1 \le i \le m_j.$$ In particular, $V_j^\perp \subset N_j$ or equivalently, $N_j^\perp \subset V_j$, and $$C_jC_j^*-\mathrm{id}_{H_j}=\left(\begin{array}{ c | c }
    (C_jC_j^*-\mathrm{id}_{H_j})|_{N_j^\perp} & 0 \\
    \hline
    0 & 0
  \end{array}\right),$$ where the second block is in $N_j.$ Thus, the first order term may actually be reduced to $$\sum_{j=1}^m p_j \mathrm{Tr}\big((C_jC_j^*-\mathrm{id})|_{N_j^\perp} A_j|_{N_j^\perp} \big),$$ and this allows to expect that one can have a first order control of $\big \|A_j|_{N_j^\perp} \big \|_{\mathrm{HS}},$ which is the ``perpendicular'' part of $A_j$ and the ``parallel'' part of $A_j$, which is $A_j|_{N_j}$, that should have order $2.$

\hfill \\
As in the non-regularized case, denoting $$V: H \rightarrow K, x \mapsto (p_1^{1/2}C_1x,\cdots,p_m^{1/2}C_mx),$$ we have that $$\begin{array}{lll}
F(A_1,\cdots,A_m) &=& 
\ds \mathrm{Tr}\Bigl[ \log(V^*AV)\Big]-\sum_{j=1}^m p_j \mathrm{Tr}[\log(A_j)] \\
&=& \ds \mathrm{Tr}\Bigl[ \log(V^*AV)\Big]-\sum_{j=1}^m p_j \mathrm{Tr}(C_j^*\log(A_j)C_j)+\sum_{j=1}^m p_j \mathrm{Tr}[(C_jC_j^*-\mathrm{id}_{H_j})\log(A_j)] \\
&=& \ds \mathrm{Tr}_S \Big[ \log(PAP)\Big|_S \Big]-\mathrm{Tr}\Bigl( P\log(A)P\Big|_S \Big)+\sum_{j=1}^m p_j \mathrm{Tr}[(C_jC_j^*-\mathrm{id}_{H_j})\log(A_j)].
\end{array}$$
\hfill \\
Moreover, we know from the previous section that $$\mathrm{Tr}_S \Big[ \log(PAP)\Big|_S \Big]-\mathrm{Tr}\Bigl( P\log(A)P\Big|_S \Big) \ge 0,$$ and $$(\mathrm{id}_{H_j}-C_jC_j^*)(Q_j-\mathrm{id}_{H_j})=0 \Longrightarrow (\mathrm{id}_{H_j}-C_jC_j^*)\log(Q_j)=0 \quad \forall 1 \le j \le m,$$ so that $$\mathrm{Tr}[(C_jC_j^*-\mathrm{id}_{H_j})\log(A_j)]=\mathrm{Tr}\Big[(\mathrm{id}_{H_j}-C_jC_j^*)\big( \log(Q_j)-\log(A_j) \big) \Big] \ge 0 \quad \forall 1 \le j \le m.$$ 

\hfill \\
Now, let us derive the Gaussian equality cases. 
We define as previously, $$\mathfrak{n}=\{ A \in \mathcal B: AS \subset S \},$$ and $$\mathfrak{t}=\{ A \in \mathcal B: A_j=P_{N_j}A_jP_{N_j} \}.$$

\begin{lemma}
    Let $0<A_j \le Q_j.$ Then, $$F(A)=0 \iff A-\mathrm{Id} \in \mathfrak{n} \cap \mathfrak{t}.$$ 
\end{lemma}

\begin{proof}
For the first expression to be $0$, we need that $A(S) \subset S.$ For the second expression to be $0$, we need that $$(\mathrm{id}_{H_j}-C_jC_j^*)^{1/2}\Big(\log(Q_j)-\log(A_j)\Big) (\mathrm{id}_{H_j}-C_jC_j^*)^{1/2}=0 \Longrightarrow \Big[\log(Q_j)-\log(A_j)\Big]\Big|_{N_j^\perp}=0.$$ This implies that $$\log(A_j)|_{N_j^\perp}=0,$$ since $Q_j|_{N_j^\perp}=\mathrm{id}_{N_j^\perp}$, or in other words, $$\log(A_j)=\left(\begin{array}{ c | c }
    \log(A_j)|_{N_j} & 0 \\
    \hline
    0 & 0
  \end{array}\right).$$
These are clearly necessary and sufficient conditions for $F(A_1,\cdots,A_m)$ to be equal to $0.$ 
\end{proof}

\hfill \\
Similarly to the previous subsection, we will show that $J:A \mapsto \mathrm{Tr}_S\Big[\log(PAP)|_S \Big]-\mathrm{Tr}\Big( P\log(A)P|_S \Big)$ controls the distance to $\mathfrak{n}$ and that $L:A \mapsto \sum_{j=1}^m p_j \mathrm{Tr}[(C_jC_j^*-\mathrm{id}_{H_j})\log(A_j)]$ controls the distance to $\mathfrak{t}.$

\begin{lemma} \label{Quadratic_SG}
    There exist constants $\rho_0>0, c_0>0$ such that, if $\mathrm{Id}+B>0$ and $\|B\|_p \le \rho_0$, then $$J(\mathrm{Id}+B) \ge c_0 \mathrm{dist}_p(B,\mathfrak{n})^2.$$
\end{lemma}

\begin{proof}
Computing the Hessian of $J$ at the identity, we get $$J(\mathrm{Id}+B)=\frac12 \|PB(\mathrm{Id}-P)\|_{\mathrm{HS}}^2+O(\|B\|_{\mathrm{HS}}^3).$$ Moreover, as $\|PB(\mathrm{Id}-P)\|_{\mathrm{HS}}=0$ implies that $B \in \mathfrak{n}$, there exists $\lambda>0$ such that $$\|PB(\mathrm{Id}-P)\|_{\mathrm{HS}}^2 \ge \lambda \|B\|_{\mathrm{HS}}^2 \text{ for all } B \in \mathfrak{n}^{\perp_p}.$$
\hfill \\
We may write $B=Y+Z$ where $Y \in \mathfrak{n}, Z \in \mathfrak{n}^{\perp_p}$. 
Notice that $J(\mathrm{Id}+Y)=0$ and $J(\mathrm{Id}+Y+sZ) \ge 0$, the first derivative of $s \mapsto J(\mathrm{Id}+Y+sZ)$ at zero is equal to $0.$
Thus, by integration by parts, $$\begin{array}{lll}
J(\mathrm{Id}+B) &=& \ds J(\mathrm{Id}+Y+Z)=\int_0^1 \frac{d}{ds}[J(\mathrm{Id}+Y+sZ)]\, ds \\ \\ 
&=& \ds -\int_0^1 \frac d{ds}(1-s) \frac{d}{ds}[J(\mathrm{Id}+Y+sZ)]\, ds \\ \\ 
&=& \ds \int_0^1 (1-s)\frac{d^2}{ds^2}[J(\mathrm{Id}+Y+sZ)]\, ds=\int_0^1 (1-s) D^2J(\mathrm{Id}+Y+sZ)[Z,Z] \, ds
\end{array}$$
\hfill \\
Since $$D^2J(\mathrm{Id})[Z,Z] \ge \lambda \|Z\|_p^2 \text{ as } Z \in \mathfrak{n}^{\perp_p},$$ and by smoothness of $J$ around the identity, we consider $\rho_0>0$ small enough such that, if $\|X\|_p \le \rho_0$, $$D^2J(\mathrm{Id}+X)[Z,Z] \ge \frac\lambda 2 \|Z\|_p^2.$$ Hence, $$\|Y+sZ\|_p^2=\|Y\|_p^2+s^2\|Z\|_p^2 \le \|Y\|_p^2+\|Z\|_p^2=\|B\|_p^2 \; \Longrightarrow \|Y+sZ\|_p \le \rho_0$$ for any $s \in [0;1]$, which implies that $$D^2J(\mathrm{Id}+Y+sZ)[Z,Z] \ge \frac{\lambda}2 \|Z\|_p^2 \quad \forall s \in [0;1].$$ Therefore, $$J(\mathrm{Id}+B) \ge \frac \lambda 4\|Z\|_p^2 \ge \frac \lambda 4 \mathrm{dist}_p(B,\mathfrak{n})^2.$$
\end{proof}

\begin{lemma} \label{Regularized_Linear_Term}
    There exist constants $\rho_1>0$ and $c_1>0$ such that, if $$A_j:=\mathrm{id}_{H_j}+B_j>0, \quad P_{N_j^\perp}A_jP_{N_j^\perp} \le \mathrm{id}_{N_j^\perp} \text{ and } \|B\|_p \le \rho_1,$$ then $$L(\mathrm{Id}+B) \ge c_1\mathrm{dist}_p(B,\mathfrak{t})^2.$$ In addition, there exists a constant $C_1>0$ such that, for each $1 \le j \le m$, $$\|\mathrm{id}_{N_j^\perp}-P_{N_j^\perp}A_jP_{N_j^\perp}\|_{\mathrm{HS}} \le C_1L(A).$$
\end{lemma}

\begin{proof}

The first claim will be proven block by block. It is quite clear that $$\mathrm{dist}_p(B,\mathfrak{t})^2=\sum_{j=1}^m p_j \Big(\|P_{N_j^\perp}B_jP_{N_j^\perp}\|_{\mathrm{HS}}^2+2\|P_{N_j^\perp}B_jP_{N_j}\|_{\mathrm{HS}}^2\Big).$$ If $N_j^\perp \ne \{0\}$, set $\delta_j>0$ the smallest eigenvalue of $\mathrm{id}_{H_j}-C_jC_j^*$ on $N_j^\perp.$ 
\hfill \\
Choose $\rho_1$ sufficiently small so that $\|B_j\|_{\mathrm{op}} \le \frac12$ for any $1 \le j \le m.$ As $$-\log(1+b) \ge -b+\frac13 b^2 \quad \forall |b| \le \frac12,$$ we get that $$-\log(\mathrm{id}_{H_j}+B_j) \ge -B_j+\frac13 B_j^2 \quad \forall 1 \le j \le m.$$ Since $P_{N_j^\perp}B_jP_{N_j^\perp} \le 0$ and $\mathrm{id}_{H_j}-C_jC_j^*$ has its only non-zero block on $N_j^\perp$, $$\mathrm{Tr}\Big[(C_jC_j^*-\mathrm{id}_{H_j})B_j\Big] \ge 0.$$ Hence, $$\begin{array}{lll}
\ds \mathrm{Tr}\Big[(C_jC_j^*-\mathrm{id}_{H_j})\log(\mathrm{id}_{H_j}+B_j) \Big] &\ge& \ds  \frac13 \mathrm{Tr}\Big[(\mathrm{id}_{H_j}-C_jC_j^*)B_j^2\Big] \ge \frac{\delta_j}{3}\mathrm{Tr}(P_{N_j^\perp}B_j^2)
\\ \\
\nm 
&\ge& \ds \frac{\delta_j}{3}\Big(\|P_{N_j^\perp}B_jP_{N_j^\perp}\|_{\mathrm{HS}}^2+\|P_{N_j^\perp}B_jP_{N_j}\|_{\mathrm{HS}}^2\Big).
\end{array}$$
This concludes the proof of the first assertion of the lemma.

\hfill \\
We move on to the proof of the second claim which is the consequence of the following: $$-\mathrm{Tr}\Big((\mathrm{id}_{H_j}-C_jC_j^*)\log(A_j) \Big) \ge \delta_j \Big\|\mathrm{id}_{N_j^\perp}-P_{N_j^\perp}A_jP_{N_j^\perp}\Big\|_{\mathrm{HS}}.$$ Indeed, this is because $-\log(A_j) \ge \mathrm{id}_{H_j}-A_j$ and $A_j|_{N_j^\perp} \le \mathrm{id}_{N_j^\perp}$, so that $$\begin{array}{lll}
\ds -\mathrm{Tr}\Big( (\mathrm{id}_{H_j}-C_jC_j^*)\log(A_j) \Big) &\ge& \ds \mathrm{Tr}\Big( (\mathrm{id}_{H_j}-C_jC_j^*)(\mathrm{id}_{H_j}-A_j) \Big) \\ 
&\ge& \ds \mathrm{Tr}\Big( (\mathrm{id}_{H_j}-C_jC_j^*)|_{N_j^\perp} (\mathrm{id}_{H_j}-A_j)|_{N_j^\perp} \Big) \\
&\ge& \ds \delta_j \mathrm{Tr}(\mathrm{id}_{N_j^\perp}-A_j|_{N_j^\perp}),
\end{array}$$ and finally using that $\mathrm{Tr}(A)^2 \ge \mathrm{Tr}(A^2)$ for $A \ge 0$. Taking $C_1 \ge \max_{j: N_j^\perp \ne \{0\}}(p_j\delta_j)^{-1}$ concludes the proof of the lemma. 
\end{proof}

\begin{lemma} \label{Sum_Distances_Compatible}
    There exists $c_2>0$ such that for every $B \in \mathcal B$, $$\mathrm{dist}_p(B,\mathfrak{n} \cap \mathfrak{t})^2 \le c_2(\mathrm{dist}_p(B,\mathfrak{n})^2+\mathrm{dist}_p(B,\mathfrak{t})^2).$$
\end{lemma}

\begin{proof}
   By contradiction, if such a statement would not hold, there would exist by homogeneity, $(B_k)_{k \ge 1} \subset \mathcal B$ such that $$\mathrm{dist}_p(B_k,\mathfrak{n} \cap \mathfrak{t})=1, \quad  \mathrm{dist}_p(B_k,\mathfrak{n})+\mathrm{dist}_p(B_k,\mathfrak{t}) \rightarrow_{k \rightarrow +\infty} 0.$$ Subtracting the orthogonal projection onto $\mathfrak{n} \cap \mathfrak{t}$, we may also assume that $B_k \perp_p \mathfrak{n} \cap \mathfrak{t}$ so that $\|B_k\|_p=1$ for any $k \ge 1.$ By compactness, $B_k \rightarrow_{k \rightarrow +\infty} B \in \mathcal B.$ The distance assumptions imply that $B \in \mathfrak{n} \cap \mathfrak{t}, 
   \mathrm{dist}_p(B,\mathfrak{n} \cap \mathfrak{t})=1$ and $\|B\|_p=1.$ Finally, since $B_k \perp_p \mathfrak{n} \cap \mathfrak{t}$, $B \perp_p \mathfrak{n} \cap \mathfrak{t}.$ This implies that $B=0$ in contradiction with $\|B\|_p=1.$
\end{proof}

\hfill \\
Finally, now that we know if $\mathrm{Id}+B$ is near from being a Gaussian extremizer, $B$ is near $\mathfrak{n} \cap \mathfrak{t}$, it still needs to be justified that $B$ must actually be close\footnote{This part is breaking the linear structure, which is why it needs a further argument.}  to $$\mathcal L=\{E \in \mathfrak{n} \cap \mathfrak{t} : -\mathrm{id}_{H_j}<E_j \le Q_j-\mathrm{id}_{H_j}\}.$$ 

\hfill \\
Let $$G_0=\mathrm{span} \{ \mathrm{Ran}(B): B \in \mathfrak{n} \cap \mathfrak{t} \},$$ and $W$ the orthogonal projection onto $G_0.$ We first show that $W$ is the appropriate partial identity along the extremizer directions.

\begin{lemma}
The projection $W$ belongs to $\mathfrak{n} \cap \mathfrak{t}$ and every $B \in \mathfrak{n} \cap \mathfrak{t}$ satisfies $$WB=BW=B.$$
\end{lemma}

\begin{proof}
    Since every $B \in \mathfrak{n} \cap \mathfrak{t}$ is block-diagonal and supported on $N_j$ in each block, $G_0$ is a block direct sum of subspaces of the $N_j$'s. Hence, $W$ is block-diagonal and supported on $\bigoplus_j N_j$, so that $W \in \mathfrak{t}.$
    
    \hfill \\
    Moreover, every $B \in \mathfrak{n}$ is self-adjoint, preserves $S$ and $S^\perp$, so that $\mathrm{Ran}(B)$ is the orthogonal direct sum of a subspace of $S$ and a subspace of $S^\perp.$ Thus, $G_0$ is reduced by the projection $P$ onto $S$, which implies that $W$ commutes with $P$, or equivalently, $W \in  \mathfrak{n}.$
    
    \hfill \\
    Finally, for $B \in \mathfrak{n} \cap \mathfrak{t}$, $\mathrm{Ran}(B) \subset G_0$ implies that $WB=B$ and since $B$ vanishes on $G_0^\perp$,\footnote{Indeed, if $x \in G_0^\perp$, since $B$ is self-adjoint, $$\forall y \in K, \langle Bx,y\rangle=\langle x,By\rangle=0 \Longrightarrow Bx=0.$$} $BW=B.$
\end{proof}

\begin{lemma} \label{Feasible_Extremizer}
    There exists a constant $\kappa<+\infty$ depending only on $R:=(Q_j-\mathrm{id}_{H_j})_{1 \le j \le m}$ and $W$, such that if $Y \in \mathfrak{n} \cap \mathfrak{t}$ and $$R-Y \ge -\delta I_K \text{ for some } 0 \le \delta<1,$$ then $$Y-\kappa \delta W \le R.$$
\end{lemma}

\begin{proof}
    Assume $\delta>0.$ Decomposing $K=WK \oplus (W K)^\perp$, we may write, by the previous lemma, $$Y=\left(\begin{array}{ c | c }
    Y_{11} & 0 \\
    \hline
    0 & 0
  \end{array}\right),(Q_j-\mathrm{id}_{H_j})_{1 \le j \le m}:=R=\left(\begin{array}{ c | c }
    R_{11} & R_{12} \\
    \hline
   R_{12}^* & R_{22}
  \end{array}\right) \ge 0,$$ where the first diagonal block is taken with respect to $WK$ and the second with respect to $(WK)^\perp.$ The assumption $R-Y 
  \ge -\delta I_K$ implies that $$R_{11}-Y_{11}+\delta \mathrm{Id}_{W(K)} \ge R_{12}(R_{22}+\delta \mathrm{Id}_{(WK)^\perp})^{-1}R_{12}^*.$$  Since $R_{22}$ is not necessarily invertible, we will make use the generalized Schur decomposition, where one considers the Moore-Penrose inverse of $R_{22}$, denoted by $R_{22}^\dagger.$ In fact, notice that $$R_{12}R_{22}^\dagger R_{12}^*-R_{12}(R_{22}+\delta \mathrm{Id}_{(WK)^\perp})^{-1}R_{12}^* \le C\delta \mathrm{Id}_{W(K)},$$
  since $R_{22}^\dagger-(R_{22}+\delta \mathrm{Id})^{-1}$ is bounded above by $C\delta \mathrm{Id}$ on $\ker(R_{22})^\perp$, with $C=\max_{\lambda \in \sigma(R_{22})\setminus \{0\}}|\lambda|^{-2}$, while on $\ker(R_{22})$, $R_{22}=0$ and by positivity of $R$, $R_{12}$ is also zero. 
  Thus, $$R_{11}-Y_{11}+(1+C)\delta \mathrm{Id}_{W(K)} \ge R_{12}R_{22}^\dagger R_{12}^*,$$ so that by generalized Schur criterion, for $\kappa \ge 1+C$, $$\begin{pmatrix}
      R_{11}-Y_{11}+\kappa \delta  \mathrm{Id}_{W(K)} & R_{12} \\
      R_{12}^* & R_{22}
  \end{pmatrix} \ge 0,$$ as wished. 
\end{proof}

\hfill \\
By combining the previous lemmas together, we can conclude that one has also quantitative stability for regularized Brascamp-Lieb data at the level of Gaussians. 

\begin{theorem}
In the generalized geometric setting, there exist constants $C<+\infty, \epsilon_0>0$, depending only on the fixed datum $(\textbf{C},p,\textbf{Q})$, so that, if $$0<A_j \le Q_j \quad \forall 1 \le j \le m \text{ and } F(A) \le \epsilon_0$$ then there is a Gaussian extremizer $(E_j)_{1 \le j \le m}, 0<E_j \le Q_j$ and $F(E)=0$, such that $$\sum_{j=1}^m p_j \frac{\|A_j-E_j\|_{\mathrm{HS}}^2}{\|A_j\|_{\mathrm{HS}}^2} \le CF(A),$$ and the corresponding centered Gaussian normalized densities satisfy $$\sum_{j=1}^m \|g_{A_j}-g_{E_j}\|_{L^1} \le CF(A)^{1/2}.$$ 
\end{theorem}

\begin{proof}
Setting $B=e^X-\mathrm{Id}$, $\|B\|_p$ is small for $\|X\|_p$ sufficiently small. Then, assuming that $$P_{N_j^\perp}e^{X_j}P_{N_j^\perp} \le \mathrm{id}_{N_j^\perp} \quad \forall 1 \le j \le m,$$ we get $$\begin{array}{lll}
\ds F(e^X) &=& \ds J(e^X)+L(e^X) \ge c\big( \mathrm{dist}_p(B,\mathfrak{n})^2+\mathrm{dist}_p(B,\mathfrak{t})^2 \Big) \\
\nm 
&\ge& \ds c'\mathrm{dist}_p(B,\mathfrak{n} \cap \mathfrak{t})^2,
\end{array}$$ using \cref{Quadratic_SG}, \cref{Regularized_Linear_Term} and \cref{Sum_Distances_Compatible}. As $$\|e^X-\mathrm{Id}-X\|_p \le C\|X\|_p^2,$$ we find that for $X \perp_p \mathfrak{n} \cap \mathfrak{t}$ with $\|X\|_p$ small enough, $$\mathrm{dist}_p(B,\mathfrak{n}\cap \mathfrak{t}) \ge \mathrm{dist}_p(X,\mathfrak{n} \cap \mathfrak{t})-\|e^X-\mathrm{Id}-X\|_p \ge \|X\|_p-C\|X\|_p^2 \ge \frac12 \|X\|_p.$$
\hfill \\
Hence, we obtain that $$F(e^X) \ge c\|X\|_p^2 \quad \text{ if } X \perp_p \mathfrak{n} \cap \mathfrak{t},\quad \|X\|_p \le \rho \text{ and } P_{N_j^\perp}e^{X_j}P_{N_j^\perp} \le \mathrm{id}_{N_j^\perp},$$ for any $1 \le j \le m.$
\hfill \\
Furthermore, for any $1 \le j \le m$, if $P_{N_j^\perp}e^{X_j}P_{N_j^\perp} \le \mathrm{id}_{N_j^\perp}$, $$P_{N_j^\perp}e^{tX_j}P_{N_j^\perp} \le \mathrm{id}_{N_j^\perp} \quad \forall t \in [0;1].$$
By convexity of $t \mapsto F(e^{tX})$, the previous local inequality extends to $$F(e^X) \ge c\min\{\|X\|_p,\|X\|_p^2\} \quad \text{ if } X \perp_p \mathfrak{n} \cap \mathfrak{t}, \quad P_{N_j^\perp}e^{X_j}P_{N_j^\perp} \le \mathrm{id}_{N_j^\perp},$$ for any $1 \le j \le m.$ 
Similarly to \cref{Orthogonal_Minimizer}, there exists $E=I+B>0$ such that $B \in \mathfrak{n} \cap \mathfrak{t}$ and $$e^X=E^{-1/2}AE^{-1/2} \text{ for some } X \perp_p \mathfrak{n} \cap \mathfrak{t}, \quad d_p(A,E)=\|X\|_p.$$ 
Then, $$F(E)=0 \text{ and } F(e^X)=F(E^{-1/2}AE^{-1/2})=F(A).$$ Moreover, $$P_{N_j^\perp}e^{X_j}P_{N_j^\perp}=P_{N_j^\perp}E_j^{-1/2}A_jE_j^{-1/2}P_{N_j^\perp}=P_{N_j^\perp}A_jP_{N_j^\perp}\le \mathrm{id}_{N_j^\perp},$$ for any $1 \le j \le m.$
In particular, $F(A) \le \epsilon_0$ implies $\|X\|_p \le C\sqrt{F(A)}.$ 

\hfill \\
Note that \begin{equation} \label{FromXtoA} e^{-C\|X\|_p}E \le A \le e^{C\|X\|_p}E \text{ and in particular, } E \le e^{C\|X\|_p}Q,\end{equation} where $Q:=(Q_j)_{1 \le j \le m}.$ Hence, $$Q-E \ge -(e^{C\|X\|_p}-1)Q \ge -C_Q \|X\|_p \, \mathrm{Id},$$ using again that $\epsilon_0$ may be chosen sufficiently small which makes $\|X\|_p$ very small as well. 

\hfill \\
Applying \cref{Feasible_Extremizer} to $Y=E-\mathrm{Id}$, $$E-uW \le Q \text{ for some } u=O_Q(\|X\|_p).$$ Considering directly $E-uW$ could remove the needed positivity assumption, so that we define $$\tilde E=\frac{E+u(\mathrm{Id}-W)}{1+u}>0.$$ Moreover, $\tilde E-\mathrm{Id} \in \mathfrak{n} \cap \mathfrak{t}$ and $$Q-\tilde E=\frac{Q-(E-uW)+u(Q-\mathrm{Id})}{1+u} \ge 0,$$ so that $\tilde E$ is also a feasible extremizer.
Then, as $E-\mathrm{Id}$ is supported on $W$ by definition, $$E=(\mathrm{Id}-W)+EW \Longrightarrow \frac1{1+u}E \le \tilde E \le E,$$ with $u=O_Q(\|X\|_p).$
Applying the previous inequality together with \cref{FromXtoA} and $\|X\|_p \le C\sqrt{F(A)}$ provides $$-C\sqrt{F(A)}A \le \tilde E-A \le C\sqrt{F(A)}A,$$ and in particular, the Hilbert-Schmidt inequality. The $L^1$ inequality is a consequence of \cref{L1Gaussian_Behaviour}. 
\end{proof}

\hfill \\ 
To include uncentered Gaussians to our analysis, we define for $A=(A_j)_{1 \le j \le m}$ and $a=(a_j)_{1 \le j \le m} \subset \prod_{j=1}^m H_j$, $$Q_A(a)=\inf_{x \in H} \sum_{j=1}^m p_j \|A_j^{1/2}(C_jx-a_j)\|^2_{H_j}.$$ Define $$\overline{x}_{A,a}=\Big(\sum_{j=1}^m p_jC_j^*A_jC_j\Big)^{-1}\sum_{j=1}^m p_j C_j^*A_ja_j.$$ Then, $$Q_A(a)=\sum_{j=1}^m p_j \|A_j^{1/2}(C_j\overline{x}_{A,a}-a_j)\|_{H_j}^2,$$ by for example checking the gradient of the function considered. In particular, \begin{equation} \label{Square_Completion_UGaussian} \begin{array}{lll}
\ds Q_A(a)=\sum_{j=1}^m p_j \langle A_jC_j\overline{x}_{A,a}-A_ja_j,C_j\overline{x}_{A,a}-a_j\rangle=\sum_{j=1}^m p_j \langle A_ja_j,a_j\rangle
+\Big\langle \sum_{j=1}^m p_jC_j^*A_jC_j\overline{x}_{A,a},\overline{x}_{A,a} \Big\rangle \\ \\
\ds -2 \Big\langle \sum_{j=1}^m p_jC_j^*A_ja_j, \overline{x}_{A,a}\Big\rangle=\sum_{j=1}^m p_j \langle A_ja_j,a_j\rangle-\Big\langle \sum_{j=1}^m p_jC_j^*A_jC_j\overline{x}_{A,a},\overline{x}_{A,a} \Big\rangle.
\end{array} \end{equation}

\hfill \\
We write $$g_{A_j,a_j}(y_j)=g_{A_j}(y_j-a_j).$$
Then, $$\int_{H}\prod_{j=1}^m g_{A_j,a_j}(C_jx)^{p_j}\, dx=\int_H \prod_{j=1}^m g_{A_j}(C_jx)^{p_j}\, dx \cdot \exp\Big(-\frac12 Q_A(a) \Big),$$ since $$\sum_{j=1}^m p_j \langle A_j(C_jx-a_j),C_jx-a_j\rangle=\Big\langle \sum_{j=1}^m p_jC_j^*A_jC_j x,x\Bigl \rangle-2 \Bigl\langle \sum_{j=1}^m p_jC_j^*A_ja_j,x\Big\rangle+\sum_{j=1}^m p_j \langle A_ja_j,a_j\rangle,$$ so that completing the square, $$\sum_{j=1}^m p_j \langle A_j(C_jx-a_j),C_jx-a_j\rangle=\Big\langle \sum_{j=1}^m p_j C_j^*A_jC_j(x-\overline{x}_{A,a}),x-\overline{x}_{A,a}\Bigl\rangle+Q_A(a),$$ by \cref{Square_Completion_UGaussian}.

\hfill \\
Furthermore, note that for any finite dimensional Hilbert space $H$, there exists a constant $C_H>0$, such that, for any positive definite matrix $A$ and any $a,b \in H$, \begin{equation} \label{Translation_SameGaussian} \|g_{A,a}-g_{A,b}\|_{L^1} \le C_H\min\{1,\|A^{1/2}(a-b)\|\}.\end{equation}

\hfill \\
Indeed, by the affine change of variables $A^{1/2}(y-b)$, it suffices to establish such a bound for $$\|g_{\mathrm{Id},h}-g_{\mathrm{Id},0}\|_{L^1},$$ where $h=A^{1/2}(a-b).$ A trivial bound is $2$ and for the other one, it suffices to recall that $$g_{\mathrm{Id}}(y-h)-g_{\mathrm{Id}}(y)=-\int_0^1 \nabla g_{\mathrm{Id}}(y-th) \cdot h\, dt,$$ so that $$\|g_{\mathrm{Id},h}-g_{\mathrm{Id},0}\|_{L^1} \le |h|\int_H |\nabla g_{\mathrm{Id}}(y)|\, dy.$$

\hfill \\
In particular, this allows to conclude Gaussian precompactness, \textit{i.e.} if $$g_{A_1^{(n)},a_1^{(n)}},\cdots,g_{A_m^{(n)},a_m^{(n)}}$$ is near-extremizing (in other words, their Brascamp-Lieb deficit converge to $0$), then there exists a sequence of extremizers $$g_{B_1^{(n)},C_1x_n},\cdots, g_{B_m^{(n)},C_mx_n},$$ for some $x_n \in H$, such that $$\sum_{j=1}^m \|g_{A_j^{(n)},a_j^{(n)}}-g_{B_j^{(n)},C_jx_n}\|_{L^1} \rightarrow_{n \rightarrow +\infty}0.$$ 
\hfill \\
Indeed, by the expression we derived for the Brascamp-Lieb deficit of uncentered Gaussians, this implies that $g_{A_1^{(n)}},\cdots,g_{A_m^{(n)}}$ is near-extremizing and $Q_{A^{(n)}}(a^{(n)}) \rightarrow_{n \rightarrow +\infty}0.$ We know that there must exist a sequence of extremizers $(g_{B_1^{(n)}},\cdots,g_{B_m^{(n)}})$ such that $$\sum_{j=1}^m \|g_{A_j^{(n)}}-g_{B_j^{(n)}}\|_{L^1} \rightarrow_{n \rightarrow +\infty}0.$$ Taking $x_n=\overline x_{A^{(n)},a^{(n)}}$, we have by \cref{Translation_SameGaussian} that $$\sum_{j=1}^m \|g_{A_j^{(n)},C_jx_n}-g_{A_j^{(n)},a_j^{(n)}}\|_{L^1} \le C\sum_{j=1}^m \big \|(A_j^{(n)})^{1/2} (a_j^{(n)}-C_jx_n) \big \|_{H_j} \rightarrow_{n \rightarrow +\infty }0.$$ Hence, by translation invariance, $$\sum_{j=1}^m \|g_{A_j^{(n)},a_j^{(n)}}-g_{B_j^{(n)},C_jx_n}\|_{L^1} \le \sum_{j=1}^m \|g_{A_j^{(n)},a_j^{(n)}}-g_{A_j^{(n)},C_jx_n}\|_{L^1}+\sum_{j=1}^m \|g_{A_j^{(n)}}-g_{B_j^{(n)}}\|_{L^1} \rightarrow_{n \rightarrow +\infty}0.$$ It is clear that also $(g_{B_j^{(n)},C_jx_n})_{1 \le j \le m}$ is a sequence of extremizers.

\hfill \\
For quantitative stability, we argue analogously using the previous remark. In the following theorem, we may also take $Q_j=+\infty$ thanks to the previous subsection, and by change of variables, any extremizable Brascamp-Lieb datum $(\textbf{C},p,\textbf{Q})$ may be considered. 

\begin{theorem}
Let $(\textbf{C},p,\textbf{Q})$ be an extremizable Brascamp-Lieb datum. 
There exist constants $\epsilon_0>0$ and $C<+\infty$ so that, if $$0<A_j \le Q_j, \quad a_j \in H_j \quad \forall 1 \le j \le m,$$ and $$\mathrm{BL}(\textbf{C},p,\textbf{Q})-\frac{\int_{\mathbb{R}^n}\prod_{j=1}^m g_{A_j,a_j}(C_jx)^{p_j}\, dx}{\prod_{j=1}^m (\int_{\mathbb{R}^{n_j}}g_{A_j,a_j})^{p_j}} \le \epsilon\le \epsilon_0,$$ then there exists a Gaussian extremizer $(E_j)_{1 \le j \le m}, 0<E_j \le Q_j$ and $x \in H$, such that $$\sum_{j=1}^m \|g_{A_j,a_j}-g_{E_j,C_jx}\|_{L^1} \le  \; C\sqrt{\epsilon}.$$ 
\end{theorem}

\section{The spectral gap}
\label{Sec:3}

In this section, we treat the general local stability around the standard Gaussian for a geometric Brascamp-Lieb datum. The main difficulty lies in identifying correctly the kernel of the quadratic form obtained after Taylor expansion, whereas the spectral gap is easily obtained.

\hfill \\
Let us first explain why we may reduce to $(\textbf{C},p)$ being a geometric Brascamp-Lieb datum. Indeed, if \cref{Complete_Gaussian_stability} holds for any geometric Brascamp-Lieb datum, we may use the notion of equivalence and the fact that any extremizable Brascamp-Lieb datum $(\textbf{C},p)$ is equivalent to a geometric one $(\textbf{C'},p)$ (see the discussion in \cref{Sec1.1}). More precisely, recall this means that there exist linear invertible transformations $$D:H \rightarrow H', D_j: H_j \rightarrow H_j' \text{ s.t. } C_j=D_j^{-1}C_j'D \quad \forall 1 \le j \le m,$$ with \begin{equation} \label{Equiv_BL_From_Geo}\mathrm{BL}(\textbf{C},p)=\frac{\prod_{j=1}^m |\det(D_j)|^{p_j}}{|\det(D)|}.\end{equation} Then, by using \cref{Complete_Gaussian_stability} for the geometric datum $(\textbf{C'},p)$ and the tuple $(f_j \circ D_j^{-1})_{1 \le j \le m}$, we obtain $(g_j)_{1 \le j \le m}$ such that $\int_{H_j'} g_j=1$ for any $1 \le j \le m$ and \begin{equation} \label{Stability_Geometric_BL_ChangeOfV} 1-\frac{\ds \int_{H'}\prod_{i=1}^m f_i(C_iD^{-1}x)^{p_i}\, dx }{\ds \prod_{i=1}^m\Big(\int_{H_i'}f_i \circ D_i^{-1}\, \Big)^{p_i}}=1-\frac{\ds \int_{H'}\prod_{i=1}^m f_i(D_i^{-1}C_i'x)^{p_i}\, dx }{\ds \prod_{i=1}^m\Big(\int_{H_i'}f_i \circ D_i^{-1}\, \Big)^{p_i}} \ge c \sum_{i=1}^m p_i \Big\|\sqrt{\frac{f_i\circ D_i^{-1}}{\|f_i \circ D_i^{-1}\|_{L^1}}}-\sqrt{g_i}\Big\|_{L^2}^2.\end{equation} Setting $\tilde g_j=|\det D_j|g_j\circ D_j$ for any $1 \le j \le m$, one has that $\int_{H_j}\tilde g_j=1$ and by \cref{Stability_Geometric_BL_ChangeOfV} together with \cref{Equiv_BL_From_Geo}, $$1-\frac{\ds \int_H \prod_{i=1}^m f_i(C_ix)^{p_i}}{\ds \mathrm{BL}(\textbf{C},p)\prod_{i=1}^m\Big(\int_{H_i}f_i\Big)^{p_i}} \ge c\sum_{i=1}^m p_i\Big\|\sqrt{\frac{f_i}{\|f_i\|_{L^1}}}-\sqrt{\tilde g_i}\Big\|_{L^2}^2.$$

\hfill \\
From now on, assume that $(\textbf{C},p)$ is geometric. Recall that $H \cong \mathbb{R}^n, H_i \cong \mathbb{R}^{n_i}$ for any $1 \le i \le m$ and as $(\textbf{C},p)$ is a geometric Brascamp-Lieb datum, we may assume that $H_i \subset H$ for any $1 \le i \le m.$ Set $$\begin{array}{lll}
F(\varphi_1,\cdots,\varphi_m) &=& \ds \sum_{i=1}^m p_i \log \Big(\int_{\mathbb{R}^{n_i}} e^{-|x_i|^2/2- \varphi_i(x_i)}  \, dx_i \Big)-\log \Big( \int_{\mathbb{R}^n} \prod_{i=1}^m \big(e^{-|C_ix|^2/2- \varphi_i(C_ix)} \big)^{p_i}\, dx \Big) \\ \\
&=& \ds \sum_{i=1}^m p_i \log \Big(\int_{\mathbb{R}^{n_i}} e^{-\varphi_i(x_i)}  \, d\gamma_{n_i}(x_i) \Big)-\log \Big( \int_{\mathbb{R}^n} e^{- \sum_{i=1}^m p_i \varphi_i(C_ix)}\, d\gamma_n(x) \Big) \ge 0.
\end{array}$$

\hfill \\
We also used here the homogeneity condition to get a cancellation of the constants introduced by considering Gaussian measures, \textit{i.e.} $$\sum_{i=1}^m p_i n_i=n.$$
Note that in the previous section, we proved that if $\varphi_i$ is quadratic, one has the desired stability statement of order $2$.

\hfill \\
For $f$ bounded, one has that $$\begin{array}{lll}
\ds \log \Big( \int_{\mathbb{R}^m} e^{-\epsilon f}\, d\gamma_m \Big) =  \log \Big( 1-\epsilon \int_{\mathbb{R}^m}f \, d\gamma_m +\frac{\epsilon^2}2 \int_{\mathbb{R}^m}f^2 \, d\gamma_m +O(\epsilon^3) \Big) \\ \\ 
= \ds -\epsilon \int_{\mathbb{R}^m}f\, d\gamma_m +\frac{\epsilon^2}2 \int_{\mathbb{R}^m}f^2 \, d\gamma_m-\frac12 \Big( -\epsilon \int_{\mathbb{R}^m}f\, d\gamma_m +\frac{\epsilon^2}2 \int_{\mathbb{R}^m}f^2 \, d\gamma_m\Big)^2+O(\epsilon^3) \\ \\ 
\ds =-\epsilon \int_{\mathbb{R}^m} f\, d\gamma_m+\frac{\epsilon^2}2 \mathrm{Var}_{\gamma_m}(f)+O(\epsilon^3). 
\end{array}$$
\hfill \\
Hence, for $\varphi_1,\cdots,\varphi_m$ bounded, $$\begin{array}{lll}
\ds F(\epsilon \varphi_1,\cdots,\epsilon \varphi_m)=-\epsilon \sum_{i=1}^m p_i \int_{\mathbb{R}^{n_i}}\varphi_i \, d\gamma_{n_i}+\frac{\epsilon^2}{2}\sum_{i=1}^m p_i \mathrm{Var}_{\gamma_{n_i}}(\varphi_i) +\epsilon \int_{\mathbb{R}^n}\sum_{i=1}^m p_i \varphi_i(C_ix)\, d\gamma_n(x) \\\\ 
\ds -\frac{\epsilon^2}2 \mathrm{Var}_{\gamma_n}\Big(\sum_{i=1}^m p_i \varphi_i \circ C_i\Big) +O(\epsilon^3).
\end{array}$$ 
\hfill \\
By geometricity, $C_i: H \rightarrow H_i$ is a projection with $H_i \subset H$ and $H_i^\perp$ its orthogonal complement.
Using Fubini (or more directly, as $(C_i)_{\#}\gamma_n=\gamma_{n_i}$), $$
\ds \int_{\mathbb{R}^n} \sum_{i=1}^m p_i \varphi_i(C_ix)\, d\gamma_n(x)=(2\pi)^{-n/2}\sum_{i=1}^m p_i \int_{H_i^\perp} \Big(\int_{H_i} \varphi_i(y)e^{-\frac{|y|^2}2}\, dy \Big) e^{-\frac{|z|^2}2}\, dz=\sum_{i=1}^m p_i \int_{\mathbb{R}^{n_i}}\varphi_i \, d\gamma_{n_i}.$$

\hfill \\
Now, let us try to understand why the inequality $$\sum_{i=1}^m p_i \mathrm{Var}_{\gamma_{n_i}}(\varphi_i) \ge \mathrm{Var}_{\gamma_n}\Big(\sum_{i=1}^m p_i \varphi_i \circ C_i \Big) \quad \forall \varphi_i \in L^2(\gamma_{n_i}),$$ holds and what are the potential extremizers. Here, the resulting variance inequality is extended to $L^2$ functions by approximation.
Again, in the Gaussian case, $\varphi_i(x_i)=\langle A_ix_i,x_i \rangle$ for any $x_i \in H_i$, so that the inequality reduces to $$\sum_{i=1}^m p_i \mathrm{Tr}(A_i^2) \ge \mathrm{Tr}\Bigg( \Big(\sum_{i=1}^m p_i C_i^*A_iC_i\Big)^2 \Bigg).$$

\hfill \\
In order to investigate this inequality, it is worth introducing Hermite polynomials. With the standard basis $(e_1,\cdots,e_n)$, $u \in \mathrm{Sym}^k(\mathbb{R}^n)$ if $$u=\sum_{j_1,\cdots,j_k=1}^n u_{j_1\cdots j_k}e_{j_1} \otimes \cdots \otimes e_{j_k},$$ for coefficients that must be invariant under permutations. 
On such a space, we may define the inner product as $$\langle u,v\rangle=\sum_{j_1,\cdots,j_k=1}^n u_{j_1\cdots j_k}v_{j_1\cdots j_k} \quad \forall u,v \in \mathrm{Sym}^k(\mathbb{R}^n).$$ 
\hfill \\
In addition, for any $u \in \mathrm{Sym}^k(\mathbb{R}^n)$, set $$I_k^{(n)}(u)(x)=\sum_{j_1,\cdots,j_k=1}^n u_{j_1\cdots j_k}(-1)^k e^{|x|^2/2} \frac{\partial^k}{\partial x_{j_1}\cdots \partial x_{j_k}}e^{-|x|^2/2}.$$ 
\hfill \\
By definition, for any $x \in \mathbb{R}^n$, $$I_k^{(n)}(t^{\otimes k})(x)=(-1)^k e^{|x|^2/2} \sum_{j_1,\cdots,j_k=1}^nt_{j_1}\cdots t_{j_k} \frac{\partial^k}{\partial x_{j_1}\cdots \partial x_{j_k}}e^{-|x|^2/2}=(-1)^k e^{|x|^2/2} (t \cdot \nabla)^k e^{-|x|^2/2},$$ for $t=(t_1,\cdots,t_n) \in \mathbb{R}^n$ and $$\frac{d^k}{ds^k}e^{-|x-st|^2/2}\Big|_{s=0}=(-1)^k (t \cdot \nabla)^k e^{-|x|^2/2},$$ which yields $$e^{|x|^2/2}e^{-|x-st|^2/2}=e^{|x|^2/2}\sum_{k=0}^\infty \frac{s^k}{k!} \Big(\frac{d^k}{ds^k}e^{-|x-st|^2/2}\Big|_{s=0}\Big)=\sum_{k=0}^\infty \frac{s^k}{k!} I_k^{(n)}(t^{\otimes k})(x) \quad \forall s \in \mathbb{R}.$$ In particular, $$e^{\langle x,t\rangle-|t|^2/2}=\sum_{k=0}^\infty \frac1{k!}I_k^{(n)}(t^{\otimes k})(x) \quad \forall x,t \in \mathbb{R}^n.$$  

\hfill \\
Observe that $$\int_{\mathbb{R}^n}e^{\langle x,t_1\rangle-|t_1|^2/2} e^{\langle x,t_2\rangle-|t_2|^2/2} \, d\gamma_n(x)=(2\pi)^{-n/2}\int_{\mathbb{R}^n}e^{-|x-(t_1+t_2)|^2/2} e^{\langle t_1,t_2\rangle}\, dx=e^{\langle t_1,t_2\rangle}.$$ Consequently, $$\sum_{k,j=0}^\infty \frac1{k!}\frac1{j!} \int I_k^{(n)}(t_1^{\otimes k})(x) I_j^{(n)}(t_2^{\otimes j})(x)\, d\gamma_n(x)=e^{\langle t_1,t_2\rangle}=\sum_{k=0}^\infty \frac1{k!}\langle t_1^{\otimes k},t_2^{\otimes k}\rangle,$$ and so, $$\int I_k^{(n)}(t_1^{\otimes k})(x) I_j^{(n)}(t_2^{\otimes j})(x)\, d\gamma_n(x)=\delta_{kj}\, k!\, \langle t_1^{\otimes k},t_2^{\otimes k}\rangle.$$ 

\hfill \\
One has that $\mathrm{Sym}^k(\mathbb{R}^n)=\mathrm{span} \{ x^{\otimes k} :x \in \mathbb{R}^n \}$ or in other words, any $k$-tensor can be expressed a linear combination of ``rank-one'' tensors $x^{\otimes k},$ which we will call the polarization principle for symmetric tensors. 
Hence, $$\int I_k^{(n)}(u)I_l^{(n)}(v)\, d\gamma_n=\delta_{kl}\, k!\, \langle u,v\rangle,$$ as this equality holds true for ``rank-one'' tensors and by the polarization principle.

\hfill \\
Furthermore, $$L^2(\gamma_n)=\bigoplus_{k=0}^{\infty} I_k^{(n)}(\mathrm{Sym}^k(\mathbb{R}^n)),$$ or in other words, every $f \in L^2(\gamma_n)$ has a unique expansion $$f=\sum_{k=0}^\infty I_k^{(n)}(u_k) \text{ in } L^2(\gamma_n),$$ and $$u_0=\int f \, d\gamma_n, \quad \|f\|_{L^2(\gamma_n)}^2=\sum_{k=0}^\infty k! \|u_k\|^2.$$

\begin{lemma} \label{Composition_Lemma}
    Let $C: \mathbb{R}^n \rightarrow \mathbb{R}^r$ such that $CC^*=\mathrm{Id}_{\mathbb{R}^r}$. Then, for any $u \in \mathrm{Sym}^k(\mathbb{R}^r)$, $$I_k^{(r)}(u) \circ C=I_k^{(n)}((C^*)^{\otimes k}u).$$
\end{lemma}

\begin{proof}
Note that, for any $t \in \mathbb{R}^r$, $$\sum_{k=0}^\infty \frac1{k!}I_k^{(r)}(t^{\otimes k})(Cx)=e^{\langle Cx,t\rangle-|t|^2/2}=e^{\langle x,C^*t\rangle-|C^*t|^2/2}=\sum_{k=0}^\infty \frac1{k!}I_k^{(n)}((C^*)^{\otimes k} t^{\otimes k})(x).$$ Thus, by uniqueness of the decomposition, we get equality of each $I_k$ term so that the polarization principle yields the desired claim. 
\end{proof}

\begin{lemma} \label{Conditional_Exp_Lemma}
    If $P: H \rightarrow V$ is a projection on a subspace $V \subset H$ and $X \sim \gamma_H$, then for any $h \in \mathrm{Sym}^k(H)$, $$\mathbb{E}\Big[I_k^{(\dim H)}(h)(X)| PX \Big]=I_k^{(\dim H)}(P^{\otimes k}h)(X).$$
\end{lemma}

\begin{proof}
    As usual, it suffices to establish the claim for ``rank-one'' tensors. Observe that $PX$ and $(\mathrm{Id}-P)X$ are independent since uncorrelated and jointly Gaussian as linear image of a Gaussian vector. Thus, $$\begin{array}{lll}
     \mathbb{E}(e^{\langle t,X\rangle-|t|^2/2}|PX) &=& \ds \mathbb{E}(e^{\langle t,PX\rangle-|Pt|^2/2} e^{\langle t,(\mathrm{Id}-P)X\rangle-|(\mathrm{Id}-P)t|^2/2}|PX) \\ \\
     &=& \ds \mathbb{E}(e^{\langle t,(\mathrm{Id}-P)X\rangle-|(\mathrm{Id}-P)t|^2/2}) e^{\langle t,PX\rangle-|Pt|^2/2} \\ \\
     &=& e^{\langle Pt,X\rangle-|Pt|^2/2},
    \end{array}$$ so that comparing the Taylor expression term by term concludes the proof of the lemma. 
\end{proof}

\hfill \\
By using Hermite polynomials, we confirm that the quadratic form is non-negative and give an abstract description of its kernel. 

\begin{theorem} \label{Quadratic_Form_GeneralCase}
    For every $\varphi_i \in L^2(\gamma_{n_i})$, \begin{equation} \label{Gaussian_Variance_BL} \sum_{i=1}^m p_i \mathrm{Var}_{\gamma_{n_i}}(\varphi_i) \ge  \mathrm{Var}_{\gamma_n}\Big(\sum_{i=1}^m p_i \varphi_i \circ C_i \Big).\end{equation}
   There is equality if and only if $$\varphi_i-\int \varphi_i \, d\gamma_{n_i}=\sum_{k=1}^\infty I_k^{(n_i)}(C_i^{\otimes k}h^{(k)}) \text{ in } L^2(\gamma_{n_i}),$$ where $h^{(k)} \in \ker(\mathrm{Id}-T_kT_k^*)$ for any $k \ge 1$ and such that $\sum_{k=1}^{\infty} k!\|h^{(k)}\|^2<+\infty.$ Here, $$T_k: \bigoplus_{i=1}^m \mathrm{Sym}^k(\mathbb{R}^{n_i}) \rightarrow \mathrm{Sym}^k(\mathbb{R}^n), \quad T_k(u)=\sum_{i=1}^m p_i (C_i^*)^{\otimes k}u_i.$$ 
\end{theorem}

\begin{proof}
    First, we may assume that $\varphi_i$ are centered for any $1 \le i \le m.$ Then, $$\varphi_i=\sum_{k=1}^\infty I_k^{(n_i)}(u_i^{(k)}) \text{ in } L^2(\gamma_{n_i}),$$ and $$\mathrm{Var}_{\gamma_{n_i}}(\varphi_i)=\|\varphi_i\|_{L^2(\gamma_{n_i})}^2=\sum_{k=1}^\infty k! \|u_i^{(k)}\|^2,$$ for any $1 \le i \le m.$ Moreover, using \cref{Composition_Lemma}, $$\mathrm{Var}_{\gamma_n}\Big(\sum_{i=1}^m p_i \varphi_i \circ C_i \Big)=\int \Big( \sum_{k=1}^\infty \sum_{i=1}^m p_i I_k^{(n)}((C_i^*)^{\otimes k}u_i^{(k)}) \Big)^2 \, d\gamma_n=\sum_{k=1}^\infty k!\Big\| \sum_{i=1}^m p_i (C_i^*)^{\otimes k}u_i^{(k)}\Big\|^2.$$ In consequence, it suffices to show that for any $k \ge 1$, $$\sum_{i=1}^m p_i \|u_i^{(k)}\|^2 \ge \Big\| \sum_{i=1}^m p_i (C_i^*)^{\otimes k}u_i^{(k)} \Big\|^2.$$ To prove this, the argument resembles the one used in \cref{QF_Degree2_Rk}. More precisely, we introduce a weighted inner product on $$U_k=\bigoplus_{i=1}^m \mathrm{Sym}^k(\mathbb{R}^{n_i}),\; \langle u,v\rangle_p=\sum_{i=1}^m p_i \langle u_i,v_i \rangle.$$ Define $$T_k: U_k \rightarrow \mathrm{Sym}^k(\mathbb{R}^n), \quad T_k(u)=\sum_{i=1}^m p_i (C_i^*)^{\otimes k}u_i.$$ Then, $$\langle T_k(u),v\rangle=\sum_{i=1}^m p_i \langle u_i,C_i^{\otimes k}v\rangle \Longrightarrow T^*_k(v)=(C_i^{\otimes k}v)_{1 \le i \le m}.$$ Moreover, since $(C_i^*C_i)^{\otimes k} \le C_i^*C_i \otimes \mathrm{Id}_{n}^{\otimes (k-1)}$, one has that $$T_kT_k^*=\sum_{i=1}^m p_i (C_i^*C_i)^{\otimes k} \le \sum_{i=1}^m p_iC_i^*C_i \otimes \mathrm{Id}_n^{\otimes (k-1)}=\mathrm{Id}_n^{\otimes k} \Longrightarrow \|T_k^*v\|_p^2 \le \|v\|^2,$$ for any $v \in \mathrm{Sym}^k(\mathbb{R}^n).$ This implies that $$\|T_ku\|^2 \le \|u\|_p^2 \quad \forall u \in U_k,$$ and yields the desired claim. 
    Note that there is equality in the variance inequality if $$\sum_{i=1}^m p_i \|u_i^{(k)}\|^2=\Big\|\sum_{i=1}^m p_i (C_i^*)^{\otimes k}u_i^{(k)}\Big\|^2,$$ for any $k \ge 1.$  This implies that $$\|T_ku^{(k)}\|=\|u^{(k)}\|_p \Longrightarrow u^{(k)}=T_k^*T_k u^{(k)},$$ for any $k \ge 1.$ In particular, setting $h^{(k)}=T_ku^{(k)}$, one has that $h^{(k)} \in \ker(\mathrm{Id}-T_kT_k^*)$ so that $u^{(k)}=T_k^*h^{(k)},$ which proves the desired claim about the extremizers. 
\end{proof}

\hfill \\
The goal is now to have a more concrete characterization of these special Hermite polynomials.  

\hfill \\
If $E \subset H$, then $$\mathrm{Sym}^k(H)=\bigoplus_{r=0}^k \Big( \mathrm{Sym}^r(E) \odot \mathrm{Sym}^{k-r}(E^\perp) \Big),$$ where $\odot$ is the symmetrized tensor product. 

\begin{lemma} \label{Description_Kernel_GeneralCase}
For every $k \ge 1$, $$\Big\{ h \in \mathrm{Sym}^k(H): \sum_{i=1}^m p_i(C_i^*C_i)^{\otimes k}h=h \Bigl\}=\ker(\mathrm{Id}-T_kT_k^*)=\bigcap_{i=1}^m \Big(\mathrm{Sym}^k(H_i) \oplus \mathrm{Sym}^k(H_i^\perp)\Big).$$
\end{lemma}

\begin{remark} \label{Case1_2_AlreadyDone}
    Note that for $k=1$, there is nothing to show (since $T_1T_1^*=\mathrm{Id}_n$) and for $k=2$, the result of this lemma was already proved in \cref{QF_Degree2_Rk}.
\end{remark}

\begin{proof}[Proof of \cref{Description_Kernel_GeneralCase}]
    Consider $h \in \mathrm{Sym}^k(H)$ and fix $1 \le i \le m.$ Then, $$h=\sum_{r=0}^k h_r, \text{ where } h_r \in \mathrm{Sym}^r(H_i) \odot \mathrm{Sym}^{k-r}(H_i^\perp).$$ We define for any $1 \le j \le k$, $$(C_i^*C_i)^{(j)}=\mathrm{Id}^{\otimes (j-1)} \otimes (C_i^*C_i) \otimes \mathrm{Id}^{\otimes (k-j)},$$ or in other words, for any $x_1,\cdots, x_k \in H$, $$(C_i^*C_i)^{(j)} x_1 \otimes \cdots \otimes x_k=x_1 \otimes \cdots \otimes C_i^*C_ix_j \otimes \cdots \otimes x_k.$$ If $x_j \in H_i$, $$(C_i^*C_i)^{(j)}x_1 \otimes \cdots \otimes x_k=x_1 \otimes \cdots \otimes x_k,$$ while if $x_j \in H_i^\perp$, $$(C_i^*C_i)^{(j)}x_1 \otimes \cdots \otimes x_k=0.$$ In particular, this implies that $$\sum_{j=1}^k (C_i^*C_i)^{(j)}h_r=rh_r.$$ By symmetry, $\langle (C_i^*C_i)^{(j)}h_r,h_r\rangle$ are equal for all $1 \le j \le k$, so that $$\langle (C_i^*C_i \otimes \mathrm{Id}^{\otimes(k-1)})h_r,h_r\rangle =\langle (C_i^*C_i)^{(1)} h_r,h_r\rangle=\frac rk \|h_r\|^2.$$ Moreover, \begin{equation} \begin{aligned} (C_i^*C_i)^{\otimes k}h_r=\begin{cases} 
    0, & r<k, \\
    h_k, & r=k.
    \end{cases} \end{aligned} \end{equation} Recall that $$\mathrm{Id}_n^{\otimes k}-T_kT_k^*=\sum_{i=1}^m p_i ((C_i^*C_i)^{(1)}-(C_i^*C_i)^{\otimes k}) \ge 0.$$ In particular, if $h \in \mathrm{ker}(\mathrm{Id}-T_kT_k^*)$ and as the spaces $\mathrm{Sym}^r(H_i) \odot \mathrm{Sym}^{k-r}(H_i^\perp)$ are orthogonal for distinct $r$, $$\Big \langle \big[ (C_i^*C_i)^{(1)}-(C_i^*C_i)^{\otimes k} \big]h_r,h_r \Big\rangle =0 \quad \forall 0 \le r \le k.$$ The reverse direction is of course also true.
    Considering $r=0$ or $r=k$ does not yield any information, as $h_0 \in \mathrm{Sym}^k(H_i^\perp)$ and $h_k \in \mathrm{Sym}^k(H_i).$ For $1 \le r \le k-1$, this gives that $$\frac rk \|h_r\|^2=0 \iff h_r=0.$$ Therefore, $h \in \ker(\mathrm{Id}-T_kT_k^*)$ if and only if $h=h_0+h_k \in \mathrm{Sym}^k(H_i) \oplus \mathrm{Sym}^k(H_i^\perp)$ for any $1 \le i \le m.$
\end{proof}

\hfill \\
Thanks to \cref{Description_Kernel_GeneralCase}, we can decompose further to independent subspaces, \textit{i.e.} the non-zero intersections of the form $$\bigcap_{i \in A} H_i \, \cap \, \bigcap_{i \in A^c}H_i^\perp, \text{ where } A \subset \{1,\cdots,m\},$$ or sometimes also written as $$\bigcap_{i=1}^m H_i^{(\epsilon_i)} \text{ where } H_i^{(\epsilon_i)}=H_i^\perp \text{ if } \epsilon_i=0, \; H_i^{(\epsilon_i)}=H_i \text{ if } \epsilon_i=1.$$ Denote by $\{ K_k: 1 \le k \le k_0\}$ the enumeration of independent subspaces of $H$ with respect to the geometric Brascamp-Lieb datum $(\textbf{C},p).$ Recall that the orthogonal complement of $\bigoplus_{k=1}^{k_0}K_k$ is $K_{\mathrm{dep}}.$

\begin{lemma} \label{Inclusion_Lemma_Kernel}
     For any $k \ge 1$, one has: $$\ker(\mathrm{Id}-T_kT_k^*) \subset \mathrm{Sym}^k(K_j) \oplus \mathrm{Sym}^k(K_j^\perp),$$ for any $1 \le j \le k_0.$ 
\end{lemma}

\begin{remark}
    Before proving the lemma, let us show that $$\bigcap_{i=1}^m \mathrm{Sym}^k(U_i)=\mathrm{Sym}^k \Big(\bigcap_{i=1}^m U_i \Big),$$ for $U_i \subset \mathbb{R}^n$, $1 \le i \le m.$  The inclusion $\supset$ is clear, so only $\subset$ needs to be shown. For $x_2,\cdots,x_k \in \mathbb{R}^n$ and $z \in \bigcap_{i=1}^m \mathrm{Sym}^k(U_i)$, define $z[x_2,\cdots,x_k] \in \mathbb{R}^n$ as the unique vector such that $$\langle z[x_2,\cdots,x_k],x_1\rangle=\langle z,x_1 \otimes \cdots \otimes x_k \rangle \text{ for any } x_1 \in \mathbb{R}^n.$$ We claim first that $z[x_2,\cdots,x_k] \in \bigcap_{i=1}^m U_i$ for any $x_2,\cdots,x_k \in \mathbb{R}^n.$ Indeed, fixing $1 \le i \le m$ and as $z \in \mathrm{Sym}^k(U_i)$, one may write $$z=\sum_{\alpha}c_\alpha u_{\alpha,1}\otimes \cdots \otimes u_{\alpha,k}, \quad u_{\alpha,j} \in U_i,$$ so that $$z[x_2,\cdots,x_k]=\sum_{\alpha}c_\alpha \Big(\prod_{j=2}^k \langle x_j,u_{\alpha,j}\rangle\Big) u_{\alpha,1} \in U_i.$$ This proves the claim.
    Therefore, if one chooses an orthonormal basis $\{e_1,\cdots,e_r\}$ of $\bigcap_{i=1}^m U_i$ that we complete by an orthonormal basis $\{e_{r+1},\cdots,e_n\}$ of  $(\bigcap_{i=1}^m U_i)^\perp$, one writes $$z=\sum_{j_1,\cdots,j_k=1}^n z_{j_1\cdots j_k}e_{j_1}\otimes \cdots \otimes e_{j_k}.$$ Then, $$z[e_{j_2},\cdots,e_{j_k}]=\sum_{j_1=1}^n z_{j_1\cdots j_k}e_{j_1} \in \bigcap_{i=1}^m U_i,$$ so that $z_{j_1\cdots j_k}=0$ for any $j_1 \ge r+1.$ As the coefficients are invariant under permutations, we conclude that $z_{j_1\cdots j_k}=0$ for any $j_i \ge r+1$, so that $$z=\sum_{j_1,\cdots,j_k=1}^r z_{j_1\cdots j_k}e_{j_1}\otimes \cdots \otimes e_{j_k} \in \mathrm{Sym}^k \Big(\bigcap_{i=1}^m U_i\Big).$$ 
\end{remark}

\begin{proof}[Proof of \cref{Inclusion_Lemma_Kernel}]
Set $1 \le j \le k_0.$
    We may decompose $h \in \ker(\mathrm{Id}-T_kT_k^*)$, as $$h=\sum_{r=0}^k h_r, \text{ where } h_r \in \mathrm{Sym}^r(K_j) \odot \mathrm{Sym}^{k-r}(K_j^\perp).$$ Let us assume that there exists $1 \le r \le k-1$ with $h_r \neq 0.$ Then, there must exist $x_1,\cdots,x_r \in K_j$ such that the $k-r$ linear form on $K_j^\perp$ $$s:(y_1,\cdots,y_{k-r}) \mapsto h(x_1,\cdots,x_r,y_1,\cdots,y_{k-r})$$ is not identically zero. The point is now that $K_j=\bigcap_{i=1}^m H_i^{(\epsilon_i)}$ for $\epsilon=(\epsilon_1,\cdots,\epsilon_m) \in \{0,1\}^m$ and by \cref{Description_Kernel_GeneralCase}, $h \in \mathrm{Sym}^k(H_i^{(\epsilon_i)})\oplus \mathrm{Sym}^k(H_i^{(1-\epsilon_i)})$ so that by definition\footnote{This is since $x_1,\cdots,x_r \in K_j \subset H_i^a.$} of $s$, $s \in \mathrm{Sym}^{k-r}(K_j^\perp \cap H_i^{(\epsilon_i)})$ for any $1 \le i \le m.$ Consequently, by the previous remark, $$s \in \bigcap_{i=1}^m \mathrm{Sym}^{k-r}(K_j^\perp \cap H_i^{(\epsilon_i)})=\mathrm{Sym}^{k-r} \Big(K_j^\perp \cap \bigcap_{i=1}^m H_i^{(\epsilon_i)} \Big)=\mathrm{Sym}^{k-r}(K_j^\perp \cap K_j)=\{0\}.$$ This is obviously a contradiction which concludes the proof of the lemma. 
\end{proof}

\hfill \\
Trivially, $\mathrm{Sym}^k(K_j) \subset \mathrm{ker}(\mathrm{Id}-T_kT_k^*).$ 
Then, applying inductively \cref{Inclusion_Lemma_Kernel} gives that $$\mathrm{ker}(\mathrm{Id}-T_kT_k^*)=\Big( \bigoplus_{j=1}^{k_0} \mathrm{Sym}^k(K_j) \Big) \oplus \mathcal{K}_k,$$ where $\mathcal K_k$ is the part that is left of $\ker(\mathrm{Id}-T_kT_k^*)$ after removing $\bigoplus_{j=1}^{k_0}\mathrm{Sym}^k(K_j).$ 

\hfill \\
In fact, we show that as $k \ge 3$, the non-independent part disappears. 

\begin{lemma} \label{Full_Characterization_Kernels}
    The following holds: \begin{equation} \begin{aligned} \mathrm{ker}(\mathrm{Id}-T_kT_k^*)= \begin{cases}
        \mathbb{R}^n, & k=1 \\
        \{Y \in \mathrm{Sym}(n) : \substack{C_i^*C_iY(\mathrm{id}_{H}-C_i^*C_i)=0 \\ \text{ for any } 1 \le i \le m}\}, & k=2 \\
        \ds \bigoplus_{j=1}^{k_0} \mathrm{Sym}^k(K_j), & k \ge 3.
    \end{cases} \end{aligned} \end{equation}
\end{lemma}

\begin{proof}
    The cases $k=1,2$ were already addressed in \cref{Case1_2_AlreadyDone}. Therefore, we focus on the case $k \ge 3.$ 
    \begin{claim}
        Let $W$ be a finite dimensional Euclidean space and let $R_i: W \rightarrow W$ be orthogonal projections for any $1 \le i \le m.$ Suppose $h \in \mathrm{Sym}^k(W)$ satisfies that $$h \in \mathrm{Sym}^k(\ker(R_i)) \oplus \mathrm{Sym}^k(\mathrm{Ran}(R_i)), \text{ for any } 1 \le i \le m.$$ If $h \ne 0$, there exists $v \in W, v \ne 0$ such that for every $1 \le i \le m, R_iv \in \{0,v\}.$ 
    \end{claim}
    Applying this claim to $C_i^*C_i|_{K_{\mathrm{dep}}}$, as the assumption would hold for any $h \in \mathcal K_k$ by \cref{Description_Kernel_GeneralCase}, gives a contradiction if $\mathcal K_k \ne \{0\}$, as there would exist $v \in K_{\mathrm{dep}}, v \ne 0$, such that $v$ also belongs to an independent subspace, which is of course impossible. 
    
    \hfill \\
    Thus, it suffices to prove the claim. Define the scalar function $$q(x)=\langle h, x\otimes \cdots \otimes x\rangle \quad \forall x \in W,$$ that is $k$-homogeneous. Moreover, $h \ne 0$ implies that $q$ is a non-zero function. Set $$v=\arg \max_{x \in W,|x|=1}|q(x)|.$$ Using that $h \in \mathrm{Sym}^k(\ker(R_i)) \oplus \mathrm{Sym}^k(\mathrm{Ran}(R_i))$, $$|q(v)| \le |q(R_iv)|+|q((I_W-R_i)v)| \le |q(v)|\Big(|R_iv|^k+|(I_W-R_i)v|^k \Big).$$ As $1=|v|^2=|R_iv|^2+|(I_W-R_i)v|^2$, this means that for each $1 \le i \le m$, either $|(\mathrm{Id}_W-R_i)v|=0$ or $|R_iv|=0$, or in other words, $R_iv \in \{0,v\}.$ 
\end{proof}

\begin{remark}
    Suppose that $(\varphi_i)_{1 \le i \le m}$ is a centered tuple, satisfying equality in \cref{Gaussian_Variance_BL}. By \cref{Quadratic_Form_GeneralCase}, $$\varphi_i=\sum_{k=1}^\infty I_k^{(n_i)}(C_i^{\otimes k}h^{(k)}),$$ where $h^{(k)} \in \ker(\mathrm{Id}-T_kT_k^*).$ Consider $X \sim \gamma_n$ and define $$F=\sum_{k=1}^\infty I_k^{(n)}(h^{(k)}).$$ Then, by \cref{Conditional_Exp_Lemma}, $$\begin{array}{lll}
    \varphi_i(C_iX) &=& \ds \sum_{k=1}^\infty I_k^{(n_i)}(C_i^{\otimes k}h^{(k)}) \circ C_i(X)=\sum_{k=1}^\infty I_k^{(n)}((C_i^*C_i)^{\otimes k}h^{(k)})(X) \\ \\
    &=& \ds \sum_{k=1}^\infty \mathbb{E}\Big[ I_k^{(n)}(h^{(k)})(X)|C_i^*C_iX \Big]=\mathbb{E}[F(X)|C_i^*C_iX] \\ \\ 
    &=& \ds \mathbb{E}[F(X)|C_iX],
    \end{array}$$ \hfill \\ since $C_i^*C_iX=C_i^*(C_iX)$ and $C_iX=C_i(C_i^*C_iX).$  
    Moreover, $$\begin{array}{lll}
    \ds \sum_{i=1}^m p_i \varphi_i(C_iX) &=& \ds \sum_{i=1}^m p_i \sum_{k=1}^\infty I_k^{(n)}((C_i^*C_i)^{\otimes k}h^{(k)})(X) \\
    &=& \ds \sum_{k=1}^\infty I_k^{(n)}\Big(\sum_{i=1}^m p_i (C_i^*C_i)^{\otimes k} h^{(k)} \Big)(X)= F(X),
    \end{array}$$ as $$T_kT_k^*=\sum_{i=1}^m p_i (C_i^*C_i)^{\otimes k} \text{ and } T_kT_k^*h^{(k)}=h^{(k)} \text{ for any } k \ge 1.$$
    \hfill \\
    Finally, since $h^{(k)} \in \ker(\mathrm{Id}-T_kT_k^*)$, $h^{(k)}=(C_i^*C_i)^{\otimes k}h^{(k)}+(\mathrm{Id}-C_i^*C_i)^{\otimes k}h^{(k)}$ so that $$F(X)=\mathbb{E}[F(X)|C_i^*C_iX]+\mathbb{E}[F(X)|(\mathrm{Id}-C_i^*C_i)X].$$
\end{remark}

\begin{theorem} \label{Concrete_Description_Kernel}
    Equality in \cref{Quadratic_Form_GeneralCase} holds if and only if there exist constants $c_i \in \mathbb{R}$, arbitrary functions $g_k: K_k \rightarrow \mathbb{R}$ in $L^2(\gamma_{K_k}), 1 \le k \le k_0$ and $b \in K_{\mathrm{dep}}$, a self-adjoint matrix $A: K_{\mathrm{dep}} \rightarrow K_{\mathrm{dep}}$, whose eigenspaces are critical, such that $$\varphi_i(x_i)=c_i+\sum_{k : K_k \subset H_i}g_k(P_{K_k}x_i)+\langle AP_{K_{\mathrm{dep}}}x_i,P_{K_{\mathrm{dep}}}x_i\rangle+\langle P_{K_{\mathrm{dep}}}x_i,b\rangle.$$
\end{theorem}

\begin{remark}
    For $(\varphi_i)_{1 \le i \le m}$ of the form described in \cref{Concrete_Description_Kernel}, assume that $I_{K_{\mathrm{dep}}}+2A>0$ and that $0<\int_{K_k}e^{-g_k}\, d\gamma_{K_k}<+\infty$ for every independent subspace $K_k, 1 \le k \le k_0.$ Then, the functions $e^{-\varphi}$ form an extremizing tuple for the Gaussian formulation of the geometric Brascamp-Lieb inequality, \textit{i.e.} $$\int_{\mathbb{R}^n}\prod_{i=1}^m e^{-p_i\varphi_i(C_ix)}\, d\gamma_n(x)=\prod_{i=1}^m \Big(\int_{\mathbb{R}^{n_i}}e^{-\varphi_i}\, d\gamma_{n_i}\Big)^{p_i}.$$ For the classification of all extremizers of the Brascamp-Lieb inequality, including possible factors with zeroes, we refer to \cite{Vald}. 
\end{remark}

\begin{proof}
    We may first focus on establishing the form of the function $F$ that was constructed, \textit{i.e.} $$F=\sum_{k=1}^\infty I_k^{(n)}(h^{(k)}),$$ for $h^{(k)} \in \ker(\mathrm{Id}-T_kT_k^*).$ From \cref{Full_Characterization_Kernels}, we find that $$F(x)=\sum_{k=1}^{k_0}g_k(P_{K_k}x)+\langle b,x\rangle+\langle Ax,x\rangle-\mathrm{Tr}(A),$$ where $A: K_{\mathrm{dep}} \rightarrow K_{\mathrm{dep}}$ is an element of $\ker(\mathrm{Id}-T_2T_2^*)$ and $b \in K_{\mathrm{dep}},$ and $g_k \in L^2(\gamma_{K_k})$ arbitrary centered. Note that we subtracted here by $\mathrm{Tr}(A)$ since $F$ is centered and moreover, from what was proven above, $A \in \ker(\mathrm{Id}-T_2T_2^*)$ implies that its eigenspaces are critical (see \cref{QF_Degree2_Rk}).
    
    \hfill \\
    Let us now compute the conditional expectation, term by term, for $x_i \in H_i$. If $K_k \subset H_i$, conditioned on $C_iX=x_i$, $$P_{K_k}X=P_{K_k}C_i^*C_iX=P_{K_k}C_i^*x_i=P_{K_k}x_i,$$ and if $K_k \subset H_i^\perp$, $C_iX$ is independent of $g_k(P_{K_k}X).$ Therefore, \begin{equation} \begin{aligned} \mathbb{E} \Big[g_k(P_{K_k}X) \, \big| C_iX=x_i \Big]= \begin{cases}
        g_k(P_{K_k}x_i), & \text{ if } K_k \subset H_i, \\
        0 , & \text{ if } K_k \subset H_i^\perp.
    \end{cases}\end{aligned} \end{equation}
    Remember that $X=C_i^*C_iX+(\mathrm{Id}-C_i^*C_i)X$, with the second random vector being independent of $C_iX.$
    Thanks to this, we clearly have for the linear term that $$\mathbb{E}(X|C_iX=x_i)=x_i \Longrightarrow \mathbb{E}(\langle b,X\rangle|C_iX=x_i)=\langle b,x_i\rangle.$$ For the quadratic term, since there is no cross term by independence and that $A$ is block diagonal with respect to $H_i \oplus H_i^\perp$, $$\begin{array}{lll}
 \mathbb{E}(\langle AX,X\rangle -\mathrm{Tr}(A)|C_iX=x_i) &=& \ds \langle AC_i^*x_i,C_i^*x_i\rangle+\mathrm{Tr}(A|_{H_i^\perp})-\mathrm{Tr}(A) \\ \\
    &=& \ds \langle Ax_i,x_i\rangle-\mathrm{Tr}(A|_{H_i}).
    \end{array}$$
Combining these expressions, we get the desired form of extremizers. On the other hand, if $\varphi_i$ are of this form, it is easy to establish that $\varphi_i$ can be written in the form of \cref{Quadratic_Form_GeneralCase}. 
\end{proof}

Now that the kernel of the quadratic form is understood and in order to establish local quantitative stability near a Gaussian, we would like to prove a spectral gap on the orthogonal complement of \begin{equation} \label{Def/Thm_Kernel_QF} \begin{array}{lll}
\ds E=\ker(Q)=\left\{ (u_i)_{1 \le i \le m} : u_i=c_i+\sum_{k=1}^\infty I_k^{(n_i)}(C_i^{\otimes k}h^{(k)}) \text{ where } \substack{h^{(k)} \in \ker(\mathrm{Id}-T_kT_k^*) \\ \sum_{k \ge 1} k!\|h^{(k)}\|^2<+\infty} \right\} \subset \prod_{i=1}^m L^2(\gamma_{n_i})\\
\\
\ds = \left\{\Big(c_i+\sum_{k : K_k \subset H_i}g_k(P_{K_k}x_i)+\langle AP_{K_{\mathrm{dep}}}x_i,P_{K_{\mathrm{dep}}}x_i\rangle+\langle P_{K_{\mathrm{dep}}}x_i,b\rangle\Big)_{1 \le i \le m} : \substack{\ds b \in K_{\mathrm{dep}}, A: K_{\mathrm{dep}} \rightarrow K_{\mathrm{dep}} \\
\ds \text{ with critical eigenspaces}}  \right\}.
\end{array}\end{equation}

\hfill \\
Set $1 \le j<+\infty.$ If $v=(v_1,\cdots,v_m) \in E^\perp$, $$v_i=\sum_{k=0}^\infty I_k^{(n_i)}(v_i^{(k)}),$$ and considering $h^{(k)}=0$ if $k \ne j$, $h^{(j)} \in \ker(\mathrm{Id}-T_jT_j^*)$, $$\begin{array}{lll}
\ds \Big\langle v_i,\sum_{k=0}^\infty I_k^{(n_i)}(C_i^{\otimes k}h^{(k)})\Big \rangle_{L^2(\gamma_{n_i})} =\big\langle v_i, I_j^{(n_i)}(C_i^{\otimes j}h^{(j)})\big\rangle_{L^2(\gamma_{n_i})}=\big\langle I_j^{(n_i)}(v_i^{(j)}),I_j^{(n_i)}(C_i^{\otimes j}h^{(j)})\big\rangle_{L^2(\gamma_{n_i})}  \\ \\
\ds = j! \langle v_i^{(j)}, C_i^{\otimes j} h^{(j)} \rangle,
\end{array}$$ so that $$
0 = \sum_{i=1}^m p_i \Big\langle v_i,\sum_{k=0}^\infty I_k^{(n_i)}(C_i^{\otimes k}h^{(k)})\Big \rangle_{L^2(\gamma_{n_i})}=j!\sum_{i=1}^m p_i \langle v_i^{(j)},C_i^{\otimes j}h^{(j)}\rangle,
$$ while $v_i^{(0)}=0$ for any $1 \le i \le m.$
\hfill \\
Consequently, $v^{(j)}=(v_1^{(j)},\cdots,v_m^{(j)})$ is orthogonal to $T_j^*\ker(\mathrm{Id}-T_jT_j^*)$, which implies that $$\|T_jv^{(j)}\|^2 \le (1-\lambda_j)\|v^{(j)}\|_p^2,$$ where $$\begin{array}{lll}
1-\lambda_j &=& \ds \sup_{v \in [T_j^*\ker(\mathrm{Id}-T_jT_j^*)]^\perp \cap \{\|\cdot\|_p=1\}} \|T_jv\|^2=\|T_j^*T_j|_{[T_j^*\ker(\mathrm{Id}-T_jT_j^*)]^\perp}\|_{\mathrm{op}} \\ 
\nm
&=& \ds \|T_jT_j^*|_{[\ker(\mathrm{Id}-T_jT_j^*)]^\perp}\|_{\mathrm{op}}<1.
\end{array}$$ Here, we use the convention $\|0|_{\{0\}}\|_{\mathrm{op}}:=0.$ Then, $$\begin{array}{lll}
\ds \sum_{j=1}^m p_j \mathrm{Var}_{\gamma_{n_j}}(v_j)-\mathrm{Var}_{\gamma_n}\Big(\sum_{i=1}^m p_i v_i \circ C_i\Big)=\sum_{j=1}^\infty j! \Big( \|v^{(j)}\|_p^2-\|T_jv^{(j)}\|^2 \Big)
\ds \ge  \sum_{j=1}^\infty j! \lambda_j \|v^{(j)}\|_p^2. 
\end{array}$$
Now, $$\sum_{j=1}^\infty j! \lambda_j \|v^{(j)}\|_p^2 =\sum_{j=1}^\infty j! \lambda_j \Big(\sum_{i=1}^m p_i |v_i^{(j)}|^2 \Big)=\sum_{i=1}^m p_i \Big(\sum_{j=0}^\infty j! \lambda_j |v_i^{(j)}|^2 \Big),$$ and if there would be a uniform lower bound of $\lambda_j$, \textit{i.e.} $c>0$ such that $\lambda_j \ge c$ for any $j \ge 1$, we would get:

\begin{theorem} \label{General_SpGap_Quadratic}
For $v \in E^\perp$, $$\sum_{j=1}^m p_j \mathrm{Var}_{\gamma_{n_j}}(v_j)-\mathrm{Var}_{\gamma_n}\Big(\sum_{i=1}^m p_i v_i \circ C_i\Big) \ge c\sum_{i=1}^m p_i \|v_i\|_{L^2(\gamma_{n_i})}^2.$$ 
\end{theorem}

\hfill \\
To prove this theorem, it suffices to show, by the reason explained previously, the following.

\begin{theorem}
    For every $j \ge 1$, $$\lambda_j \ge \min\{ \lambda_2,\lambda_3\}.$$
\end{theorem}

\begin{remark}
    The spectral gap theorem is coherent with \cref{Full_Characterization_Kernels}, where the structure for $j \ge 4$ is ``inherited'' from the case $j=3.$ Moreover, one has trivially that $\lambda_1=1$, since $T_1T_1^*=\mathrm{id}_H$, so that every spectral gap is controlled from below by the ones at degree $2$ and $3.$
\end{remark}

\begin{proof}
For $j=1,2,3$, there is nothing to show so that we assume $j \ge 4.$
    Recall that $$T_jT_j^*=\sum_{i=1}^m p_i (C_i^*C_i)^{\otimes j},$$ and $T_j^*\ker(\mathrm{Id}-T_jT_j^*)$ is the space that preserves the norm of $T_j.$
    First, note that $C_i^*C_i$ preserves independent subspaces $K_k$ for any $1 \le k \le k_0$, since either $K_k \subset H_i$ so that $C_i^*C_iK_k \subset K_k$ or $K_k \subset H_i^\perp$ so that $C_i^*C_iK_k=\{0\}.$
    \hfill \\
    Write $\epsilon^{(k)}=(\epsilon_1^{(k)},\cdots,\epsilon_m^{(k)}) \in \{0,1\}^m$ and $A_k=\{ i: \epsilon_i^{(k)}=1\}$ for any $1 \le k \le k_0$, defined such that $$K_k=\bigcap_{i=1}^m H_i^{(\epsilon_i^{(k)})}, \text{ where } H_i^{(\epsilon_i^{(k)})}=H_i^\perp \text{ if } \epsilon_i^{(k)}=0, H_i^{(\epsilon_i^{(k)})}=H_i \text{ if } \epsilon_i^{(k)}=1.$$ In particular, $$K_k=\bigcap_{i \in A_k}H_i \cap \bigcap_{i \in A_k^c}H_i^\perp,$$ and $$C_i^*C_i=\bigoplus_{k=1}^{k_0} \epsilon_i^{(k)} \mathrm{Id}_{K_k} \oplus C_i^*C_i|_{K_{\mathrm{dep}}}.$$ For $x \in K_k$, $$\Big(\sum_{i \in A_k}p_i\Big)x =\sum_{i=1}^m p_i C_i^*C_ix=x \Longrightarrow \sum_{i \in A_k}p_i=1.$$ Finally, for any integers $r$ and $a:=(a_1,\cdots,a_{k_0}) \in \mathbb{N}^{k_0}$ such that $$r+\sum_{k=1}^{k_0}a_k=j,$$ define $$\mathcal H_{r,a}=\mathrm{Sym}^r(K_{\mathrm{dep}}) \odot \bigodot_{k=1}^{k_0} \mathrm{Sym}^{a_k}(K_k), \quad S_a=\bigcap_{k \text{ s.t. } a_k>0} A_k.$$ Then, $$\mathrm{Sym}^j(H)=\bigoplus_{r+\sum_{k=1}^{k_0}a_k=j} \mathcal H_{r,a},$$ and $(C_i^*C_i)^{\otimes j}$ on $\mathcal H_{r,a}$ is identically $0$ if $i \not \in S_a$, else it acts like the identity on the $\bigodot_{k=1}^{k_0}\mathrm{Sym}^{a_k}(K_k)$-part and like $(C_i^*C_i|_{K_{\mathrm{dep}}})^{\otimes r}$ on the $\mathrm{Sym}^r(K_{\mathrm{dep}})$-part. Indeed, if $i \notin S_a$, there exists $1 \le k \le k_0$ such that $a_k>0$ and $\epsilon_i^{(k)}=0 $ which implies $ K_k \subset H_i^\perp$ and as previously observed, $C_i^*C_iK_k=\{0\}$ in that case. 
    
\hfill \\
Hence, $T_jT_j^*$ acts on $\mathcal H_{r,a}$ as $$\sum_{i \in S_a}p_i (C_i^*C_i|_{K_{\mathrm{dep}}})^{\otimes r}.$$ First, assume that there exists $1 \le j_1, j_2 \le k_0, j_1 \ne j_2$ so that $a_{j_1},a_{j_2}>0.$ Then, $S_a \subset A_{j_1} \cap A_{j_2}$ and $$\sum_{i \in S_a}p_i(C_i^*C_i|_{K_{\mathrm{dep}}})^{\otimes r} \le \Big(\sum_{i \in A_{j_1} \cap A_{j_2}} p_i \Big) \mathrm{Id}.$$ Now, as $A_{j_1} \ne A_{j_2}$, $$\sum_{i \in A_{j_1} \cap A_{j_2}}p_i<\sum_{i \in A_{j_1}}p_i=1,$$ and notice also that on $K_{j_1} \odot K_{j_2}$, $T_2T_2^*$ acts like $\Big(\ds \sum_{i \in A_{j_1} \cap A_{j_2}}p_i\Big)\mathrm{Id}.$ Thus, $$\|T_jT_j^*|_{\mathcal H_{r,a}}\| \le \|T_2T_2^*|_{[\ker(\mathrm{Id}-T_2T_2^*)]^\perp}\|.$$
Else, assume there exists a unique $1 \le k \le k_0$ such that $a_k>0$ and $r \ge 1.$ Then, $S_a=A_k$ and $T_jT_j^*$ acts on $\mathcal H_{r,a}$, as $$\sum_{i \in A_k}p_i(C_i^*C_i|_{K_{\mathrm{dep}}})^{\otimes r} \le \Big(\sum_{i \in A_k}p_iC_i^*C_i|_{K_{\mathrm{dep}}}\Big) \otimes \mathrm{Id}^{\otimes(r-1)},$$ which implies in particular that $$\Big\|\sum_{i \in A_k}p_i(C_i^*C_i|_{K_{\mathrm{dep}}})^{\otimes r}\Big\| \le \Big\|\sum_{i \in A_k}p_iC_i^*C_i|_{K_{\mathrm{dep}}}\Big\|,$$ where the second operator is the action of $T_2T_2^*$ on $K_k \odot K_{\mathrm{dep}}.$ 
Hence, we conclude that $$\|T_jT_j^*|_{\mathcal H_{r,a}}\| \le \|T_2T_2^*|_{[\ker(\mathrm{Id}-T_2T_2^*)]^\perp}\|.$$
Finally, assume that $a_k=0$ for all $1 \le k \le k_0.$ Then, $\mathcal H_{r,a}=\mathrm{Sym}^j(K_{\mathrm{dep}})$, so that $T_jT_j^*$ is $$\sum_{i=1}^m p_i(C_i^*C_i|_{K_{\mathrm{dep}}})^{\otimes j} \le \Big[\sum_{i=1}^m p_i(C_i^*C_i|_{K_{\mathrm{dep}}})^{\otimes 3} \Big] \otimes \mathrm{Id}^{\otimes (j-3)},$$ so that $$\|T_jT_j^*|_{\mathrm{Sym}^j(K_{\mathrm{dep}})}\| \le \|T_3T_3^*|_{\mathrm{Sym}^3(K_{\mathrm{dep}})}\|.$$ 
\hfill \\
Since these are the only three cases encountered when working on $[\ker(\mathrm{Id}-T_jT_j^*)]^\perp$; the only other possibility being $\mathrm{Sym}^j(K_k)$, which is precisely ruled out by the fact that we are considering $[\ker(\mathrm{Id}-T_jT_j^*)]^\perp$ for $j \ge 3$, this allows to show that $$1-\lambda_j \le \max\{1-\lambda_2,1-\lambda_3\} \Longrightarrow \lambda_j \ge \min\{ \lambda_2,\lambda_3\},$$ as wished. 
\end{proof}

Establishing the spectral gap is important in the proof of our main result but it is only useful to treat functions that would be small with respect to the $L^\infty$ norm. Hence, the main difficulty in the next section will be to deal, at the same time, with the regions where some functions are quite ``large''.

\section{Proof of \cref{Complete_Gaussian_stability}}
\label{Sec:4}

In this section, we may assume that $(\textbf{C},p)$ is a geometric Brascamp-Lieb datum. A useful notion of distance that we will employ is the Hellinger distance between two probability measures, which is defined as follows: for $P(x)=p(x)\, dx, Q(y)=q(y)\, dy$, then $$H_{\text{Hellinger}}^2(P,Q)=\frac12 \int \Big(\sqrt{p(x)}-\sqrt{q(x)} \Big )^2\, dx=1-\int \sqrt{p(x)q(x)}\, dx.$$ 

\hfill \\
Now, define the logarithmic deficit $$\mathrm{Def}(f_1,\cdots,f_m)=\sum_{i=1}^m p_i \log \Big(\int_{\mathbb{R}^{n_i}}f_i(x_i)\, d\gamma_{n_i}(x_i) \Big)-\log \Big(\int_{\mathbb{R}^n} \prod_{i=1}^m f_i(C_ix)^{p_i}\, d\gamma_n(x) \Big) \ge 0$$ and for $f_i \ge 0$ with $\int f_i\, d\gamma_{n_i}=1$ for any $1 \le i \le m$, $$A\big[(f_i)_{1 \le i \le m}\big]=\inf_{(\eta_i)_{1 \le i \le m} \text{ normalized extremizer}} \sum_{i=1}^m p_i H_{\text{Hellinger}}^2(f_i\, d\gamma_{n_i},\eta_i).$$ 

\hfill \\
Let us slightly reformulate \cref{Complete_Gaussian_stability} in the geometric setting. 

\begin{theorem} \label{Hellinger_Stab}
Let $(\textbf{C},p)$ be a geometric Brascamp-Lieb datum. There exists a constant $c=c(\textbf{C},p)>0$ such that, for any $(f_i)_{1 \le i \le m}$ with $f_i \ge 0$ and $\int f_i\, d\gamma_{n_i}=1$ for any $1 \le i \le m$, $$1-\int_{\mathbb{R}^n}\prod_{i=1}^m f_i(C_ix)^{p_i}\, d\gamma_n \ge cA\big[(f_i)_{1 \le i \le m}\big].$$
\end{theorem}

\begin{proof}[Proof of \cref{Hellinger_Stab}]

\hfill \\
\textbf{Step 1 (Local to global):}
We claim that it suffices to prove the local version of the inequality, \textit{i.e.} there exist $c,\epsilon_0>0$ such that $$A \big[ (f_i)_{1 \le i \le m}\big] \le \epsilon_0^2$$ implies that $$1-\int_{\mathbb{R}^n}\prod_{i=1}^m f_i(C_ix)^{p_i}\, d\gamma_n \ge cA\big[(f_i)_{1 \le i \le m}\big].$$

\hfill \\
Indeed, considering the Ornstein-Uhlenbeck flow $(P_t)_{t \ge 0}$, it is well-known (see \cite{Barthe_Cordero}) that $$t \mapsto \mathrm{Def}(f_{1,t},\cdots,f_{m,t}) \text{ decreases to } 0,$$ where we denote $f_{i,t}=P_t^{(n_i)}f_i$ for any $t \ge 0.$ Moreover, $t \mapsto A\big[(f_{i,t})_{1 \le i \le m} \big]$ is continuous and converges to $0$ as $t \rightarrow +\infty$, as a simple consequence of the triangle inequality, $(\sqrt a -\sqrt b)^2 \le |a-b|$ for any $a,b \ge 0$ and the fact that the Ornstein-Uhlenbeck semigroup is strongly continuous on $L^1(\gamma_n)$.

\hfill \\
Consequently, assuming the local quantitative stability, if $$A\big[(f_i)_{1 \le i \le m} \big]>\epsilon_0^2,$$ then there exists $T>0$ such that $A\big[(f_{i,T})_{1 \le i \le m} \big]=\epsilon_0^2$ by continuity, so that $$1-e^{-\mathrm{Def}(f_1,\cdots,f_m)} \ge 1-e^{-\mathrm{Def}(f_{1,T},\cdots,f_{m,T})} \ge c \epsilon_0^2,$$ by local stability. Hence, since $A\big[(f_i)_{1 \le i \le m}\big] \le \sum_{i=1}^m p_i$, $$1-\int_{\mathbb{R}^n}\prod_{i=1}^m f_i(C_ix)^{p_i}\, d\gamma_n \ge c \min\Big \{\frac{\epsilon_0^2}{\sum_{i=1}^m p_i},1\Big\} A\big[ (f_i)_{1 \le i \le m}\big].$$

\hfill \\
\textbf{Step 2 (Local perturbation around the standard Gaussian):}
In this step, we will show that we can focus on the local case around the standard Gaussian (\textit{i.e.} in Hellinger distance, the nearest extremizer to $f_i\, d\gamma_{n_i}$ is $d\gamma_{n_i}$) and that we may assume orthogonality. 

\hfill \\
To assume almost everywhere positivity of the inputs $(f_i)_{1 \le i \le m}$, we may replace $f_i$ by its heat semigroup $P_t^{(n_i)}f_i$ for $t>0$ as in Step 1 and eventually let $t \downarrow 0.$ 

\hfill \\
\textit{Step 2.1(Change of variables):}
Recall that $$H=\bigoplus_{k=1}^{k_0}K_k \oplus K_{\mathrm{dep}},$$ so that for any $x \in H$, there exists a unique representation of the form: $$x=\sum_{k=1}^{k_0}x_k+\tilde x, \text{ where } x_k \in K_k \text{ for any } 1 \le k \le k_0, \tilde x \in K_{\mathrm{dep}}.$$ If we consider $(g_i^*)_{1 \le i \le m}$ to be a normalized extremizer, set $$d\mu^*=\prod_{i=1}^m g_i^*(C_ix)^{p_i}\, d\gamma_n(x),$$ and by the characterization of extremizers of \cite{Vald}, we have $$\mu^*=\bigotimes_{k=1}^{k_0} \nu_k \otimes \mathcal N(b,\Sigma),$$ where  $\nu_k$ are probability measures on $K_k$, that are absolutely continuous with respect to $\gamma_{K_k}$, $b \in K_{\mathrm{dep}}$ and $\Sigma \in \mathrm{Sym}_{>0}(\dim(K_{\mathrm{dep}}))$ commutes with $C_i^*C_i|_{K_{\mathrm{dep}}}$ for any $1 \le i \le m.$ For $0<\tau<1$, set $$\nu_{k,\tau}=(1-\tau)\nu_k+\tau \gamma_{K_k},$$ so that $$\frac{d\nu_{k,\tau}}{d\gamma_{K_k}}=(1-\tau)\frac{d\nu_k}{d\gamma_{K_k}}+\tau \ge \tau>0.$$ Notice that $$\Big\|\frac{d\nu_{k,\tau}}{d\gamma_{K_k}}-\frac{d\nu_k}{d\gamma_{K_k}}\Big\|_{L^1(\gamma_{K_k})}=\tau\Big\|1-\frac{d\nu_k}{d\gamma_{K_k}}\Big\|_{L^1(\gamma_{K_k})} \le 2\tau,$$ which implies $$\Big\|\sqrt{\frac{d\nu_{k,\tau}}{d\gamma_{K_k}}}-\sqrt{\frac{d\nu_k}{d\gamma_{K_k}}}\Big\|_{L^2(\gamma_{K_k})}^2 \le 2\tau.$$ Therefore, by considering instead $\mu^*_\tau=\bigotimes_{k=1}^{k_0}\nu_{k,\tau} \otimes \mathcal N(b,\Sigma)$, the extremizer may be chosen positive, as the common independent factors may be chosen positive, and arbitrarily close to the distance infimum $A[(f_i)_{1 \le i \le m}].$ We thus assume that $g_i^*>0$ for any $1 \le i \le m.$ 

\hfill \\
Define $R_k: K_k \rightarrow K_k$ invertible a.e. such that $(R_k)_{\#}\gamma_{K_k}=\nu_k$ for any $1 \le k \le k_0$ and $$S\Big(\sum_{k=1}^{k_0}x_k+\tilde x\Big)=\sum_{k=1}^{k_0}R_k(x_k)+\Sigma^{1/2}\tilde x+b.$$ Then, by definition, $$S_{\#}\gamma_n=\mu^*.$$

\hfill \\
Since $H_i=\oplus_{k:K_k \subset H_i}K_k \oplus C_i^*C_iK_{\mathrm{dep}}$, we may define $$\tilde{S}_i\Big(\sum_{k:K_k \subset H_i}x_k+\tilde x_i\Big)=\sum_{k:K_k\subset H_i}R_k(x_k)+\Sigma^{1/2}\tilde x_i+C_i^*C_ib,$$ for $x_k \in K_k, \tilde x_i \in C_i^*C_iK_{\mathrm{dep}}$ for any $1 \le i \le m.$ Then, if $x=\sum_{k=1}^{k_0}x_k+\tilde x$, $$C_i^*C_ix=\sum_{k:K_k\subset H_i}x_k+C_i^*C_i\tilde x,$$ which implies that $$\begin{array}{lll}
\ds \tilde S_i(C_i^*C_ix) &=& \ds \sum_{k:K_k \subset H_i}R_k(x_k)+\Sigma^{1/2}C_i^*C_i\tilde x+C_i^*C_ib=C_i^*C_i \Big(\sum_{k=1}^{k_0}R_k(x_k)+\Sigma^{1/2}\tilde x+b\Big) \\ \\ 
&=& \ds C_i^*C_iS(x),
\end{array}$$ using also that $\Sigma^{1/2}$ commutes with $C_i^*C_i|_{K_\mathrm{dep}}.$ Setting $S_i=C_i\tilde{S}_iC_i^*$, this yields $$S_iC_i=C_iS \quad \forall 1 \le i \le m,$$ and in particular, $${S_i}_{\#}\gamma_{n_i}=(S_iC_i)_{\#}\gamma_n=(C_i)_{\#}\mu^*=g_i^* \, d\gamma_{n_i}.$$ Introducing $\tilde f_i=\frac{f_i}{g_i^*}\circ S_i,$ we have $$\int \tilde f_i \, d\gamma_{n_i}=\int f_i\, d\gamma_{n_i}=1,$$ and $$\int \prod_{i=1}^m \tilde f_i(C_ix)^{p_i}\, d\gamma_n=\int \prod_{i=1}^m \bigg(\frac{f_i(C_ix)}{g_i^*(C_ix)} \bigg)^{p_i}\, d\mu^*=\int \prod_{i=1}^m f_i(C_ix)^{p_i}\, d\gamma_n,$$ where we used that $S_iC_i=C_iS$ for any $1 \le i \le m.$  By the same argument, $f \mapsto (f/\sqrt{g_i^*})\circ S_i$ is unitary, so that the Hellinger distance between $(f_i\, d\gamma_{n_i})_{1 \le i \le m}$ and the normalized extremizers coincide with the one between $(\tilde f_i d\gamma_{n_i})_{1 \le i \le m}$ and the normalized extremizers, for any $1 \le i \le m.$ 

\hfill \\
\textit{Step 2.2(Concentration-compactness):}
Let us check that the infimum in the definition of $A$ is attained. In fact, it suffices to use the special structure of the extremizers and weak compactness.

\hfill \\
Assume that $(f_i)_{1 \le i \le m}$ is quite close to the standard Gaussian (up to the previous sub-step) in Hellinger distance and let $g^{(l)}$ be a minimizing sequence, \textit{i.e.} a tuple of square roots of densities. Then, $$\|g^{(l)}-(1,\cdots,1)\|_{2,p} \le \|g^{(l)}-(\sqrt{f_i})_{1 \le i \le m}\|_{2,p}+\|(\sqrt{f_i})_{1 \le i \le m}-(1,\cdots,1)\|_{2,p},$$ is uniformly small.
Here, $$\|u\|_{2,p}^2:=\sum_{i=1}^m p_i\|u_i\|_{L^2(\gamma_{n_i})}^2.$$
By definition, $$g_i^{(l)}=\prod_{k: K_k \subset H_i}\tilde g_k^{(l)} \cdot \sqrt{g_{B_l,b_l}},$$ where $\|\tilde g_k^{(l)}\|_{L^2(\gamma_{K_k})}=1$, with $b_l \in K_{\mathrm{dep}}, B_l \in \mathrm{Sym}_{>0}(K_{\mathrm{dep}})$ such that $C_i^*C_iB_l=B_lC_i^*C_i$ for any $1 \le i \le m.$  The mixed part of $\|g^{(l)}-(1,\cdots,1)\|_{2,p}^2$ is $\langle g^{(l)},(1,\cdots,1)\rangle_{2,p}$ and every term of this sum is of the form $$p_i \langle g_i^{(l)},1\rangle_{L^2(\gamma_{n_i})}=p_i\prod_{k: K_k \subset H_i}\langle \tilde g_k^{(l)},1\rangle_{L^2(\gamma_{K_k})} \cdot \langle \sqrt{g_{B_l,b_l}},1\rangle_{L^2(\gamma_{C_i^*C_iK_{\mathrm{dep}}})},$$ which must be bounded below, so that every factor is uniformly bounded from below. 

\hfill \\
It is also easy to compute that $$-\log\Big(1-H^2_{\text{Hellinger}}(\gamma_n, \mathcal N(m,A))\Big)=\frac12 \log \det \Big(\frac{I_n+A}{2}\Big)-\frac14 \log \det(A)+\frac14 m^T(I_n+A)^{-1}m,$$ and note that the covariance part yields the following contribution: $$\frac12 \sum_{\lambda \in \sigma(A)}\log \cosh\Big(\frac12\log \lambda \Big).$$ As the negative logarithm must be bounded above, the covariance eigenvalues are therefore bounded above and away from zero and the mean must be bounded and therefore need to converge up to subsequence, say $B$ and $b$. This implies that the Gaussian part of $g_i^{(l)}$ converges strongly in $L^2$. Up to another subsequence, $\tilde g_k^{(l)} \rightharpoonup \tilde g_k \text{ in } L^2(\gamma_{K_k}),$ with $\|\tilde g_k\|_{L^2(\gamma_{K_k})} \le 1.$ Since the factors live on different spaces, $$g_i^{(l)} \rightharpoonup \prod_{k:K_k \subset H_i}\tilde g_k \cdot \sqrt{g_{B,b}} \text{ in } L^2,$$ as the weak convergence can be tested on products of functions of separate coordinates and sums of such functions are dense in $L^2(\gamma_n).$ 
Since minimizing the Hellinger distance is equivalent to maximizing $$\sum_{i=1}^m p_i \langle \sqrt{f_i},g_i\rangle_{L^2},$$ and $\tilde g_k$ may not be normalized, we normalize them, which increases the sum and gives a desired nearest extremizer. 

\hfill \\
Let us in fact prove that every component of the nearest extremizer $(g_i^*)_{1 \le i \le m}$ constructed is positive almost everywhere. Using the notation of Step 2.1, we fix an independent subspace $K_k$ and $$u_k=\sqrt{\frac{d\nu_k}{d\gamma_{K_k}}} \ge 0, \quad \|u_k\|_{L^2(\gamma_{K_k})}=1.$$ For $1 \le i \le m$ such that $K_k \subset H_i$, $H_i=K_k \oplus (H_i \cap K_k^\perp)$ and so, $$\sqrt{g_i^*(x_k+y)}=u_k(x_k)v_{i,k}(y),$$ where $v_{i,k}$ is the product of all other square-root factors with $v_{i,k} \ge 0$ and $\|v_{i,k}\|_{L^2(\gamma_{H_i \cap K_k^\perp})}=1.$ Then, $$\sum_{i=1}^m p_i \int_{\mathbb{R}^{n_i}}\sqrt{f_ig_i^*}\, d\gamma_{n_i}=\langle u_k,F_k\rangle_{L^2(\gamma_{K_k})}+\sum_{i: K_k \not \subset H_i}p_i \int_{H_i}\sqrt{f_ig_i^*}\, d\gamma_{H_i},$$ where $$F_k(x_k)=\sum_{i: K_k \subset H_i}p_i\int_{H_i \cap K_k^\perp}\sqrt{f_i(x_k+y)}v_{i,k}(y)\, d\gamma_{H_i \cap K_k^\perp}(y).$$ Then, $F_k \in L^2(\gamma_{K_k})$ since by Cauchy-Schwarz, $$\Big(\int_{H_i \cap K_k^\perp}\sqrt{f_i(x_k+y)}v_{i,k}(y)\, d\gamma_{H_i \cap K_k^\perp}(y)\Big)^2 \le \int_{H_i \cap K_k^\perp} f_i(x_k+y)\, d\gamma_{H_i \cap K_k^\perp}(y).$$ 
\hfill \\
Moreover, $$F_k>0 \quad \gamma_{K_k}\text{-almost everywhere},$$ as $f_i>0$ almost everywhere, while $v_{i,k}$ is non-negative and non-zero, which by Fubini makes the corresponding integral strictly positive for almost every $x_k.$ 

\hfill \\
The independent factor being arbitrary and $(g_i^*)_{1 \le i \le m}$ being a nearest extremizer, we have that $$u_k=\arg \max_{\|v_k\|_{L^2(\gamma_{K_k})}=1} \langle v_k,F_k\rangle_{L^2(\gamma_{K_k})}=\frac{F_k}{\|F_k\|_{L^2(\gamma_{K_k})}}>0,$$ $\gamma_{K_k}$-almost everywhere. Repeating this argument for each independent factor proves that $g_i^*$ is positive almost everywhere, as the dependent Gaussian factor is clearly already positive. 

\hfill \\
In conclusion, for $(f_i\, d\gamma_{n_i})_{1 \le i \le m}$, we find $(g_i^*)_{1 \le i \le m}$ to be the nearest extremizer and denote again $\tilde f_i=\frac{f_i}{g_i^*}\circ S_i,$ where $S_i$ is defined as in Step 2.1. Then, $(d\gamma_{n_i})_{1 \le i \le m}$ is the nearest extremizer to $(\tilde f_i\, d\gamma_{n_i})_{1 \le i \le m}$, so that by a standard first variation argument,  $$(\sqrt{\tilde f_i}-1)_{1 \le i \le m}  \text{ is orthogonal to the centered part of the kernel of the quadratic form},$$ and $$\sum_{i=1}^m p_i\|\sqrt{\tilde f_i}-1\|_{L^2(\gamma_{n_i})}^2=2A\big[(f_i)_{1 \le i \le m} \big].$$ 

\hfill \\
\textbf{Step 3 (Separation in small and large values):}
In this final step, we prove the local stability around the standard Gaussian, by decomposing the deficit in the quadratic part, where the spectral gap will be employed, and the remainder part, which must be shown to be negligible compared to $A\big[(f_i)_{1 \le i \le m} \big].$ 

\hfill \\
\textit{Step 3.1 (The remainder and the quadratic part):}
Denote $\epsilon=\sqrt{2A\big[(f_i)_{1 \le i \le m}\big]} \in [0;\sqrt 2\epsilon_0].$ As proven previously, we may assume that the nearest Gaussian extremizer to $(f_i\, d\gamma_{n_i})_{1 \le i \le m}$ is $(d\gamma_{n_i})_{1 \le i \le m}.$ Denote $$h_i=\sqrt{f_i}-1, \quad \int (1+h_i)^2\, d\gamma_{n_i}=1 \Longrightarrow \int h_i\, d\gamma_{n_i}=-\frac12 \|h_i\|_{L^2(\gamma_{n_i})}^2,$$ and $$\epsilon^2=\sum_{i=1}^m p_i \|h_i\|_{L^2(\gamma_{n_i})}^2.$$ Moreover, $(h_i)_{1 \le i \le m}$ is orthogonal to the centered part of the kernel of the quadratic form, composed of $$\Big(\sum_{k : K_k \subset H_i}g_k(P_{K_k}x_i)+\langle AP_{K_{\mathrm{dep}}}x_i,P_{K_{\mathrm{dep}}}x_i\rangle+\langle P_{K_{\mathrm{dep}}}x_i,b\rangle-\mathrm{Tr}(A|_{H_i})\Big)_{1 \le i \le m},$$ 
    where $A: K_{\mathrm{dep}} \rightarrow K_{\mathrm{dep}}$ is an element of $\ker(\mathrm{Id}-T_2T_2^*), b \in K_{\mathrm{dep}}$ and $g_k \in L^2(\gamma_{K_k})$ are arbitrarily centered functions.
    Define $$\tilde Q((h_i)_{1 \le i \le m})=\sum_{i=1}^m p_i \|h_i\|_{L^2(\gamma_{n_i})}^2-\|Th\|_{L^2(\gamma_n)}^2,$$ where we recall that $Th(x)=\sum_{i=1}^m p_ih_i(C_ix).$ Then, $$R(x)=\Big(1+\sum_{i=1}^m p_ix_i\Big)^2-\prod_{i=1}^m (1+x_i)^{2p_i},$$ for any $x=(x_1,\cdots,x_m) \in [-1;+\infty)^m.$ Notice that $$\int \Big(1+\sum_{i=1}^m p_ih_i(C_ix)\Big)^2 \, d\gamma_n(x)=1+2\sum_{i=1}^m p_i \int h_i \, d\gamma_{n_i}+\|Th\|_{L^2(\gamma_n)}^2=1-\tilde Q(h),$$ so that $$1-\int \prod_{i=1}^m f_i(C_ix)^{p_i}\, d\gamma_n(x)=\int_{\mathbb{R}^n}R(h_1(C_1x),\cdots,h_m(C_mx))\, d\gamma_n(x)+\tilde Q(h),$$ which we claim to be larger than $\tilde Q(h)-o(\epsilon^2).$ 
    Before establishing the claim, we record the Taylor expansion of $R$, which is given by $$R(x)=\sum_{i=1}^m p_ix_i^2-\Big(\sum_{i=1}^m p_ix_i\Big)^2+O \Big(a \sum_{i=1}^m p_ix_i^2 \Big), \text{ if } \max_{1 \le i \le m}|x_i| \le a,$$ since $(1+x_i)^{2p_i}=1+2p_ix_i+p_i(2p_i-1)x_i^2+O(|x_i|^3)$, and define the quadratic expansion as $q(x_1,\cdots,x_m).$ 

\hfill \\
\textit{Step 3.2 (The remainder is negligible):}
For $0<\eta<a<\frac12$, define $$M_a=\{ x: \max_{1 \le i \le m}|h_i(C_ix)|>a\}, \quad S(x)=\{1 \le i \le m: |h_i(C_ix)|>\eta\}.$$ Clearly, if $x \in M_a$, $S(x) \ne \emptyset.$

\hfill \\
Assume that $S=S(x)$ and $$\bigcap_{i \in S}C_i^*H_i=\{0\}.$$
The following lemma is standard (see \cite{Chen_Dafnis_Paouris}).

\begin{lemma} \label{Holder_AM_GM}
    For $S \subset \{1,\cdots,m\}$ and $s_i>0$ for any $i \in S$, such that $\sum_{i \in S} s_iC_i^*C_i \le \mathrm{Id}_H$, then $$\int \prod_{i \in S}|w_i(C_ix)| \, d\gamma_n(x) \le \prod_{i \in S}\|w_i\|_{L^{1/s_i}(\gamma_{n_i})}.$$ 
\end{lemma}

\hfill \\
Now, let us show that there exists $s_i \in [p_i;1] \text{ for } i \in S$ such that $$\sum_{i \in S}s_iC_i^*C_i \le \mathrm{id}_H \text{ and } \sum_{i \in S}s_i>1.$$
\hfill \\
Indeed, set $$\lambda_S=\Big\|\sum_{i \in S}p_iC_i^*C_i \Big\|_{\mathrm{op}}.$$ Then, $p_i \le \lambda_S \le 1$ and $\lambda_S<\sum_{i \in S}p_i.$ Indeed, $ \lambda_S \le \sum_{i \in S}p_i$ is clear and equality would give that there exists a non-zero unit vector $x$, such that $$\sum_{i \in S}p_i |C_ix|^2 =\sum_{i \in S}p_i \Longrightarrow |C_ix|^2=1,$$ so that every $C_i^*C_i$ for $i \in S$, would fix $x$, which contradicts the assumption $\bigcap_{i \in S}C_i^*H_i=\{0\}.$ It suffices to then consider $s_i=\frac{p_i}{\lambda_S}$ for $i \in S.$ 

\hfill \\
On $\{|h_i|>\eta\}$, $1+h_i \le (1+\eta^{-1})|h_i|$ which means that $$\int (1+h_i)^{2p_i/s}\mathbbm{1}_{\{|h_i|>\eta\}}\, d\gamma_{n_i} \le C_\eta \int h_i^2 \, d\gamma_{n_i} \le C_\eta \epsilon^2,$$ for any $s \in [p_i;1].$ 
\hfill \\
Applying \cref{Holder_AM_GM}, we obtain that $$\int \prod_{i \in S}(1+h_i)^{2p_i}\mathbbm{1}_{\{|h_i|>\eta\}}\, d\gamma_n \le C_\eta \epsilon^{2\sum_{i \in S}s_i}=o_\eta(\epsilon^2).$$ In addition, all other factors satisfy $|h_j(C_jx)| \le \eta$, for any $j \in S(x)^c.$

\hfill \\
Now, assume that $S=S(x)$ and $$\bigcap_{i \in S}C_i^*H_i \ne \{0\},$$ and consider $v$ a unit vector in this set. Then, $$1=\sum_{i=1}^m p_i|C_iv|^2 \ge \sum_{i \in S}p_i,$$ so that by the arithmetic-geometric mean inequality, $$\prod_{i \in S}(1+x_i)^{p_i} \le 1+\sum_{i \in S}p_ix_i.$$ Then, $$R(x) \ge \Big(1+\sum_{i \in S}p_ix_i+\sum_{j \not \in S}p_jx_j \Big)^2-\Big(1+\sum_{i \in S}p_ix_i\Big)^2 \prod_{j \notin S}(1+x_j)^{2p_j}.$$ Moreover, $$\Big| \sum_{j \notin S}p_j h_j(C_jx) \Big| \le \eta \sum_{j \notin S}p_j \text{ and } \prod_{j \notin S}\Big(1+h_j(C_jx)\Big)^{2p_j}=1+2\sum_{j \notin S}p_j h_j(C_jx)+O(\eta^2).$$ Hence, $$\begin{array}{lll}
\ds R\Big[ (h_i(C_ix))_{1 \le i \le m} \Big] &\ge& \ds 2\Big(1+\sum_{i \in S}p_ih_i(C_ix)\Big)\sum_{j \notin S}p_jh_j(C_jx)+\Big(\sum_{j \notin S}p_jh_j(C_jx)\Big)^2+\Big(1+\sum_{i \in S}p_ih_i(C_ix)\Big)^2  \\
&& \ds -\Big(1+\sum_{i \in S}p_ih_i(C_ix)\Big)^2 \Big(1+2\sum_{j \notin S}p_jh_j(C_jx)+C\eta^2\Big)  \\
&\ge& \ds -2\sum_{j \notin S}p_jh_j(C_jx) \Big(1+\sum_{i \in S}p_ih_i(C_ix) \Big) \sum_{i \in S}p_ih_i(C_ix) +\Big(\sum_{j \notin S}p_jh_j(C_jx) \Big)^2 \\ 
&& \ds -C\eta^2 \Big(1+\sum_{i \in S}p_ih_i(C_ix) \Big)^2 \\ 
&\ge& \ds -C\eta \Big( \Big(\sum_{i \in S}p_ih_i(C_ix)^2 \Big)^{1/2}+\sum_{i \in S}p_i h_i(C_ix)^2 \Big)-C\eta^2 \Big(1+\sum_{i \in S}p_i h_i(C_ix)^2 \Big),
\end{array}$$ since $|\sum_{i \in S}p_ih_i(C_ix)| \le (\sum_{i \in S}p_i h_i(C_ix)^2)^{1/2}.$
On $M_a$, $\sum_{i=1}^m p_i h_i(C_ix)^2 \ge \min_{1 \le i \le m}p_i \cdot a^2$, which implies that $$R\Big[ (h_i(C_ix))_{1 \le i \le m} \Big] \ge -C_a\eta \sum_{i=1}^m p_ih_i(C_ix)^2.$$
\hfill \\
We conclude that $$\int_{M_a}R\big[ (h_i(C_ix))_{1 \le i \le m} \big]\, d\gamma_n \ge -C_a \eta \epsilon^2 -C_\eta \epsilon^{2+\alpha},$$ where $\alpha>0$ only depends on the datum $(\textbf{C},p)$. If we define $h_i^{(a)}=h_i \mathbbm{1}_{\{|h_i| \le a\}}$, for any $r \in (0;1)$, there exists $a \in [r^2;r]$ such that $$\int_{M_a}\sum_{i=1}^m p_i h_i^{(a)}(C_ix)^2 \, d\gamma_n(x) \le \frac C{\log(1/r)}\epsilon^2.$$ Indeed, $$\int_0^{+\infty} \mathbbm{1}_{\{\max_j |x_j|>a\}} \sum_{i=1}^m p_i x_i^2 \mathbbm{1}_{\{|x_i| \le a\}} \frac{da}{a}=\sum_{i:x_i \ne 0} p_i x_i^2 \log \Big(\frac{\max_j |x_j|}{|x_i|}\Big),$$ and since $x \mapsto x^2 \log(1/x)$ is bounded on $[0;1]$, the right-hand-side is at most $$C\sum_{i=1}^m p_i x_i^2,$$ also making use of the straightforward inequality $$\min_{j}p_j \cdot \max_{j}|x_j|^2 \le \sum_{i=1}^m p_i x_i^2.$$ In particular, $$\int_0^{+\infty} \bigg[\int_{M_a}  \sum_{i=1}^m h_i^{(a)}(C_ix)^2 \, d\gamma_n(x)\bigg] \frac{da}{a} \le C\epsilon^2,$$ so that there must exist $a \in [r^2;r]$, as above. 
Now, let us show that $$\begin{array}{lll}
\ds \int R\big[(h_i(C_ix))_{1 \le i \le m}\big]\, d\gamma_n &\ge& \ds \tilde Q \big((h_i^{(a)})_{1 \le i \le m} \big)-\int_{M_a}q((h_i^{(a)} \circ C_i)_{1 \le i \le m})\, d\gamma_n \\ 
&& \ds -Ca\epsilon^2+\int_{M_a}R\big[(h_i(C_ix))_{1 \le i \le m}\big]\, d\gamma_n.
\end{array}$$
Indeed, for $x \in M_a^c$, $|h_i(C_ix)| \le a$ for any $1 \le i \le m$, which implies that $$R\big[ (h_i(C_ix))_{1 \le i \le m} \big] \ge q\big( (h_i^{(a)}(C_ix))_{1 \le i \le m}\big)-Ca\sum_{i=1}^m p_ih_i(C_ix)^2,$$ and so, $$\int_{M_a^c} R\big[(h_i(C_ix))_{1 \le i \le m}\big]\, d\gamma_n \ge \int_{M_a^c}q\big[ (h_i^{(a)}(C_ix))_{1 \le i \le m} \big]\, d\gamma_n-Ca\epsilon^2.$$ To deduce the claim, it suffices to remember that $$\tilde Q\big( (h_i^{(a)})_{1 \le i \le m}\big)=\int q\big( (h_i^{(a)}(C_ix))_{1 \le i \le m}\big) \, d\gamma_n.$$ 

\hfill \\
To conclude, we combine what was previously derived. 
\hfill \\
Since trivially, $q \big( (h_i^{(a)}(C_ix))_{1 \le i \le m}\big) \le \sum_{i=1}^m p_i h_i^{(a)}(C_ix)^2$, its integral on $M_a$ is controlled by $\frac C{\log(1/r)}\epsilon^2.$ 
In addition, $$\begin{array}{lll}
\ds \tilde Q((h_i^{(a)})_{1 \le i \le m})=\sum_{i=1}^m p_i \mathrm{Var}_{\gamma_{n_i}}(h_i^{(a)})-\mathrm{Var}_{\gamma_n}\Big(\sum_{i=1}^m p_i h_i^{(a)} \circ C_i\Big) \\ \ds +\sum_{i=1}^m p_i \Big(\int h_i^{(a)}\, d\gamma_{n_i}\Big)^2-\Big(\sum_{i=1}^m p_i \int h_i^{(a)}\, d\gamma_{n_i}\Big)^2 \ge -C\Big(\frac12+\frac 1a\Big)^2 \epsilon^4,   
\end{array}$$ since the first part is positive and using in the second part, that $\sum_{i=1}^m p_i|\int h_i^{(a)}\, d\gamma_{n_i}| \le (\frac12+\frac1a)\epsilon^2$, since $\int |h_i|\mathbbm{1}_{|h_i| \ge a}\, d\gamma_{n_i} \le a^{-1}\int h_i^2 \, d\gamma_{n_i}.$

\hfill \\
Considering $\eta=r^6, a \in [r^2;r]$, we find $$\int R\big[(h_i(C_ix))_{1 \le i \le m} \big]\, d\gamma_n \ge -C\Big(r+\frac1{\log(1/r)}\Big)\epsilon^2-C_r(\epsilon^4+\epsilon^{2+\alpha}),$$ so that for fixed $\theta>0$, we can choose $r>0$ small with $$C\Big(r+\frac1{\log(1/r)}\Big)<\frac \theta 2,$$ and $\epsilon_\theta>0$ such that $C_r(\epsilon_\theta^2+\epsilon_\theta^\alpha)<\frac\theta 2$, $$\int R\big[ (h_i(C_ix))_{1 \le i \le m} \big]\, d\gamma_n \ge -\theta \epsilon^2,$$ if $\epsilon \in (0;\epsilon_\theta).$ 

\hfill \\
Now that the claim is shown, it suffices to use the spectral gap. More precisely, since $(h_i)_{1 \le i \le m}$ is orthogonal to the centered kernel of the quadratic form, $(h_i-\int h_i\, d\gamma_{n_i})_{1 \le i \le m}$ is orthogonal to the whole kernel of the quadratic form, that was denoted by $E$ in \eqref{Def/Thm_Kernel_QF}. Moreover, $$\Big \|\big(h_i-\int h_i \, d\gamma_{n_i}\big)_{1 \le i \le m}\Big\|_{2,p}^2=\epsilon^2-O(\epsilon^4) \Longrightarrow \tilde Q(h) \ge c\epsilon^2-C\epsilon^4,$$ so that considering $\theta=\frac c 4$ and taking possibly $\epsilon_0$ even smaller, we get that $$1-\int_{\mathbb{R}^n} \prod_{i=1}^m f_i(C_ix)^{p_i}\, d\gamma_n(x) \ge \frac c2 \epsilon^2.$$
\end{proof}

Let us now conclude the proof of \cref{Complete_Gaussian_stability}. As observed at the beginning of \cref{Sec:3}, it suffices to prove the main theorem in the case of a geometric Brascamp-Lieb datum. 

\begin{proof}[Proof of \cref{Complete_Gaussian_stability}]
For any $(F_i)_{1 \le i \le m}$ integrable, there exists $(G_i)_{1 \le i \le m}$, with $\int G_i\, d\gamma_{n_i}=1$ for any $1 \le i \le m$ and $$\int_{\mathbb{R}^n}\prod_{i=1}^m G_i(C_ix)^{p_i}\, d\gamma_n(x)=1,$$ such that $$A[(F_i)_{1 \le i \le m}] \ge \frac14\sum_{i=1}^m p_i\Bigg \|\sqrt{\frac{F_i}{\|F_i\|_{L^1(\gamma_{n_i})}}}-\sqrt{G_i} \, \Bigg\|_{L^2(\gamma_{n_i})}^2.$$
\hfill \\
By \cref{Hellinger_Stab}, we deduce that for any $(F_i)_{1 \le i \le m}$ integrable, there exists $(G_i)_{1 \le i \le m}$, with $\int G_i\, d\gamma_{n_i}=1$ for any $1 \le i \le m$ and $$\int_{\mathbb{R}^n}\prod_{i=1}^m G_i(C_ix)^{p_i}\, d\gamma_n(x)=1,$$ such that $$
 1-\frac{\int_{\mathbb{R}^n}\prod_{i=1}^m F_i(C_ix)^{p_i}\, d\gamma_n(x)}{\prod_{i=1}^m \Big(\int_{\mathbb{R}^{n_i}}F_i\, d\gamma_{n_i}\big)^{p_i}} \ge \ds \frac c4\sum_{i=1}^m p_i \Bigg \|\sqrt{\frac{F_i}{\|F_i\|_{L^1(\gamma_{n_i})}}}-\sqrt{G_i} \, \Bigg\|_{L^2(\gamma_{n_i})}^2. $$ 
In particular, note that by considering $F_i=\frac{f_i}{\gamma_{n_i}}$ and $g_i=G_i\gamma_{n_i}$, we get the Lebesgue version of the stability of Brascamp-Lieb inequality: there exists $(g_i)_{1 \le i \le m}$, with $\int_{\mathbb{R}^{n_i}} g_i\, dx_i=1$ and $$\int_{\mathbb{R}^n}\prod_{i=1}^m g_i(C_ix)^{p_i}\, dx=1,$$ such that $$1-\frac{\int_{\mathbb{R}^n} \prod_{i=1}^m f_i(C_ix)^{p_i}\, dx}{\prod_{i=1}^m \Big(\int_{\mathbb{R}^{n_i}}f_i \, dx\Big)^{p_i}} \ge c\sum_{i=1}^m p_i \Bigg\|\sqrt{\frac{f_i}{\|f_i\|_{L^1}}}-\sqrt{g_i}\, \Bigg\|_{L^2}^2.$$ 
\end{proof}

In fact, it is possible to argue alternatively to prove the global quantitative stability in the special case $\frac12<p_i<1$ for any $1 \le i \le m$, and more generally, for any $p_i \in (0;1)$ such that $$H_{i_1} \cap H_{i_2}=\{0\} \quad \forall i_1 \ne i_2,$$ in the setting of \cref{Hellinger_Stab}. We used here that if $\frac12<p_i<1$ for all $1 \le i \le m$, then for $x \in H_{i_1} \cap H_{i_2}$, $$\sum_{j=1}^m p_j C_j^*C_jx=x \Longrightarrow |x|^2=(p_{i_1}+p_{i_2})|x|^2+\sum_{j \notin \{i_1,i_2\}}p_j|C_jx|^2,$$ so that since $p_{i_1}+p_{i_2}>1$, $x=0.$ In particular, if the image spaces are disjoint, we have only Gaussian extremizers. 

\hfill \\
The main reason is because in the general case, one must treat differently the non-linear remainder, as not every interaction involving a ``large'' factor is negligible. This fact is true when $(C_i^*H_i)_{1 \le i \le m}$ do not overlap, whereas it was needed to collect all the interactions into a single remainder and show that its negative contribution is negligible for the general case.

\begin{proof}
    It suffices to consider the local Hellinger stability around the standard Gaussian. 
    
    \hfill \\
    We may also assume that $0<p_i<1$ for all $1 \le i \le m$, since if $p_i=1$, $H_i$ is orthogonal to all other $H_j$ and is not relevant to the inequality. 
    \hfill \\
    Set $$h_i=\sqrt{f_i}-1, \quad \int_{\mathbb{R}^{n_i}} h_i \, d\gamma_{n_i}=-\frac12 \|h_i\|_{L^2(\gamma_{n_i})}^2,$$ and denoting $\epsilon=\sqrt{2A\big[ (f_i)_{1 \le i \le m}\big]} \in [0;\sqrt 2\epsilon_0]$, we have $$\epsilon^2 =\sum_{i=1}^m p_i \|h_i\|_{L^2(\gamma_{n_i})}^2,$$ as the nearest Gaussian extremizer to $(f_i\, d\gamma_{n_i})_{1 \le i \le m}$ is $(d\gamma_{n_i})_{1 \le i \le m}.$ Note that $(h_i)_{1 \le i \le m}$ is orthogonal to the centered tangent space to the extremizer family at the standard Gaussian $$\tilde E=\Big\{ \big (\langle C_ib,x_i\rangle+\langle AC_i^*x_i,C_i^*x_i\rangle-\mathrm{Tr}(C_iAC_i^*) \big)_{1 \le i \le m} : b \in H, A \text{ self-adjoint with } AC_i^*C_i=C_i^*C_iA \Big\},$$ since $K_{\mathrm{dep}}=H$.
    
    \hfill \\
    We will decompose $h_i$ in small and large values. In order to do so, introduce for $0<p<1$, $$k^{(1)}_p(x)=(1+x)^{2p}-1, \quad k_p^{(2)}(x)=p(1+x)^2+1-p-(1+x)^{2p}.$$ Notice that $$k_p^{(2)}(x) \ge \min\{p,1-p\} x^2,$$ since if $p \le \frac12$, $t=1+x$, then by concavity\footnote{We use here that a differentiable function $f:I \rightarrow \mathbb{R}$ is concave on an interval $I$ if $$f(x) \le f(y)+f'(y)(x-y) \quad \forall x,y \in I.$$}, $$t^{2p} \le 1-2p+2pt \Longrightarrow pt^2+1-p-t^{2p} \ge p(t-1)^2,$$ and if $p \ge \frac12$, by concavity of $t^{2p-1}$, $$t^{2p} \le (2p-1)t^2+2(1-p)t \Longrightarrow pt^2+1-p-t^{2p} \ge (1-p)(t-1)^2.$$ In addition, $$k_p^{(1)}(x)=2px+O(x^2), \quad k_p^{(2)}(x)=2p(1-p)x^2+O(x^3),$$ and since $k_p^{(2)}(x)=p(1+x)^2-p-k_p^{(1)}(x)$, $$-\int_{\mathbb{R}^{n_i}} k_{p_i}^{(1)}(h_i)\, d\gamma_{n_i}=\int k_{p_i}^{(2)}(h_i) \, d\gamma_{n_i}.$$  

\hfill \\
There are now two important consequences from disjointness: 
\begin{itemize}
    \item[(i)] For any $S \subset \{1,\cdots,m\}$ with $|S| \ge 2$, there exists $s_i \ge p_i, i \in S$ such that $$\sum_{i \in S}s_iC_i^*C_i \le \mathrm{id}_H \text{ and } \sum_{i \in S}s_i>1,$$ which was shown in Step 3 of the proof of \cref{Hellinger_Stab}. 

\hfill \\ 
    \item[(ii)] For each ordered pair $i \ne j$, there exists $a>\max\{p_i,\frac12\}$ and $b \in (0;1/2)$ such that $a+b>1$ and $aC_i^*C_i+bC_j^*C_j \le \mathrm{id}_H.$ 
    \hfill \\
    Indeed, we choose $a \in (\max\{p_i,\frac12\};1)$ and by the previous argument, $$\|aC_i^*C_i+(1-a)C_j^*C_j\|<1,$$ so that we may take $b>1-a$ such that $b<\frac12$ and $$\|aC_i^*C_i+bC_j^*C_j\| \le 1.$$
\end{itemize}

\hfill \\
Now, define $$A_i=\{ |h_i|>\eta\}, h_i^{(\eta)}=h_i \mathbbm{1}_{A_i^c} \text{ and } r_i=h_i \mathbbm{1}_{A_i} \text{ for any } 1 \le i \le m.$$ Then, $$h_i=h_i^{(\eta)}+r_i \text{ and } \|h_i\|_{L^2(\gamma_{n_i})}^2=\|h_i^{(\eta)}\|_{L^2(\gamma_{n_i})}^2+\|r_i\|_{L^2(\gamma_{n_i})}^2,$$ and define further $$v_i=k_{p_i}^{(1)}(h_i)\mathbbm{1}_{A_i^c}, w_i=k_{p_i}^{(1)}(h_i) \mathbbm{1}_{A_i} \quad \forall 1 \le i \le m.$$ First, $|x| \ge \eta \Longrightarrow |k_{p_i}^{(1)}(x)| \le C_\eta |x|^{2p_i}$, which means that $$\begin{array}{lll}
\ds \gamma_{n_i}(A_i) \le \eta^{-2}\|h_i\|_{L^2(\gamma_{n_i})}^2 \le  C\eta^{-2}\epsilon^2 \Longrightarrow \int |w_i|^{1/s}\, d\gamma_{n_i} \le C_{\eta,s} \int_{A_i}|h_i|^{2p_i/s}\, d\gamma_{n_i} \\ \\ 
\ds \le C_{\eta,s} \|h_i\|_{L^2(\gamma_{n_i})}^{2p_i/s} \gamma_{n_i}(A_i)^{1-\frac{p_i}{s}} \le C_{\eta,s}\epsilon^2,
\end{array}$$ for any $s \in [p_i;1].$ 
Moreover, $$\|v_i\|_{L^\infty} \le C\eta, \|v_i\|_{L^2(\gamma_{n_i})} \le C\epsilon \Longrightarrow \|v_i\|_{L^{1/s}(\gamma_{n_i})} \le C_{\eta,s} \epsilon^{\min\{1,2s\}},$$ for any $0<s \le 1.$
\hfill \\
Now, $$f_i^{p_i}=(1+h_i)^{2p_i}=1+k_{p_i}^{(1)}(h_i),$$ so that $$\int_{\mathbb{R}^n}\prod_{i=1}^m f_i(C_ix)^{p_i}\, d\gamma_n(x)=1+\sum_{i=1}^m \int_{\mathbb{R}^{n_i}} k_{p_i}^{(1)}(h_i)\, d\gamma_{n_i}+\sum_{\substack{ S \subset \{1,\cdots,m\} \\ |S| \ge 2}} \int \prod_{i \in S}(v_i+w_i)(C_ix)\, d\gamma_n(x).$$ 
\hfill \\
For terms containing at least two large factors, all the small terms are bounded, so that we have $|S| \ge 2$ and $$\int \prod_{i \in S}|w_i(C_ix)| \, d\gamma_n \le \prod_{i \in S}\|w_i\|_{L^{1/s_i}(\gamma_{n_i})} \le C_\eta \epsilon^{2\sum_{i \in S}s_i},$$ and recall by what we established above, that $\sum_{i \in S}s_i>1$, so that this integral has smaller order than $\epsilon^2.$
\hfill \\
For terms containing exactly one large factor, all the other small terms are bounded, so that we are left with $i \ne j$ and an integral of the form $$\int_{\mathbb{R}^n} |w_i(C_ix) v_j(C_jx)|\, d\gamma_n(x) \le \|w_i\|_{\frac1a} \|v_j\|_{\frac1b} \le C\epsilon^{2(a+b)},$$ since $a \ge p_i$ and $b<\frac12,$ with $a+b>1.$ 
\hfill \\
Hence, there exists $\alpha>0$ depending only on the datum $(\textbf{C},p)$, such that the absolute value of the sum of all terms containing a large factor is less than $C\epsilon^{2+\alpha}.$ 
\hfill \\
For terms containing the small factors, if at least three occur, then we can bound the contribution by $C\eta \epsilon^2$, and if there are only two, then $$v_i=2p_i h_i^{(\eta)}+e_i, \text{ where } |e_i| \le C\eta|h_i^{(\eta)}| \text{ is the remainder,}$$ which means that $$\int_{\mathbb{R}^n}v_i(C_ix)v_j(C_jx)\, d\gamma_n(x)=4p_ip_j \int_{\mathbb{R}^n}h_i^{(\eta)}(C_ix)h_j^{(\eta)}(C_jx)\, d\gamma_n(x)+O(\eta \epsilon^2).$$ This implies that $$1-\int_{\mathbb{R}^n}\prod_{i=1}^m f_i(C_ix)^{p_i}\, d\gamma_n(x)=\sum_{i=1}^m \int_{\mathbb{R}^{n_i}} k_{p_i}^{(2)}(h_i)\, d\gamma_{n_i}-4\sum_{i<j}p_ip_j \int_{\mathbb{R}^n} h_i^{(\eta)}(C_ix)h_j^{(\eta)}(C_jx)\, d\gamma_n(x)+\mathbf{R},$$ where $$|\mathbf{R}| \le C\eta \epsilon^2+C_\eta \epsilon^{2+\alpha}.$$ 
\hfill \\
Notice that $$\|b\|_{2,p}^2-\|Tb\|_{2}^2=\sum_{i=1}^m p_i(1-p_i)\|h_i^{(\eta)}\|_{L^2(\gamma_{n_i})}^2-2\sum_{i<j}p_ip_j \int_{\mathbb{R}^n}h_i^{(\eta)}(C_ix)h_j^{(\eta)}(C_jx)\, d\gamma_n(x),$$ where we note $b=(h_i^{(\eta)})_{1 \le i \le m}.$ Furthermore, $$\bigg|\int_{\mathbb{R}^{n_i}}h_i^{(\eta)}\, d\gamma_{n_i}\bigg|= \bigg|\int_{\mathbb{R}^{n_i}}h_i\, d\gamma_{n_i}-\int_{\mathbb{R}^{n_i}}r_i\, d\gamma_{n_i}\bigg|\le \Big(\frac12+\frac1\eta\Big) \|h_i\|_{L^2(\gamma_{n_i})}^2 \le C_\eta \epsilon^2,$$ so that replacing $(h_i^{(\eta)})_{1 \le i \le m}$ by $(h_i^{(\eta)}-\int h_i^{(\eta)}\, d\gamma_{n_i})_{1 \le i \le m}$ changes the quadratic expression $\|b\|_{2,p}^2-\|Tb\|_2^2$ by $O_\eta(\epsilon^4)$, as we need to add $\sum_{i=1}^m p_i (\int h_i^{(\eta)}\, d\gamma_{n_i})^2-(\sum_{i=1}^m p_i \int h_i^{(\eta)} \, d\gamma_{n_i})^2$.  

\hfill \\
Now, for $\{u^{(1)},\cdots,u^{(M)}\}$ an orthonormal basis for $\tilde E$, which is composed of centered polynomials of degree at most two,  $(h_i)_{1 \le i \le m} \perp \tilde E$ implies that $$\Big\langle \Big(h_i^{(\eta)}-\int h_i^{(\eta)} \, d\gamma_{n_i} \Big)_{1 \le i \le m},(u^{(j)}_i)_{1 \le i \le m}\Bigl\rangle_p=-\langle (r_i)_{1 \le i \le m}, (u_i^{(j)})_{1 \le i \le m} \rangle_p,$$ and $$\begin{array}{lll}
\ds \big|\langle (r_i)_{1 \le i \le m}, (u_i^{(j)})_{1 \le i \le m} \rangle_p \big | &\le& \ds \sum_{i=1}^m p_i \|r_i\|_{L^2(\gamma_{n_i})}\|u_i^{(j)}\mathbbm{1}_{A_i}\|_{L^2(\gamma_{n_i})} \\ \\
 &\le& \ds \sum_{i=1}^m p_i \|r_i\|_{L^2(\gamma_{n_i})}\|u_i^{(j)}\|_{L^4(\gamma_{n_i})} \gamma_{n_i}(A_i)^{1/4} \le C\eta^{-1/2}\epsilon^{3/2}.
\end{array}$$ Therefore, $$\Big \|\mathrm{Proj}_{\tilde E}\Big((h_i^{(\eta)}-\int h_i^{(\eta)}\, d\gamma_{n_i})_{1 \le i \le m} \Big) \Big\|_{2,p}^2 \le C\eta^{-1}\epsilon^3.$$ 
\hfill \\
Finally employing the spectral gap with $(h_i^{(\eta)}-\int h_i^{(\eta)}\, d\gamma_{n_i})_{1 \le i \le m}$, we get $$\|b\|_{2,p}^2-\|Tb\|_2^2 \ge c \|b\|_{2,p}^2-C_\eta(\epsilon^3+\epsilon^4).$$ Writing $p_*=\min_i \frac{\min\{p_i,1-p_i\}}{p_i}>0$, the local expansion of $k_{p_i}^{(2)}$ yields $$\sum_{i=1}^m \int k_{p_i}^{(2)}(h_i)\, d\gamma_{n_i} \ge 2\sum_{i=1}^m p_i(1-p_i) \|h_i^{(\eta)}\|_{L^2(\gamma_{n_i})}^2+p_*\|r\|_{2,p}^2-C \eta \epsilon^2.$$
Putting altogether, we have $$1-\int_{\mathbb{R}^n}\prod_{i=1}^m f_i(C_ix)^{p_i}\, d\gamma_n(x) \ge 2c\|b\|_{2,p}^2+p_*\|r\|_{2,p}^2-C\eta \epsilon^2-C_\eta(\epsilon^{2+\alpha}+\epsilon^3+\epsilon^4).$$ Recall that $\|b\|_{2,p}^2+\|r\|_{2,p}^2=\epsilon^2$ and $c_*=\min\{2c,p_*\}>0.$ Finally, choose $\eta>0$ such that $C\eta \le c_*/4$ and then, $\epsilon_0>0$ such that $$C_\eta(\epsilon_0^\alpha+\epsilon_0+\epsilon_0^2) \le c_*/4.$$ Thanks to these choices, we deduce that $$\mathrm{Def}(f_1,\cdots,f_m) \ge 1-e^{-\mathrm{Def}(f_1,\cdots,f_m)} \ge \frac{c_*}4 \epsilon^2.$$  
\end{proof}

To conclude this final section, we prove that the exponent of the Hellinger distance in \cref{Hellinger_Stab} is sharp. Let us first show the following consequence of the proof of \cref{Hellinger_Stab}. 

\begin{proposition} \label{Sharpness_Extremizable}
Assume that a curve of normalized non-negative tuples $(f_i^{(\epsilon)})_{1 \le i \le m}$ satisfies $$\sqrt{f_i^{(\epsilon)}}=1+\epsilon h_i+o_{L^2}(\epsilon), \quad \tilde Q((h_i)_{1 \le i \le m})>0,$$ and $$1-\int_{\mathbb{R}^n}\prod_{i=1}^m f_i^{(\epsilon)}(C_ix)^{p_i}\, d\gamma_n(x)=O(\epsilon^2).$$ Then, $$A\big[(f_i^{(\epsilon)})_{1 \le i \le m}\big] \asymp \epsilon^2.$$ Consequently, for every $\beta>1$, $$\frac{A\big[(f_i^{(\epsilon)})_{1 \le i \le m}\big]}{\ds \Big(1-\int_{\mathbb{R}^n}\prod_{i=1}^m f_i^{(\epsilon)}(C_ix)^{p_i}\, d\gamma_n(x)\Big)^\beta} \rightarrow_{\epsilon \rightarrow 0}+\infty.$$ 
\end{proposition}

\begin{proof}
    Recall that $$\tilde Q((\tilde h_i)_{1 \le i \le m})=\tilde Q(\tilde h):=\sum_{i=1}^m p_i\|\tilde h_i\|_{L^2(\gamma_{n_i})}^2-\Big\|\sum_{i=1}^m p_i \tilde h_i \circ C_i\Big\|_{L^2(\gamma_n)}^2,$$ where $\tilde h=(\tilde h_i)_{1 \le i \le m}.$ For any normalized non-negative tuple $(g_i)_{1 \le i \le m}$, setting $\tilde h_i=\sqrt{g_i}-1$ for any $1 \le i \le m$, Step 3.1 and the remainder estimate of Step 3.2 give \begin{equation} \label{Remainder_Reloaded} 1-\int_{\R^n}\prod_{i=1}^m g_i(C_ix)^{p_i}\, d\gamma_n(x) \ge \tilde Q((\tilde h_i)_{1 \le i \le m})-\omega(\|\tilde h\|_{2,p})\|\tilde h\|_{2,p}^2, \quad \omega(\|\tilde h\|_{2,p}) \rightarrow 0 \text{ as } \|\tilde h\|_{2,p} \rightarrow 0.\end{equation} Note that this does not require the fact that $h$ is orthogonal to \cref{Def/Thm_Kernel_QF}, which is afterwards used to apply the spectral gap.
    \hfill \\
    First, $$A\big[(f_i^{(\epsilon)})_{1 \le i \le m}\big] \le \frac12 \sum_{i=1}^m p_i\|\sqrt{f_i^{(\epsilon)}}-1\|_{L^2}^2=O(\epsilon^2).$$ Suppose by contradiction, that the lower bound fails. Along a sequence $\epsilon_k \rightarrow_{k \rightarrow +\infty} 0$, we would have $$A\big[(f_i^{(\epsilon_k)})_{1 \le i \le m}\big]=o(\epsilon_k^2).$$ By definition of $A$ as an infimum, we may choose normalized extremizers $(g_i^{(\epsilon_k)})_{1 \le i \le m}$ such that $$\Big\|\big(\sqrt{g_i^{(\epsilon_k)}}-\sqrt{f_i^{(\epsilon_k)}}\big)_{1 \le i \le m}\Big\|_{2,p}=o(\epsilon_k).$$ Hence, $$h_i^{(\epsilon_k)}=\sqrt{g_i^{(\epsilon_k)}}-1=\epsilon_k h_i+o_{L^2}(\epsilon_k),$$ so that by continuity of $\tilde Q$ and \cref{Remainder_Reloaded}, $$\begin{array}{lll}
    \ds 0 = 1-\int_{\mathbb{R}^n}\prod_{i=1}^m g_i^{(\epsilon_k)}(C_ix)^{p_i}\, d\gamma_n(x) &\ge& \ds \tilde Q\big( (h_i^{(\epsilon_k)})_{1 \le i \le m}\big)-\omega(\|h^{(\epsilon_k)}\|_{2,p})\|h^{(\epsilon_k)}\|_{2,p}^2 \\
    \nm 
    &\ge& \ds \epsilon_k^2 \tilde Q((h_i)_{1 \le i \le m})+o(\epsilon_k^2)>0,
    \end{array}$$ for sufficiently large $k$, which is a contradiction. Therefore, $A[(f_i)_{1 \le i \le m}] \ge c\epsilon^2$ as desired. 
\end{proof}

\hfill \\
Now, we may for example slightly perturb the standard Gaussian in a ``preferred'' direction $f_i^{(\epsilon)}=1+\epsilon v_i+O(\epsilon^2)$, which we will take to be in $I_k^{(n_i)}(\mathrm{Sym}^k(H_i))$, and more precisely an eigenvector that does not belong to the kernel of $\mathrm{Id}-T_k^*T_k$. Here, we assume that some $p_i$ are strictly less than $1$, since if $p_i=1$ for all $1 \le i \le m$, the Brascamp-Lieb inequality is an identity and both the deficit and $A$ are zero.

\hfill \\
Hence, for a smooth curve of normalized positive densities, $$f_{i,k}^{(\epsilon)}=\frac{1+\epsilon v_{i,k}(x)+\epsilon^2 v_{i,k}(x)^2}{1+\epsilon^2 \|v_{i,k}\|_{L^2(\gamma_{n_i})}^2}=1+\epsilon v_{i,k}+\epsilon^2(v_{i,k}^2-\|v_{i,k}\|_{L^2(\gamma_{n_i})}^2)+o_{L^2}(\epsilon^2) \text{ as } \epsilon \rightarrow 0,$$ with $v_{i,k}=I_k^{(n_i)}(u_{i,k}) \in I_k^{(n_i)}(\mathrm{Sym}^k(H_i))$ and $w_{i,k}=v_{i,k}^2-\|v_{i,k}\|_{L^2(\gamma_{n_i})}^2.$  Let us, for now, drop the dependence on the index $k \ge 1$ considered to be fixed. 
\hfill \\
We now want to expand $$\int \prod_{i=1}^m f_i^{(\epsilon)}(C_ix)^{p_i}\, d\gamma_n(x).$$ Recall that $$(1+\epsilon v+\epsilon^2 w)^p=1+p \epsilon v+\epsilon^2 \Big(\frac{p(p-1)}{2} v^2+pw \Big)+o(\epsilon^2),$$ so that $$\begin{array}{lll}
\ds \prod_{i=1}^m f_i^{(\epsilon)}(C_ix)^{p_i}=1+\epsilon \sum_{i=1}^m p_i v_i(C_ix)+\epsilon^2 \sum_{i=1}^m p_i w_i(C_ix) \\ \\
\ds +\frac{\epsilon^2}2 \Big[\Big(\sum_{i=1}^m p_i v_i(C_ix)\Big)^2-\sum_{i=1}^m p_i v_i(C_ix)^2\Big]+o(\epsilon^2).
\end{array}$$
Moreover, $$\int_{\mathbb{R}^n} v_i(C_ix)\, d\gamma_n=\int_{\mathbb{R}^{n_i}}v_i \, d\gamma_{n_i}=0,$$ by the previously established properties of Hermite polynomials and similarly, $$\int_{\mathbb{R}^n}w_i(C_ix)\, d\gamma_n(x)=0.$$ We conclude that $$\int \prod_{i=1}^m f_i^{(\epsilon)}(C_ix)^{p_i}\, d\gamma_n(x) =1+\frac{\epsilon^2}{2}k! \Big(\|T_ku\|^2-\|u\|_p^2\Big)+o(\epsilon^2),$$ so that $$1-\int \prod_{i=1}^m f_i^{(\epsilon)}(C_ix)^{p_i}\, d\gamma_n(x)=\frac{\epsilon^2}{2}k!\Big(\|u\|_p^2-\|T_ku\|^2\Big)+o(\epsilon^2).$$  
Suppose now that $T_k^*T_ku=\lambda u$ with $\lambda \in [0;1).$ In this case, $$\|T_ku\|^2=\langle T_k^*T_ku,u\rangle_p=\lambda \sum_{i=1}^m p_i \|u_i\|^2.$$ 
\hfill \\
Thus, $$\begin{array}{lll}
\ds 1-\int \prod_{i=1}^m f_i^{(\epsilon)}(C_ix)^{p_i}\, d\gamma_n(x) &=& \ds \frac{\epsilon^2}2 (1-\lambda)k!\sum_{i=1}^m p_i \|u_i\|^2+o(\epsilon^2) \\ \\
&=& \ds \frac{\epsilon^2}2 (1-\lambda)\sum_{i=1}^m p_i \|v_i\|_{L^2(\gamma_{n_i})}^2+o(\epsilon^2)=O(\epsilon^2).
\end{array}$$ 
\hfill \\
Since $$\tilde Q((v_i)_{1 \le i \le m})=\frac{\epsilon^2}{2}(1-\lambda)\sum_{i=1}^m p_i\|v_i\|_{L^2(\gamma_{n_i})}^2+o(\epsilon^2)>0,$$ we may apply \cref{Sharpness_Extremizable} which proves the sharpness.

\hfill \\
\textbf{AI Statement:}
The mathematical work in this paper was carried out by the author without the use of AI. After completing this work, a commercial AI tool was used to check the arguments and correct minor typos.

\section*{Acknowledgements}

The author is grateful to his supervisor Alessio Figalli for his support and useful discussions throughout the elaboration of this work.

\enlargethispage{\baselineskip}
\printbibliography

\end{document}